\documentclass[11pt]{article}

\usepackage{amsmath}
\usepackage{amssymb,amsthm}
\usepackage{boxedminipage}
\usepackage{makeidx}
\usepackage{epsfig,wrapfig}
\usepackage{psfrag}
\usepackage{xspace}
\usepackage{color}

\usepackage[bookmarks]{hyperref}

\numberwithin{equation}{section}

\input myhyphen.sty

\def\Rd{\color[rgb]{0.0,0,0}}

\makeatletter
\renewcommand{\@seccntformat}[1]{\@nameuse{the#1}.\quad}
\makeatother

\gdef\hyp#1{^{{\Rd (#1)}}}

\def\girth{{\rm girth}}

\def\dense{\itemsep=2pt\parskip=0pt    
 }

\def\dori{\hfill\break\phantom{~}\hskip\parindent }
\def\dorib{\hfill\break\hbox{\quad} $\bullet$\hbox{~~} }

\def\emR{\em\Rd  }

\def\aa{{\Rd  {\mathbf a}}}
\def\bb{{\Rd {\mathbf b}}}

\def\xx{{\Rd  \mathbf x}}
\def\yy{{\Rd \mathbf y}}

\def\jj{{\bf j}}

\def\C{{\cal C}}

\def\M{{\Rd \mathcal M}}

\def\TT{{\Rd\cal T}}
\def\QQ{{\Rd \mathbb Q}}

\def\UU{{\Rd\mathbb U}}
\def\LL{{\Rd  \mathcal L}}

\def\ZZ{{\Rd\mathbb Z}}

\def\Ck{{\Rd C_k}}
\def\Gm{{\Rd  G_m}}
\def\Gn{{\Rd  G_n}}
\def\Hn{{\Rd  H_n}}
\def\Sn{{\Rd S_n}}
\def\Pk{{\Rd  P_k}}
\def\Tnp{{\Rd T_{n,p}}}

\def\Tk{\,{\Rd  T_k}}

\def\embed{{\hookrightarrow\, }}

\def\TembG{\Tk\embed\Gn}

\def\ti{\tilde}

\def\Hnu{{\Rd  H_\nu}}

\def\mindeg{d_{\sf min}}
\def\maxdeg{d_{\sf max}}

\def\proof. {\medbreak\noindent{\bf \Rd  Proof. \enspace}}
\def\Proof#1. {\medbreak\noindent{\bf\Rd  Proof #1.\enspace}}

\def\Qed{\hfill$\Box$}

\def \Tur(#1){\left(1-{1\over #1}\right)}

\def\rema#1. {\medskip\noindent{\bf\Rd #1.\enspace}}

\gdef\ext#1){{\Rd  {\bf ex}#1)}}

\gdef\extcon#1){{\Rd  {\bf ex_{\rm con}}#1)}}

\gdef\extconn#1){{\Rd  {\bf ex_{\rm 2-con}}#1)}}

\gdef\EXT#1){{\Rd  {\bf EX}#1)}}

\def\beq#1{\begin{equation}\label{#1}}
\def\eeq{\end{equation}}

\def\beqn{\begin{eqnarray}}
\def\eeqn{\end{eqnarray}}

\def\beqns{\begin{eqnarray*}}
\def\eeqns{\end{eqnarray*}}

\def\text#1{~\mbox{#1}~}
\def\Text#1{\qquad\mbox{#1}\qquad}

\def\half{\frac{1}{2}}
\def\reci#1 {{1\over #1}}

\def\HH{{\Rd\mathbb H}}
\def\HHH#1#2{{\Rd\mathcal H_{#2}^{(#1)}}}

\def\eps{{\Rd \varepsilon}}

\newtheorem{theorem}{\Rd Theorem}[section]

\newtheorem{claim}[theorem]{\Rd Claim}

\newtheorem{lemma}[theorem]{\Rd Lemma}

\newtheorem{conjecture}[theorem]{\Rd Conjecture}
\newtheorem{construction}[theorem]{Construction}
\newtheorem{corollary}[theorem]{Corollary}

\theoremstyle{definition}

\newtheorem{definition}[theorem]{\Rd\bf  Definition}

\newtheorem{problem}[theorem]{\Rd Problem}

\newtheorem{remark}[theorem]{\Rd  Remark}
\newtheorem{remarks}[theorem]{\Rd  Remarks}

\def\SepaRef{ 
\noindent{\bf Further sources to read:} }

\gdef\BalAbra#1{\begin{wrapfigure}{L}{30mm}\epsfig{file=#1,width=30mm}
\end{wrapfigure}}

\gdef\balAbra#1{\begin{wrapfigure}{L}{23mm}\epsfig{file=#1,width=24mm}
\end{wrapfigure}}

\gdef\BBalAbra#1#2{\begin{wrapfigure}{L}{#2mm}\epsfig{file=#1.eps,width=#2mm}
\end{wrapfigure}}

\newcount\picsize
\newtheorem{thm}{Theorem}

\begin{document}

\def\subsub#1//{\medskip\noindent {\bf\Rd  #1~}}

\gdef\Otimes{{\otimes}}

\gdef\Rendez#1{
\marginpar{\epsfig{file=konyv.eps,width=8mm}~\fbox{#1}}}

\def\HK33{K_{3,3}}

\def\extbip{{\Rd\bf ex_{\rm bip}}}

\def\exto{{\Rd \bf ext^*_{\nearrow}}}
\def\exto{{\Rd \bf ext_{mat}}}
\def\bM{{\bf M}}
\def\bL{{\bf L}}
\def\bP{{\bf P}}
\def\G{{\Rd\mathcal G}}

\def\B{{\Rd\cal B}}
\def\D{{\Rd\cal D}}

\def\Cg#1{{\Rd\mathcal C_{\ge #1}}}

\def\MR{~\hfill\fbox{\bf\Rd MR}}

\def \bibin#1{\medskip ---\medskip \hskip-10mm \fbox{#1}~}

\def \bibin#1#2\par{\par}

\newcount\bibnum

\def\bibi#1{\global\advance\bibnum by1
\bibitem{#1}{\Rd\bf [N]\the\bibnum.~ }}

\makeatletter
\renewcommand{\@seccntformat}[1]{\@nameuse{the#1}.\quad}
\makeatother

\gdef\girth{{\bf girth}}
\gdef\EE{{\mathbb E}}

\def\Terv#1#2{\marginpar{\fbox{\bf #1}$\longrightarrow$\\ #2}}

\def\Lt#1{L_{#1}(L,\psi)}

\gdef\diam{{\bf diam}}

\def\Boxit#1{\medskip\vbox{\hrule\hbox{\vrule\kern3pt\vbox{\kern3pt
\vbox{\advance\hsize -6.8pt{
#1}\par}\kern3pt}\kern3pt\vrule}\hrule}\medskip}

\newcount\abranum
\newcount\abraref
\newcount\mpgw
\newcount\mpgh

\gdef\figref#1{\abraref\abranum
\global\advance\abraref by #1 Fig.~\the\abraref}

\def\Unfinished{
\global\advance \sracnum by 1
{\Rd \fbox{\tiny \the\sracnum}
\epsfig{file=Rdanger.eps,height=4mm} }
 \marginpar{\epsfig{file=Rdanger.eps,height=9mm}
{\Rd \margofont
\fbox{\the\sracnum}\\ Unfinished}}}

\def\doublebox#1#2#3#4{
\hbox{
\begin{minipage}{#3cm}
#1
\end{minipage}
\qquad
\begin{minipage}{#4cm}
\epsfig{file=#2.eps,width=#4cm}
\\ {\sf An almost extremal graph?}
\end{minipage}}}

\def\YFig#1/#2/#3/#4/     
{\medskip
\leftskip=-4mm\par
\noindent\vbox{\epsfig{file=#1.eps,width=34mm}
\par\leftskip=3mm
\parindent=0pt\par
\medskip{\quad \sf #4}}
~\smallskip~
\leftskip=0pt
\vskip -#2cm
\hangindent42mm\hangafter-#3
}

\def\ZFig#1/#2/#3
{\medskip
\marginpar{\fbox{#1}}
\begin{minipage}{4cm}
\begin{center}
\epsfig{file=#1.eps,width=4cm,height=3cm}
\par
{\small\it #2}
\end{center}
\end{minipage}
\leftskip=0pt\par
\vskip -37mm
\hangindent52mm\hangafter-#3
}

\def\WFig#1/#2/#3/#4/#5\par
{\medskip
\marginpar{\fbox{\small fi=#1: ha=#3 up=#4}}
\vbox{\epsfig{file=#1.eps,width=5cm,height=4cm}
\par\leftskip=3mm
\parindent=0pt\par
\medskip{\small\it #2}}
\leftskip=0pt\par
\vskip -#4mm
{\hangindent62mm\hangafter-#3
#5\par}
\par
}

\def\WTHFig#1/#2/
{\medskip
\marginpar{\fbox{#1}}
\vbox{\epsfig{file=#1.eps,width=5cm,height=4cm}
\par\leftskip=3mm
\parindent=0pt\par
\medskip{\small\it #2}}
\leftskip=0pt\par
\vskip -8cm
\hangindent62mm\hangafter-12
}

\def\LTH#1{
\vskip -#1cm
\hangindent65mm\hangafter-13}

\def\pic#1 {\fbox{\bf Figure #1}~\marginpar{\quad\fbox{Fig}}}

\def\subsubsub#1{\subsubsection*{#1}}



\def\Erdos1{{\Prp Erd\H{o}s}}
\def\Furedi1{{\Prp F\"uredi}}
\def\Sos1{{\Prp S\'os}}


\def\turtwo#1{\left[#1^2\over 4\right]}

\def\th{^{\rm th}}

\def\apr{$\sim1_k$-}

\def\ga{{\Rd \gamma}}

\def\be{\beta }
\def\bex{{1\over 5}}
\def\betaTP{{1\over 100}}
\def\gaX{{1\over 1000} }
\def\gax{{1\over 10000} }

\def\dep{\partial}
\def\deS{\Delta_G }
\def\ka{\kappa }

\def\gar{\ga^*} 
\def\de{{\Rd \delta}}
\def\om{\omega}
\def\OS{{\Rd \Omega^*}}
\def\OSS{{\Rd \Omega^{**}}}

\def\odd1{{\Prp\sc Odd}}
\def\even1{{\Prp\sc Even}}

\def\densecase{{\Rd ``Dense Case''}}

\def\sparsecase{{\Rd ``Sparse Case''}}

\def\K#1{{\Rd K_{#1}}}
\def\T#1{{\Rd T^{(#1)}}}
\def\GG{{\Rd\cal G}}
\def\UiVi{{\Rd (U_i,V_i)}}
\def\GUiVi#1{{\Rd G^{#1}[U_i,V_i]}}

\def\KK#1#2{{\Rd K_{#1,#2}}}

\def\sqe{\sqrt{\eta}}
\def\sqek{\sqrt\eta k}
\def\sqep{\sqrt\varepsilon k}
\def\del{\sqrt{\delta}}

\def\sectionx#1{\section*{#1}\marginpar{x}}
\def\subsectionx#1{\subsection*{#1}\marginpar{x}}

\def\Abra#1\par{\vskip1cm\hrule\vskip1cm{Figure Missing:
#1}\vskip1cm\hrule\vskip1cm}

\def\abox#1{\hbox{~\vbox{\hrule\hbox{\vrule\kern3pt\vbox
{\kern3pt\vbox{\advance\hsize
 -6.8pt \hbox{{~#1~}}}\kern3pt}\kern3pt\vrule}\hrule}~}}

\def\chk{$^*$\marginpar{** ?? **}}

\def\Case#1:
{\bigskip{\bf Case #1:\enspace}}

\def\kerdes{\marginpar{$\leftarrow$\hskip-1.6mm\fbox{\bf ?}}}

\newtheorem{lem}{\Rd Lemma}  
\newtheorem{histrem}{\Rd Historical Remark}  

\def\negyes#1#2#3#4{${#1#2\choose #3#4}$}

\def\BiB{\fbox{Big Bang}\enspace}

\def\ext#1){{\Rd {\bf ex}#1)}}

\def\0{\emptyset}

\def\A{{\cal A}}
\def\bA{{\mathbf A}}
\def\cA{{\cal A}}
\def\bB{{\mathbf B}}
\def\cB{{\cal B}}

\def\AA{{\Rd \mathbb A}}
\def\BB{{\Rd \mathbb B}}
\def\CC{{\Rd \mathbb C}}
\def\DD{{\Rd \mathbb D}}
\def\F{{\cal F}}
\def\SS{{\Rd \mathbb S}}
\def\SSs{{\Rd \mathbb S^*}}
\def\VS{{\Rd (\mathbb V-\mathbb S)}}
\def\VsS{{\Rd (\mathbb V^*\hskip-1.4mm-\mathbb S)}}
\def\Vs{{\Rd \mathbb V-\mathbb S}}
\def\VSs{{\Rd \mathbb X}}
\def\FTutte{{\mathcal F_{\bf Tutte}}}
\def\FVS{{\Rd \mathcal F_{V-S}}}
\def\NN{{\mathbb N}}
\def\MM{{\mathcal M}}   
\def\O{{\Rd \Omega}}
\def\U{{\bf U}}
\def\WW{{\mathbb W}}

\def\oldal#1{[$\rightarrow$]\marginpar{SZ/M: #1}}

\def\CoupUnc1{{\Prp \bf Coupled-Uncoupled}}
\def\matching1{{\Prp \sf matching}}

\def\skeleton{{\bf SKELETON}}

\def\xilent#1\par{#1\par}

\def\suru{\itemsep=0pt\parskip=0pt}
\def\notation. {\medbreak\noindent{\bf Notation. \enspace}}

\def\Hova{\marginpar{Hova?}}

\def\Expand{\marginpar{\fbox{\bf Expand?}}}

\newcount\remcount
\newcount\Sracnum

\def\atem{\global \advance \remcount by 1
\vskip1pt\hrule\vskip2pt
\noindent(\the\remcount)~~ }

\title{\bf The history of degenerate
 (bipartite) extremal graph problems}

\author{Zolt\'an F\"uredi \, and \, Mikl\'os Simonovits}
\date{May 15, 2013\\
Alfr\'ed R\'enyi Institute of Mathematics, Budapest,
Hungary
\\
{\tt z-furedi@illinois.edu}
\quad and \quad {\tt
simonovits.miklos@renyi.mta.hu
}}

\pagestyle{myheadings}\markboth{F\"uredi-Simonovits: Degenerate [A]/\jobname}
{{} 
${}$ 
 F\"uredi-Simonovits:
 Degenerate (bipartite) extremal graph problems
 }

\maketitle

   \begin{abstract}
This paper is a survey on Extremal Graph Theory, primarily focusing
on the case when one of the excluded graphs is bipartite.
On one hand we give an introduction to this field and also
  describe many important results, methods, problems, and constructions. \footnote{
Research supported in part by the Hungarian National Science Foundation
OTKA 104343,
and by the European Research Council Advanced Investigators Grant 267195
(ZF)
and by the Hungarian National Science Foundation OTKA 101536,
and by the European Research Council Advanced Investigators Grant 321104.
(MS).
}

   \end{abstract}

\tableofcontents

\gdef\NZ{\marginpar{\fbox{\margofont NZ}}}
\gdef\New{\Prp\sf }
\gdef\PP{{\Rd\mathcal P}}

\newpage


\section{Introduction}

This survey describes the theory of {\emR \index{degenerate}
Degenerate Extremal Graph Problems},
the main results of the field, and its connection to the surrounding areas.

Extremal graph problems we consider here are often called Tur\'an type
extremal problems, because Tur\'an's theorem and his questions were the most
important roots of this huge area  \cite{TuranML}, \cite{TuranColloq}.

Generally, we have a {\emR \index{Universe} Universe} of graphs,
$\UU$, where this universe may be the family of ordinary graphs, or
digraphs, or hypergraphs, or ordered graphs, or bipartite
graphs, etc and a property $\PP$, saying, e.g., that $G\in\UU$ does
not contain some subgraphs $L\in\LL$, or that it is Hamiltonian, or it
is at most 3-chromatic, and we have some parameters on $\UU$, say
$v(G)$ and $e(G)$, the number of vertices and edges. Our aim is to
maximize the second parameter under the condition that $G$ has
property $\PP$ and its first parameter is given.

\begin{quote}
We call such a problem {\emR Tur\'an type extremal problem}
 if we are given a family $\LL$ of
graphs from our universe, $\Gn$ is a graph of $n$ vertices, $e(\Gn)$
denotes the number of edges of $\Gn$ and we try to maximize $e(\Gn)$
under the condition that $\Gn$ contains no $L\in\LL$, where
``contains'' means ``not necessarily induced subgraph''. (Here graph
may equally mean digraph, or multigraph, or hypergraph).

The maximum will be denoted by $\ext(n,\LL)$ and the graphs attaining this
maximum without containing subgraphs from $\LL$ are called {\emR extremal
  graphs}. The family of \index{extremal graph} extremal graphs is denoted by
$\EXT(n,\LL)$ and $\ext(n, \LL)$ is called the {\emR Tur\'an number} of the
family $\LL$.

Speaking of $\ext(n,L)$ we shall always assume that
 $n\ge |V(L)|$, otherwise the problem is trivial.
\end{quote}

\begin{definition}

If the Universe $\UU$ is the family of $r$-uniform
hypergraphs\footnote{$r=2$ included, moreover, mostly we
  think of $r=2$.}, then we
shall call the problem {\emR degenerate} if the maximum,
$$\ext(n,\LL)=o(n^r).$$
Otherwise we shall call it {\emR non-degenerate}
\end{definition}

Below we shall mention several open problems.
Yet to get more problems, we refer the reader to the

\begin{center}Erd\H{o}s
homepage: ~~www.renyi.hu/\~{}p\_erdos\end{center}
where the papers of Erd\H{o}s can be
found  \cite{Erdhomepage}.
Also, many open problems can be found in Chung-Graham
 \cite{ChungGrahamErdLegacy}.


\subsection{Some central theorems of the field}\label{CentralExamplesS}

We start with some typical theorems of the field and two conjectures.  The
aim of this ``fast introduction'' is to give a feeling for what are the
crucial types of results here.

\begin{theorem}[K\H ov\'ari--T. S\'os--Tur\'an,  \cite{KovSosTur}]
\label{KovSosTurTh}
 Let $K_{a,b}$ denote the complete bipartite graph with $a$ and
$b$ vertices in its color-classes.
Then 
$$ \ext(n,K_{a,b}) \le
\half\root a\of{b-1}\cdot n^{2-(1/a)} +O(n).
$$ 
\end{theorem}

We use this theorem with $a\le b$, since that way we get a better estimate.

\begin{theorem}[Koll\'ar-Alon--R\'onyai--Szab\'o
 \cite{KollRonySzab}, \cite{AlonRonySzab}]\label{AlonRonySzab}

If $b>(a-1)!$, then $$\ext(n,\KK ab )>c_an^{2-(1/a)}.$$
\end{theorem}

\begin{theorem}[Erd\H{o}s, Bondy and Simonovits  \cite{BondySim}]

$$\ext(n,C_{2k})\le 100kn^{1+(1/k)}.$$
\end{theorem}

\begin{theorem}[Erd\H{o}s--Simonovits, Cube Theorem  \cite{ErdSimCube}]\label{CubeTh}

Let $Q_8$ denote the cube graph defined by the vertices and edges of a
3-dimensional cube. Then $$\ext(n,Q_8)=O(n^{8/5}).$$
\index{Cube Theorem}
\end{theorem}

\begin{conjecture}[Erd\H{o}s and Simonovits, Rational exponents]\label{RacConj}

For any finite family $\LL$ of graphs, if there is a bipartite $L\in\LL$, then there exists a rational $\alpha\in[0,1)$ and a $c>0$ such that
$${\ext(n,\LL)\over n^{1+\alpha}}\to c.$$
\end{conjecture}

\begin{theorem}[F\"uredi  \cite{FureC4JCT96},  \cite{FureC4Prep}]\label{FureC483}

If $q\ne 1,7,9,11,13$, and $n=q^2+q+1$, then
$$\ext(n,C_4)\le\half q(q+1)^2.$$
Moreover, if $q$ is a power of a prime,
then
$$\ext(n,C_4)= \half q(q+1)^2.$$
\end{theorem}

\begin{conjecture}[Erd\H{o}s]\label{ErdosC3C4}
\footnote{This conjecture is mentioned in  \cite{Erd1975-42} but it is definitely older, see e.g. Brown,  \cite{BrownNonExist} .}

$$\ext(n,\{C_3,C_4 \})={1\over 2\sqrt 2}n^{3/2}+o(n^{3/2}).$$
\end{conjecture}

We close this part with a famous result of \index{Ruzsa-Szemeredi} Ruzsa and Szemer\'edi:

\begin{theorem}[Solution of the (6,3) problem,  \cite{RuzsaSzemer}]\label{RuzsaSzemTh}
 If $\HHH3n$ is a 3-uniform hypergraph not containing 6
vertices determining (at least) 3 hyperedges, then this
hypergraph has $o(n^2)$ hyperedges.
\end{theorem}

The above theorems will be discussed in more details below.

\subsection{The structure of this paper}  

\begin{figure}[ht]
\label{Area}
\epsfig{file=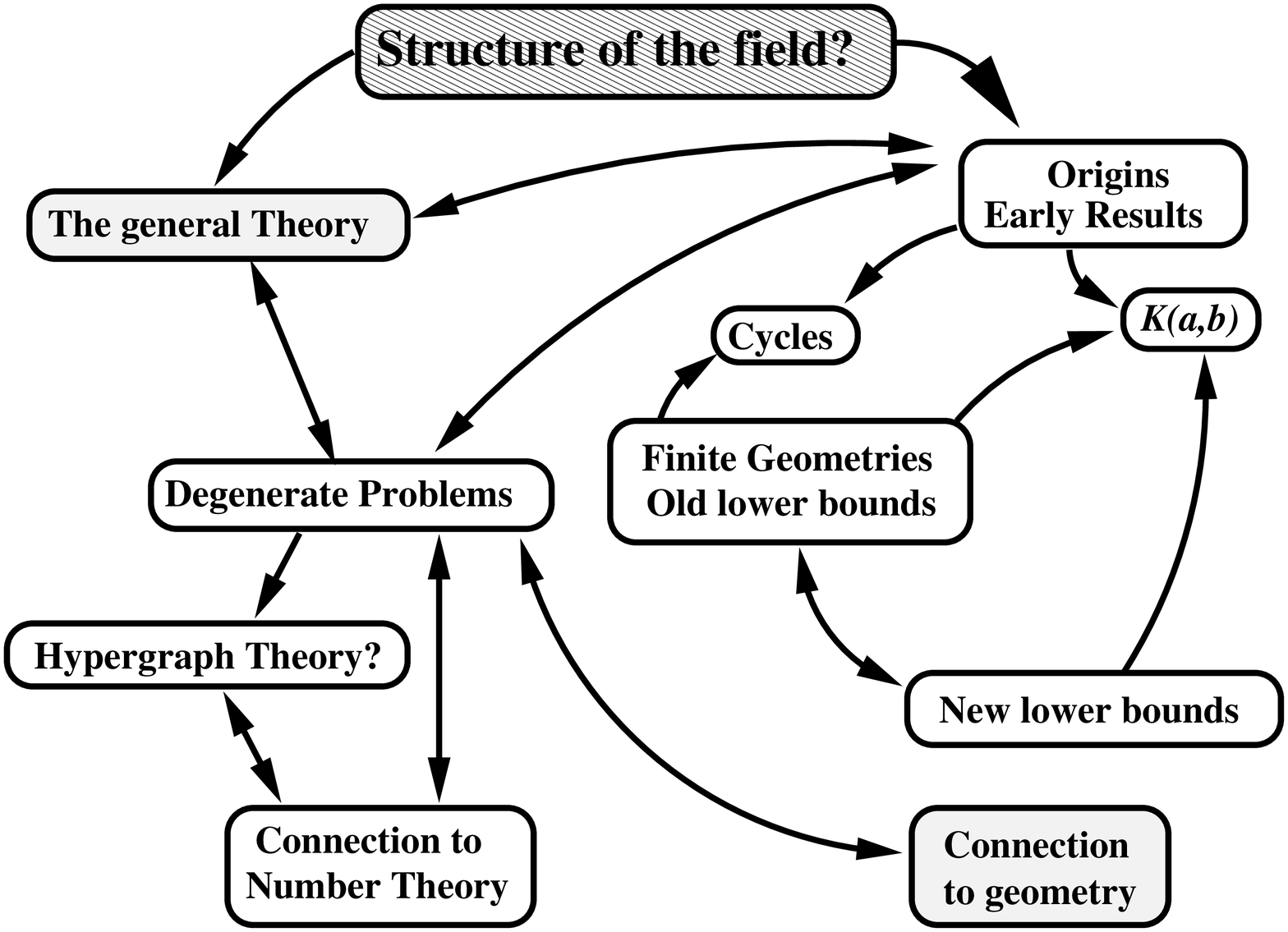,width=10cm,height=6cm}
\caption{Area Map}
\end{figure}

The area is fairly involved.
Figure 1 shows a complicated -- but not complete -- map of the
interactions of some subfields of the field discussed here.
We start with describing the Extremal problems in general, then move
to the Degenerate problems, also describing why they are
important.
Among the most important degenerate extremal problems we mention are
the extremal problem of $K_{a,b}$, and also $C_{2k}$, where --
to classify the extremal problems -- we shall need the Random Graph
Method to get a lower bound in the problem of $C_{2k}$. These
results are enough to give a good \index{classification} {\emR
 classification} of degenerate extremal graph problems.

\begin{enumerate}\dense
\item The two lowest boxes of Figure 1 show that this whole area has a
 strong connection to Geometry and Number Theory. This will be
 explained in Sections \ref{AppliGeom}, \ref{HistorA}, and \ref{AppliNumb}.
\item The origins of this field are
\begin{enumerate}\dense
\item an early, singular result of Mantel,
\item the multiplicative Sidon Problem (see Section \ref{HistorA})
\item Tur\'an's theorem and his systematic approach to the field.
\end{enumerate}
So the origins come -- in some sense -- from Number Theory, are
strongly connected to Finite Geometry, and in this way also to ordinary
geometry (Tur\'an's theorem comes from the Erd\H{o}s-Szekeres version
of Ramsey Theorem, which they invented to solve the Esther Klein
problem from Geometry.)
\item To understand the
  field we start with a very short description of the general theory, and then
  -- skipping most of the hypergraph theory -- we move to the main area of
  this paper: to the questions where we consider ordinary extremal graphs, and
  exclude some bipartite $L$: therefore, by Theorem \ref{KovSosTurTh}, we have
  $\ext(n,\LL)=O(n^{2-c})$.
\item One important phenomenon is that many extremal graph problems can be
  ``reduced'' to some {\emR degenerate} extremal graph problems that we also
  call sometimes {\emR bipartite extremal problems}.
\item The upper bounds in the simpler cases are obtained by some
 double counting, Jensen type inequalities, or applying some
\index{supersaturated graph}
supersaturated graph theorems.\footnote{Lov\'asz and Szegedy had
 a beautiful conjecture, which we formulate here only in a
 restricted form: Any (valid) extremal theorem can be proven by
 applying the Cauchy--Schwarz inequality finitely many times. This
 conjecture was killed in this strong form -- by Hatami and Norine
  \cite{HatamiNorine} -- but proven in a weaker,
 ``approximation-form''. }
\item There are also much more complicated cases, where the above
 simple approach is not enough, we need some finer arguments.
 Perhaps the first such case was treated by F\"uredi:
Section \ref{M11Excluded} and
  \cite{FurediL11}. Also such an approach is the application of the
 general Dependent Random Choice Method, (see
 the survey of Fox and Sudakov
 \cite{FoxSudDRCSurv}).
\item The lower bounds are sometimes provided by random graphs (see
 Section~\ref{ProbLowerS}) but these are often too weak. So we often
 use some finite geometric constructions, (see Section
 \ref{FinGeomA}) or their generalizations -- coming from
 commutative algebras (see Sections \ref{NormGraphS},
 \ref{LuboPhilSarnakS}, \ref{EigenValuesS}), etc., and they
 occasionally provide matching upper and lower bounds. Again, there
 is an important general method with many important results, which we
 shall call the Lazebnik-Ustimenko method but will treat only very
 superficially in Section \ref{NormGraphS}.
\end{enumerate}

\subsection{Extremal problems}\label{s:Expose}

We shall almost entirely restrict ourselves to Tur\'an
type extremal problems for ordinary simple graphs, i.e. loops,
multiple edges are excluded.

To show the relation of the areas described here, we start with a list
of some related areas.

\begin{enumerate}\dense
\item Ramsey Theory
\begin{itemize}\dense
\item Problems not connected to density problems; in some sense these
 are the real Ramsey Problems
\item Problems connected to density problems, i.e. cases, where we do not really use the Ramsey Condition, only that some color class is large.
\end{itemize}
\item Ordinary extremal graph theory
\begin{itemize}\dense
\item \label{Degex} Excluding bipartite graphs (degenerate problems)
\item Excluding topological subgraphs (very degenerate extremal problems)
\item Matrix problems, ordered and not ordered;
\item \label{NonDeg} Non-degenerate case, and its relation to
 degenerate problems
\end{itemize}
\item Theory of extremal digraph problems
\item Ramsey-Tur\'an Problems
\item Connection to Random Graphs
\item Hypergraph extremal problems
\item Connection of Number Theoretical problems to Extremal Graph Theory
\item Continuous problems
\item Applications
\end{enumerate}

There are several surveys on these fields, see e.g., T. S\'os
 \cite{SosLincei}, F\"uredi  \cite{FureLondon},  \cite{FureZurich},
Simonovits  \cite{SimFra},  \cite{SimPragNew},  \cite{SimErdosInflu},
 \cite{SimTurInflu}, Simonovits and S\'os  \cite{SimSosRT},
 \cite{KeevashHyperSurv}.
Perhaps the first survey on this topic was Vera S\'os' paper  \cite{SosLincei},
discussing connections between extremal graph problems, finite geometries,
block designs, etc.
 and, perhaps, the nearest to this survey is
 \cite{SimWat}, Bollob\'as, \cite{BolloHB} Sidorenko  \cite{SidorSurv},
 \cite{FoxSudDRCSurv}, and also some books, e.g., Bollob\'as
 \cite{BolloExtreBook}.
Of course, a lot of information is hidden
in the papers of Erd\H{o}s, among others, in  \cite{Erd1975-42},  \cite{ErdMost},
 \cite{ErdFavourTheor}.

So, here we shall concentrate on Case \ref{Degex}, but to position this
area we shall start with some related fields, among others, with the
general asymptotic in Case~\ref{NonDeg}.

\begin{problem}[General Host-graphs]

In a more general setting we have a sequence $(\Hn)$ of ``host'' graphs
and the question is, how many edges can a subgraph
$\Gn\subset\Hn$ have under the condition that it does not contain any
forbidden subgraph $L\in\LL$.
The maximum will be denoted by $\ext(\Hn,\LL)$.
\end{problem}

For $\Hn=K_n$ we get back the ordinary extremal graph problems.
There are several further important subcases of this question:

(a) when $\Hn=\KK ab $ for $a\approx n/2$;

(b) when the host-graph is the $d$-dimensional cube, $n=2^d$; see Section~\ref{S:Hypercube}.

(c) when $\Hn$ is a \index{random host graph} random graph on $n$ vertices,
see e.g. \cite{RodlSchacht}.

\rema Notation. Given some graphs $\Gn$, $\Tnp$, $\Tk$, $\Hnu$,... the
(first) subscript will almost always denote the number of
vertices.\footnote{One important exception is the complete bipartite
 graph $K(a,b)=\KK ab$, see below.
Another exception is, when we list
 some excluded subgraphs, like $L_1,\dots,L_\nu$.} So $\K p$ is the
complete graph on $p$ vertices, $\Pk$ the path on $k$ vertices, $\Ck$
is the cycle of $k$ vertices, while $\Cg k$ will denote the family of
cycles of length at least $k$. $\de(x)$ denotes the degree of the
vertex $x$.

The complete bipartite graph $K_{a,b}$ with $a$ vertices in its first
class and $b$ in its 2nd class will be crucial in this paper. Often,
we shall denote it by $K(a,b)$, and its $p$-partite generalization by
$K_p(a_1,\dots,a_p)$. If $\sum a_i=n$ and $|n_i-n_j|\le1$, then
$K_p(a_1,\dots,a_p)$ is the Tur\'an graph $\Tnp$ on $n$ vertices and $p$
classes.

Given two graphs $U$ and $W$, their product graph is the one obtained
from vertex-disjoint copies of these two graphs by joining each vertex
of $U$ to each vertex of $W$.
This will be denoted by $U\Otimes W$.\footnote{This product is often called also the joint of the two
 graphs.}

Given a graph $H$, $v(H)$ is its number of vertices, $e(H)$ its number of
edges and $\chi(H)$ its chromatic number, $\mindeg(G)$ and $\maxdeg(G)$ denote
the minimum and maximum degrees of $G$, respectively.

We shall write $f(x)\approx g(x)$ if $f(x)/g(x)\to 1.$
Occasionally $[n]$ denote the set of first $n$ integers, $[n]:= \{ 1, 2, \dots, n\}$.

\rema The Overlap. Some twenty years ago Simonovits wrote a survey
 \cite{SimErdosInflu} on the influence of Paul Erd\H{o}s in the areas
described above, Many-many features of these areas changed drastically
since that. Jarik Ne\v set\v ril and Ron Graham were the editors of
that survey-volume, and now they decided to republish it.
Fortunately, the authors had the option to slightly rewrite their
original papers.
Simonovits has rewritten his original paper  \cite{SimPragNew},
basically keeping everything he could but {\emR indicating} many new
developments, and adding remarks and many new references to it.

To make {\emR this} paper readable and self-contained, we shall touch on some
basic areas also treated there, or in other survey papers of ours,
Here, however, we shall explain many-many results and phenomena just
mentioned in other survey papers.

\begin{remark}
 There is also a third new survey to be mentioned here:
Simonovits gave a lecture at the conference on Tur\'an's $100\th$
anniversary, in 2011.
His lecture tried to cover the whole
influence of Paul Tur\'an in Discrete Mathematics. In the volume of
this conference Simonovits wrote a survey  \cite{SimTur2013} covering
his lecture, except that \dorib the area called Statistical Group
Theory is discussed in a survey of P\'alfy and Szalay
 \cite{PalfySzalay} and \dorib some parts of the applications of
Extremal Graph Theory, primarily in Probability Theory are
covered by Katona  \cite{KatonaTurSurv}, in the same volume.
\end{remark}

\subsection{Other types of extremal graph problems}\label{OtherTypes}

Above we still tried to maximize the number of edges, hyperedges, etc.
More generally, instead of maximizing
$e(\Gn)$, we may maximize something else:

\begin{enumerate}\dense
\item {\emR Min-degree problems} (or Dirac type problems):
How large min-degree can $\Gn$ have without containing subgraphs from $\LL$.
\item {\emR Median problems} which will be called here Loebl-Koml\'os-S\'os
 type problems: Given a graph $\Gn$, which $m$ and $d$ ensure that if
 $\Gn$ has at least $m$ vertices of degree $\ge d$, then $\Gn$
 contains some $L\in\LL$.
\item {\emR Eigenvalue-extremal problems}
\footnote{As usual, given a graph $\Gn$, an
 $n\times n$ matrix is associated to it, having $n$
 eigenvalues. The largest and the second largest is what we are
 mostly interested in.}: maximize the maximum
 eigenvalue $\lambda(\Gn)$
 under the condition that $\Gn$ does not
 contain any $L\in\LL$. (These are sharper forms of some extremal
 results, since the maximum eigenvalue
$$\lambda(\Gn)\ge{2e(\Gn)\over n},$$
see Section \ref{EigenValuesS}.)
\item \label{SubGrCount}
 {\emR Subgraph count inequalities,} which assert that if $\Gn$
contains many copies of some subgraphs $L_1,\dots,L_\lambda$, then we
 have at least one (or maybe ``many'') subgraphs $L$.
\item {\emR Diameter-extremal problems}. Here we mention just a
 subcase: if
$$\diam(\Gn)\le d\Text{ and }\maxdeg(\Gn)<M,$$ at least how
 many edges must $\Gn$ have. The Erd\H{o}s-R\'enyi paper
  \cite{ErdRenyiDiam} is of importance here
 and also some
 related papers, like F\"uredi  \cite{FureDiam3},  \cite{FureDiamMin2}.
\item {\emR Combined extremal problems:} There are many--many further types of extremal problems.
Here we mention, as an example, the results of Balister, Bollob\'as, Riordan
and Schelp  \cite{BaliBollRioSch}, where an odd cycle is excluded, and at the
same time an upper bound is fixed on the degrees
and the number of edges are to be maximized.
\end{enumerate}

 The approach \ref{SubGrCount} is very popular in the theory of
Graph Limits  \cite{LovaszLimitBook}.
We mention a breakthrough in this area in connection with
 Erd\H{o}s' combinatorial problems, of type \ref{SubGrCount}.
A famous conjecture of Erd\H{o}s was

\begin{conjecture}[Erd\H{o}s  \cite{ErdosCubeConj}]

A $K_3$-free $\Gn$ contains at most $({n\over5})^5$ copies of $C_5$'s.
\end{conjecture}

The motivation of this is that the blownup\footnote{Given
 a graph $H$, its blownup version $H[t]$ is defined as follows: we
 replace each vertex $x$ of $H$ by $t$ independent new vertices and
 we join two new vertices coming from distinct vertices $x,y$ iff
 $xy$ was an edge of $H$.} $C_5$, i.e. $C_5[n/5]$ has no triangles
and has $({n\over 5})^5$ copies of $C_5$. Erd\H{o}s conjectured that
no triangle-free $\Gn$ can have more $C_5$'s than this.
The first ``approximation'' was due to Ervin Gy\H{o}ri:

\begin{theorem}[Gy\H{o}ri  \cite{GyoriC3C5}]

A $K_3$-free $\Gn$ contains at most $1.03({n\over5})^5$ $C_5$'s.
\end{theorem}

Next F\"uredi improved the constant to 1.001 (unpublished)
 \cite{FurediC5}, and finally independently Grzesik  \cite{GrzesikC5},
and Hatami, Hladk\'y, Kr\'al, Norine,
and Razborov  \cite{HataRazboC5C3} proved the conjecture.

\subsection{Historical remarks}\label{HistorA}

Erd\H{o}s in 1938  \cite{ErdTomsk} considered the following
``multiplicative Sidon Problem''\footnote{For a longer
 description of the number theoretical parts see
  \cite{SimPragNew}. Erd\H{o}s also refers to his ``blindness''
 overlooking the general problem in  \cite{Erd1975-42}.}.

\begin{problem} \label{MultipSidon}
 How many integers, $a_1,\dots, a_m\in[1,n]$ can we find so
that $a_ia_j=a_ka_\ell$ does not hold for any $i,j,k,\ell$, unless
$\{i,j\} =\{k,\ell\}$.

\end{problem}

To get an upper bound in Problem \ref{MultipSidon},
Erd\H{o}s proved

\begin{theorem}\label{TomskC4}

Let $G[n,n]$ be a bipartite graph with $n$ vertices in both classes.
If it does not contain $C_4$, then $e(G[n,n])<3n\sqrt{n}$,
\end{theorem}

Much later this problem was asked in a more general setting: find an
upper bound on $e(G[n,n])$ if $\KK ab\not\subset G[n,n] $.
Zarankiewicz  \cite{Zarank} posed the following question:

\begin{problem}[Zarankiewicz problem]\label{ZaraProb}

 Determine the largest integer
$Z(m,n,a,b)$ for which there is an $m\times n$ 0-1 matrix containing $Z(m,n,a,b)$ $1$'s
 without an $a\times b$ submatrix consisting entirely of $1$'s.
\end{problem}

Hartman, Mycielski and Ryll-Nardzevski  \cite{HartmanRyll} gave upper
and lower bounds for the case $a=b=2$, weaker than the
Erd\H{o}s-Klein\footnote{In  \cite{Erd1975-42} Erd\H{o}s
 (again) attributes the finite geometric construction to Eszter (Esther)
 Klein.} result, and K\H{o}v\'ari, T. S\'os and Tur\'an (see Theorem
\ref{KovSosTurTh}) provided a more general upper bound. We shall
discuss these problems and results in details in Sections
\ref{ComplBipUp} and \ref{FinGeomA}.

While exact values of $Z(m,n,a,b)$ are known for infinitely many
parameter values, mostly only asymptotic bounds are known in the
general case. Even $Z(m,n,2,2)$ is not sufficiently well known.

\section{The general theory, classification}

In many ordinary extremal problems the minimum chromatic number plays
a decisive role. The {\bf subchromatic number} $p(\LL)$ of $\LL$ is
defined by \beq{SubChrom}p(\LL)=\min\{\chi(L): L\in\LL\}-1.\eeq
Recall that the Tur\'an graph $\Tnp$ is the largest graph on $n$ vertices
 and $p$ classes.


\begin{claim}\label{GenLowerC}
\beq{GenLower}\ext(n,\LL)\ge e(\Tnp)=\left(1-{1\over
p}\right){n\choose 2}+o(n^2). \eeq
\end{claim}

Indeed, $\Tnp$ does not contain any $L\in\LL$.
An easy consequence of the Erd\H{o}s-Stone theorem  \cite{ErdStone}
provides the asymptotic value of $\ext(n,\LL)$, at least if
$p(\LL)>1$.

\begin{theorem}[Erd\H{o}s-Simonovits  \cite{ErdSimLim}]\label{ErdSimLimTh}
 If $\LL$ is a family of
graphs with subchromatic number $p>0$, then $$\ext(n,\LL)=\left(1-{1\over
p}\right){n\choose 2}+o(n^2). $$
\end{theorem}

\BBalAbra{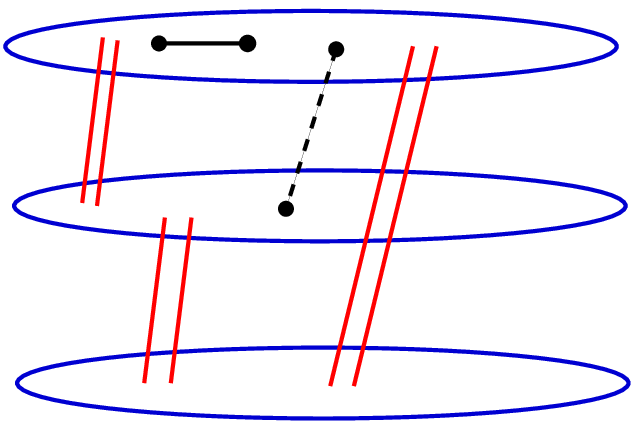}{33}
This means that $\ext(n,\LL)$ depends only very loosely on $\LL$; up
to an error term of order $o(n^2)$;
it is already determined by $p(\LL)$.
\footnote{Better error terms are proved e.g. in
  \cite{SimTihany}, however, this will not be discussed here.} The
question is whether the structure of the extremal graphs is also
almost determined by $p(\LL)$, and (therefore) it must be very similar
to that of $\Tnp$
\footnote{Actually, this was the original question;
   Theorem \ref{ErdSimLimTh} was a partial answer to it.}.
The answer is {\bf YES}.
This is expressed by the following results of Erd\H{o}s and
Simonovits  \cite{ErdRome},  \cite{ErdTihany},  \cite{SimTihany}:

\begin{theorem}[The Asymptotic Structure Theorem]
\label{AsymptStrTh}

Let $\LL$ be a family of forbidden graphs with subchromatic number
$p$. If $\Sn\in\EXT(n,\LL)$, (i.e, $\Sn$ is extremal for $\LL)$, then
it can be obtained from $\Tnp$ by deleting and adding $o(n^2)$
edges. Furthermore, if $\LL$ is finite, then the minimum degree
$$\mindeg(\Sn)=\Tur(p)n+o(n).$$ \end{theorem}

Further, the almost-extremal graphs are similar to $\Tnp$.

\begin{theorem}[The First Stability Theorem]
 Let $\LL$ be a family of forbidden graphs with subchromatic
number $p$. For every $\eps>0$, there exist a $\de>0$ and
an $n_{\eps}$ such that, if $\Gn$ contains no $L\in\LL$, and if,
for $n>n_{\eps}$, \beq{StabA}e(\Gn)>\ext(n,\LL)-\de n^2,\eeq
then $\Gn$ can be obtained from $\Tnp$ by
changing\footnote{deleting and adding} at most $\eps n^2$
edges. \end{theorem}

\begin{remark}

For ordinary graphs ($r=2$) we often call the
degenerate extremal graph problems {\emR bipartite extremal problems.}
This is the case when $\LL$ contains some bipartite graphs.
There is a slight problem here: we shall also consider the
case when not only some $L\in\LL$ is bipartite but $\chi(\Gn)=2$
is as well.

\end{remark}

\subsection{The importance of the Degenerate Case}

There are several results showing that if we know sufficiently well
the extremal graphs for the degenerate cases, then we can also reduce the
non-degenerate cases to these problems.

\subsubsection*{Exact Tur\'an numbers, product conjecture}

 We start with an illustration. Let $O_6=K(2,2,2)$ be the octahedron graph.
Erd\H{o}s and Simonovits proved that

\begin{theorem}[Octahedron Theorem  \cite{ErdSimOcta}]

 If $\Sn$ is an extremal graph for the
octahedron $O_6$ for $n$ sufficiently large, then there exist extremal
graphs $G_1$ and $G_2$ for the circuit $C_4$ and the path $P_3$ such that
$\Sn=G_1\Otimes G_2$ and $|V(G_i)|=\half n+o(n)$, $i=1,2$.

If $G_1$ does not contain $C_4$ and $G_2$ does not contain $P_3$, then
$G_1\Otimes G_2$ does not contain $O_6$. Thus, if we replace $G_1$ by
any $H_1$ in $\EXT(v(G_1),C_4)$ and $G_2$ by any $H_2$ in
$\EXT(v(G_2),P_3)$, then $H_1\Otimes H_2$ is also extremal for $O_6$.
\end{theorem}

More generally, 

\begin{theorem}[Erd\H{o}s--Simonovits  \cite{ErdSimOcta}]\label{ErdSimOctaTh}

 Let $L$ be a complete $(p+1)$-partite graph, $L:=K(a,b,r_3,r_4,...,r_{p+1})$,
where $r_{p+1}\ge r_p\ge\dots\ge r_3\ge b\ge a$ and $a=2,3$. There exists an $n_0=n_0(a,b,\dots,r_{p+1})$ such that
if $n>n_0$ and $S_n\in\EXT(n,L)$, then $S_n=U_1\Otimes U_2\Otimes \dots\Otimes U_p$, where
\begin{enumerate} \dense
\item $v(U_i)=n/p+o(n)$, for $i=1,\dots,p$.
\item $U_1$ is extremal for $\KK ab$
\item $U_2,U_3,\dots,U_p\in \EXT(n,K(1,r_3))$.
\end{enumerate}
 \end{theorem}

It follows that this theorem is indeed a reduction theorem.

\begin{conjecture}[The Product Conjecture, Simonovits]\label{ProdConj}

Assume that $p(\LL)=\min_{L\in\LL}\chi(L)-1>1$.
If for some constants $c>0$ and $\eps\in(0,1)$
\beq{SupLinLow}\ext(n,\LL)>e(\Tnp)+cn^{1+\eps},\eeq
then
there exist $p$ forbidden families $\MM_i$, with
$$
p(\MM_i)=1\Text{and}
\max_{M\in\M_i} v(M)\le \max_{L\in\LL} v(L),$$
such that
for any $\Sn\in\EXT(n,\LL)$, $\Sn=G_1\Otimes\dots\Otimes G_p$, where
$G_i$ are extremal for $\MM_i$.
\end{conjecture}

This means that the extremal graphs $\Sn$ are
products of extremal graphs for some degenerate extremal problems (for
$\MM_i$), and therefore we may reduce the general case
to degenerate extremal problems.

\begin{remarks}

(a) If we allow infinite families $\LL$, then one can easily find
 counterexamples to this conjecture.

(b) If we allow linear error-terms, i.e. do not assume
 \eqref{SupLinLow}, then one can also find counterexamples, using a
 general theorem of Simonovits  \cite{SimSymm77}; however, this is not
 trivial at all, see  \cite{SimBirk}.

(c) A weakening of the above conjecture would be the following: for
 arbitrary large $n$, in Conjecture \ref{ProdConj} there are several
 extremal graphs, and for each $n>n_\LL$, some of them are of product
 form, (but maybe not all of them) and the families $\MM_i$ also may depend
 on $n$ a little.

\end{remarks}

\SepaRef Griggs, Simonovits, and Thomas  \cite{GriggsSim}, Simonovits,  \cite{SimStirin}.

\subsection{The asymmetric case of Excluded Bipartite graphs}\label{AsymQueS}

The degenerate extremal graph problems have three different forms:

\begin{problem}[Three versions]\label{OurProblems}

(a) Ordinary extremal graph problems, where some bipartite or
 non-bipartite sample graphs are excluded, and we try to maximize
 $e(\Gn)$ under this conditions.

(b) The {\emR bipartite case}, where the host graph is
$K(m,n)$ and we maximize $e(G_{n+m})$ under the
 conditions that $G_{n+m}\subseteq K(m,n)$ and $G_{n+m}$ contains no
 $L\in\LL$. (Here we may assume that all
 $L\in\LL$ are bipartite.) In this case we often use the notation
 $\ext(m,n,\LL)$. If $m\le n$ but $m>c n$ for some constant $c>0$,
 then the answer to this problems and to the problem of $\ext(n,\LL)$
 are the same, up to a constant. If, however, we assume that $n$ is
 much larger than $m$, then some surprising new phenomena occur, see
 Section~\ref{AppliNumb}.

(c) The {\emR asymmetric case.}
Color the vertices of the sample graphs $L$ in RED-BLUE
 and exclude only those $\Gn\subseteq K(m,n)$ where the RED
 vertices of some $L\in\LL$ are in the FIRST class of $K(m,n)$:
 maximize $e(G_{n+m})$ over the remaining graphs $G_{n+m}\subseteq
 K(m,n)$.

\end{problem}

Denote the maximum number of edges in this third case by
$\ext^*(m,n,\LL)$.

\begin{remark}
 We have seen Zarankiewicz' problem (i.e. Problem
\ref{ZaraProb}). That corresponds to an asymmetric graph problem, (c). If
we exclude in an $m\times n$ matrix both a $a\times b$ and an
$b\times a$ submatrices of 1's, that will correspond to a bipartite graph
problem, (b).
\end{remark}

\begin{conjecture} [Erd\H{o}s, Simonovits  \cite{SimWat}]
$$\ext^*(n,n,\LL)=O(\ext(n,\LL)).$$
\end{conjecture}

The simplest 
 case when we cannot prove this is $L=K(4,5)$.

\begin{remark}[Matrix problems]
 Case (c) has also a popular matrix form where we
consider 0-1 matrices and consider an $m\times n$ matrix not
containing a submatrix $A$. The question is: how many $1$-s can be in
such a matrix.
This problem has (at least) two forms: the unordered and ordered one.
We return to the Ordered Case in Subsection \ref{MatrixTardos}.
\end{remark}

\subsection{Reductions: Host graphs}

The following simple but important observation shows that there is not
much difference between considering any graph as a ``host'' graph or
only bipartite graphs.

\begin{lemma}[Erd\H{o}s' bipartite subgraph lemma]\label{ErdLemHalf}

Every graph $\Gn$ contains a bipartite subgraph $\Hn$ with $e(\Hn)\ge \half e(\Gn).$
\end{lemma}

This lemma shows that there is not much difference between
considering $K_{2n}$ 
 or $K_{n,n}$
 as a host graph. 

\begin{corollary}

If $\ext_B(n,\LL)$ denotes
the maximum number of edges in an $\LL$-free bipartite graph,
then $$\ext_B(n,\LL)\le\ext(n,\LL)\le 2\,\ext_B(n,\LL).$$
\end{corollary}

Assume now that we wish to have an upper bound on $\ext(m,n,\LL)$,
where $n\gg m$. One way to get such an upper bound is to partition the
$n$ vertices into subsets of size $\approx m$. If, e.g, we know that
$\ext(m,m,\LL)\le cm^{1+\ga}$, then
we obtain that
\beq{SubPsrtEstim} \ext(m,n,\LL)\le {n\over m}\cdot \ext(m,m,\LL)\le c
 nm^\gamma .\eeq
This often helps, however, occasionally it is too weak.
Erd\H{o}s formulated

\begin{conjecture}
If $n>m^2$ then $\ext(m,n,C_6)=O(n).$
\end{conjecture}

Later this conjecture was made more precise, by Erd\H{o}s, A. S\'ark\"ozy and T. S\'os,
and proved by Gy\H{o}ri, see Section \ref{AppliNumb} and  \cite{GyoriC6}.

We start with a trivial lemma.

\begin{lemma}

Let $d$ be the average degree in $\Gn$, i.e. $d:=2e(\Gn)/n$.
Then $\Gn$ contains a $\Gm$ with $\mindeg(\Gm) \ge d/2$.
\end{lemma}

To solve the cube-problem, Erd\H{o}s and Simonovits used two
reductions.\index{Cube Theorem} The first one was a reduction to
bipartite graphs, see Section \ref{ErdLemHalf}. The other one
eliminates the degrees are much higher than the average.

\begin{definition}[$\Delta$-almost-regularity]

$G$ is $\Delta$-almost-regular if $\maxdeg(G)<\Delta\cdot \mindeg(G)$.
\end{definition}

\begin{theorem}[$\Delta$-almost-regularization  \cite{ErdSimCube}]
\label{DeltaReg}

Let $e(\Gn)>n^{1+\alpha}$, and $\Delta=20\cdot 2^{(1/\alpha)^2}$. Then there is a $\Delta$-almost-regular $\Gm\subset\Gn$
for which
$$e(\Gm)>{2\over 5}m^{1+\alpha},
\Text{where} m>n^{\alpha{1-\alpha\over 1+\alpha}} ,$$
unless $n$ is too small.
\end{theorem}

This means that whenever we wish to prove that
$\ext(n,\LL)=O(n^{1+\alpha})$, we may restrict ourselves to
bipartite $\Delta$-almost-regular graphs.

It would be interesting to understand the limitations of this lemma better.
The next remark and problem are in this direction.

\begin{remark}

By the method of random graphs one can show  \cite{ErdSimCube} that for
every $\Delta$ and $n$ and $\eps>0$, there is $\Gn$ with $e(\Gn)=\lfloor
n^{3/2}\rfloor$ which does not have a $\Delta$-almost-regular subgraph $\Gm$
with $e(\Gm)> \eps\sqrt n m$.
\end{remark}

\begin{problem}[Erd\H{o}s-Simnovits \cite{ErdSimCube}]
Is it true that for every $\Delta$ there exists an $\eps>0$
 such that every $\Gn$, with $e(\Gn)=\lfloor n \log n\rfloor$, contains a
 $\Delta$-almost-regular subgraph $\Gm$, with $e(\Gm) > \eps m \log m$ where
 $m\to\infty$ when $n\to\infty$?
\end{problem}

\subsection{Excluding complete bipartite graphs} \label{ComplBipUp}

Certain questions from topology (actually, Kuratowski theorem on planar graphs)
led to Zarankiewicz problem  \cite{Zarank}. After some weaker results
K\H ov\'ari, T. S\'os and Tur\'an proved the following theorem, already
mentioned in Section \ref{CentralExamplesS}.

\begin{theorem}[K\H ov\'ari--T. S\'os--Tur\'an,  \cite{KovSosTur}]
\label{KovSosTurThB}
 Let $K_{a,b}$ denote the complete bipartite graph with $a$ and
$b$ vertices in its color-classes. Then
\beq{KSTshf}
\ext(n,K_{a,b}) \le \half\root a\of{b-1}\cdot n^{2-(1/a)} +{a-1\over 2}n.\eeq
\end{theorem}

\begin{remarks}
 (a) If $a\neq b$ then \eqref{KSTshf} is better if we apply it
with $a<b$.

(b) We know from Theorem \ref{ErdSimLimTh} that $\ext(n,\LL)=o(n^2)$
if and only if $\LL$ contains a bipartite $L$. Actually Claim \ref{GenLowerC}
and Theorem \ref{KovSosTurThB} show that if
$\ext(n,\LL)=o(n^2)$ then
$\ext(n,\LL)=O(n^{2-c})$, for some constant $c=c_{\LL}>0$.


\end{remarks}

\begin{conjecture}[\cite{KovSosTur}, see also e.g.  \cite{Erd1975-42}]

The upper bound in Theorem \ref{KovSosTurTh} is sharp: $$\ext(n,\KK ab )>c_{a,b}n^{2-(1/a)}.$$
\end{conjecture}

{\emR Sketch of proof of Theorem \ref{KovSosTurThB}. }
The number of $a$-stars $K_{a,1}$ in a graph $\Gn$ is $\sum {d_i\choose a}$
 where $d_1,\dots, d_n$ are the degrees in $\Gn$.
If $\Gn$ contains no $K_{a,b}$ then at most $b-1$ of these $a$-stars
 can share the same set of endpoints.
We obtain
\begin{equation}\label{eq:KST}
\sum {d_i\choose a} ={\rm the\, number\, of\, }a{\rm -stars} \le  (b-1){n\choose a}.
 \end{equation}
Extending $n\choose a$ to all $x>0$ by
 $${x\choose a}:=\begin{cases}{x (x-1)\dots(x-a+1)\over a!}& \text{for}
 x\ge a-1,\\ 0&\text{otherwise}\\
\end{cases}
$$ we have a convex function.
Then Jensen's Inequality implies that, the
 left hand side in \eqref{eq:KST} is at least
 $n{2e(G)/n\choose a},$ and the result follows by an easy calculation.
\Qed

\begin{remark}
 Slightly changing the above proof we get analogous upper
bounds on $e(\Gn)$ in all three cases of Problem \ref{OurProblems}.
\end{remark}

We shall return to these questions in Sections \ref{C4ZaraS},
\ref{Zara2S}
where
we shall discuss some improvements of the upper bound and also some
lower bounds.

\SepaRef Guy  \cite{GuyTihany},
Zn\'am:  \cite{ZnamA},  \cite{ZnamB}, Guy--Zn\'am  \cite{GuyZnam}.

\subsection{Probabilistic lower bound}\label{ProbLowerS}

The theory of random graphs is an interesting, important, and rapidly
developing subject. The reader wishing to learn more about it should either
read the original papers of Erd\H{o}s, e.g.,  \cite{ErdGraphProbA},
 \cite{ErdGraphProbB}, Erd\H{o}s and R\'enyi, e.g.,  \cite{ErdRenyiEvol}, or
some books, e.g., Bollob\'as,  \cite{BolloRGbook}, Janson, \L{u}czak and
Ruci\'nski,  \cite{JansonLuczRuc}, Molloy and Reed  \cite{MolloyReedBook}.

\begin{theorem}[Erd\H{o}s-R\'enyi First Moment method]
\label{ErdRenyEv}
 Let
$\LL=\{L_1,\dots,L_t\}$ be a family of graphs, and let
\beq{FirstMoment}
c=\max_j \min_{H\subseteq L_j} {v(H)\over e(H)},
\qquad
\ga=\max_j \min_{H\subseteq L_j} {v(H)-2\over e(H)-1}
,\eeq where the minimum is taken only for subgraphs $H$ where the denominator
is positive.
\dori (a) Let $\Gn$ be a graph of order $n$ chosen uniformly, at random, from
graphs with
$E_n$
edges.
For every $\eps>0$ there exists a $\de>0$ such that if
$E_n<\de n^{2-c}$, then
the probability that $\Gn$
contains at least one $L\in\LL$ is at most $\eps$.
\dori
(b) If we know only $E_n<\eps n^{2-\gamma}$, then the probability that
$\Gn$ contains at least $\half E_n$ copies of $L\in\LL$ is at most $\eps$.
\end{theorem}

This implies that \beq{LowerExt} \ext(n,\LL)>c_\LL n^{2-\gamma}\ge c_\LL
n^{2-c}\eeq with $c\ge\gamma>0$ defined above.

\begin{remarks}

(a) A graph $L$ is called balanced if the minimum in \eqref{FirstMoment}, for $c$,
 is achieved for $H=L$. Erd\H{o}s and R\'enyi formulated their result containing
Theorem~\ref{ErdRenyEv}(a) only for balanced graphs $L$. The part we use is
trivial from their proof.

(b) Later Bollob\'as extended the Erd\H{o}s-R\'enyi theorem to arbitrary $\LL$.

(c) Gy\H{o}ri, Rothschild and
  Ruci\'nski achieved the generality by embedding any graph into a
  balanced graph  \cite{GyoriRothRuc}.

\end{remarks}

\begin{corollary}
\label{SuperLinX}

If a finite $\LL$ contains no trees,\footnote{neither forests} then
for some $c_\LL>0$, $\ext(n,\LL)\ge c_\LL n^{1+c}$.
\end{corollary}

Mostly the weaker Theorem~\ref{ErdRenyEv}(a) implies Corollary
\ref{SuperLinX}: it does, whenever $\LL$ is finite and each $L\in\LL$ contains
at least two cycles in the same component. However, for cycles we need
the stronger Theorem \ref{ErdRenyEv}(b).

For example, for $L=K_{a,b}$ we have $c=a^{-1}+b^{-1}$.
Then, for $c_0$ sufficiently small, the probability that a graph $\Gn$ with
$c_0n^{2-c}$ edges does not contain $K_{a,b}$ is positive. Hence
$$\ext(n,K_{a,b}) \ge c_0n^{2-(1/a)-(1/b)}.$$

Comparing this with the K\H ov\'ari-T. S\'os-Tur\'an theorem (Theorem
\ref{KovSosTurThB}), we see that the exponent is sharp there, in some sense,
if $a$ is fixed while $b\to\infty$.

\Proof of Theorem \ref{ErdRenyEv}. Consider the random graph $\Gn$
with $n$ labeled vertices, in which each edge occurs independently,
with the same probability $p$. For each $L_j$, choose a subgraph $H_j$
which attains the inner minimum for $\gamma$, in
\eqref{FirstMoment}. Let $h_j:=v(H_j)$, $e_j:=e(H_j)$, and let
$\alpha_j$ denote the number of copies of $H_j$ in $K_{h_j}$, and
$\beta_j$ denote the expected number of copies of $H_j$ in $\Gn$.

Clearly, $K_n$ contains $\alpha_j{n\choose h_j}$ copies of $H_j$. For
each copy $H$ of $H_j$, define a random ``indicator'' variable
$k_H=k_H(\Gn)=1$ if $H\subseteq \Gn$, and $0$ otherwise. Since the
number of copies of $H_j$ in $\Gn$ is just $\sum_{H\subseteq K_n}k_H$,
therefore, if $\EE$ denote the
expected value, then
$$\beta_j= 
\sum_{H\subseteq K_n} \EE(k_H)=\alpha_j{n\choose h_j}p^{e_j}.$$
Summing over $j$ and taking $p=c_1n^{-c}$, (for some $c_1\in(0,1)$) we get
$$\sum_j \beta_j \le t \max
\alpha_j{n\choose h_j}p^{e_j}\le t \max c_1n^{h_j-ce_j}=tc_1n^{2-c}.$$
Now let $\eta(\Gn)=e(\Gn)-\sum_j\beta_j$. Then, for $c_1$ sufficiently small, the
expected value is $$\EE(\eta(\Gn))>\half {n\choose 2}p > {1\over
5}c_1n^{2-c}.$$ Hence there exists a $\Gn$ with
$\eta(\Gn)>{1\over 5}c_1n^{2-c}$. Delete an edge from each $H_j$ in this
$\Gn$. The resulting graph contains no $L_j$, and has at least
${1\over 5}{n\choose 2}p\ge {1\over 11} c_1n^{2-c}$ edges, completing the
proof. \Qed

\begin{remarks}[How did the probabilistic methods start?]\label{ProbabStarts}

Mostly we write that applications of the Random Graphs (probabilistic
method) started when Erd\H{o}s (disproving a conjecture of
Tur\'an on the Ramsey Numbers) proved the existence of graphs $\Gn$
without complete subgraphs of order $2\log n$ and without independent
sets of size $2\log n$.
\begin{enumerate}\dense
\item Erd\H{o}s himself remarks (e.g., in  \cite{ErdSmoleProb}) that perhaps
 Szele was the first who applied this method in Graph Theory. (Erd\H{o}s --
 in his birthday volume  \cite{ErdFavourTheor} -- also mentions an even
 earlier application of J. Er\H{o}d but we did not succeed in locating that
 source.)
\item Perhaps the earliest case of applying probabilistic methods
 was that of Paul Tur\'an's proof of the Hardy-Ramanujan Theorem
 \cite{TuranHardyRam},
 where -- reading the paper -- it is obvious that Tur\'an gave
 a probabilistic proof of a beautiful and important theorem, using the
 Chebishev inequality. However, either Tur\'an did not realize that
 this is an application of the probabilistic method or he did not
 wish to burden the reader with that.
\item An important application of the probabilistic methods was that of Claude Shannon, when he constructed random codes.
\end{enumerate}

\end{remarks}

Applying Theorem \ref{ErdRenyEv} to some families of cycles we obtain

\begin{corollary}\label{CycleRandCons}
 For some constant $c_m>0$, $$\ext(n,\{C_3,\dots,C_m\})
\ge c_mn^{1+{1\over m-1}}.$$
\end{corollary}

Erd\H{o}s' even cycles theorem asserts that
$\ext(n,C_{2t})=O(n^{1+(1/t)})$, and this upper bound is probably sharp.
\footnote{The reference is missing here, since Erd\H{o}s did formulate
 this theorem but never have published a proof of it, as far as we
 know.} The random method (that is, Theorem \ref{ErdRenyEv}) yields
a lower bound of $cn^{1+{1\over 2t-1}}$ , a weaker result. Simonovits thinks
that it is unlikely that Theorem \ref{ErdRenyEv} ever yields a sharp
bound for a finite family. \footnote{Some related results of
 G. Margulis, and A. Lubotzky, R. Phillips and P. Sarnak will be
 discussed in Section \ref{LuboPhilSarnakS}.}

Corollary \ref{CycleRandCons} is used in the next section
 to prove that $\ext(n,\LL)=O(n)$ if and only if contains
a tree or forest.

\subsection{Classification of extremal problems}\label{subs:Classification}

The extremal graph problems can be classified in several ways. Here
we shall speak of (a) {\emR non-degenerate}, (b) {\emR degenerate} and
(c) {\emR linear} extremal problems.

For Case (a) Theorem \ref{AsymptStrTh} provides an appropriately good
description of the situation. In Case (b) $p(\LL)=1$. Here the ``main
term'' disappears, $(1-{1\over p})=0$;
therefore ``the error terms dominate''.
Case (c) will be discussed here shortly and in Sections \ref{TreeExtremal}
and \ref{InfiFami} in more details.

The classification immediately follows from the following theorems:

\begin{theorem}\label{DegenCriter}
  $\ext(n,\LL)=o(n^2)$ if and only if $\LL$ contains a bipartite
 graph. Actually, if $\LL$ contains a bipartite graph then
 $\ext(n,\LL)=O(n^{2-c})$ for, e.g., $c=2/v(L)$ for any bipartite
 $L\in\LL$. If $\LL$ does not contain bipartite graphs, then
 $\ext(n,\LL)\ge \turtwo n$.
\end{theorem}

\begin{theorem}\label{TreeUpper}

For finite $\LL$,
$\ext(n,\LL)=O(n)$ if and only if $\LL$ contains a tree, or a forest.
If $L\in\LL$ is a tree or a forest, then, for $v(L)\ge 3$,
\beq{TreeFolk} \ext(n,\LL)<(v(L)-2)n.
\eeq
\end{theorem}

\begin{theorem}[Erd\H{o}s \cite{ErdGraphProbA,ErdGraphProbB}]\label{ClassiBX}

If $\LL$ is finite and no $L\in\LL$ is a tree, then
 $\ext(n,\LL)>n^{1+c_\LL}$ for some $c_{\LL}>0$.
\end{theorem}

\begin{theorem}[Erd\H{o}s \cite{ErdGraphProbB}, Bondy and Simonovits \cite{BondySim}]\label{Girth}

Given an integer $k$, for some constants $c_k,\tilde c_k>0$,
\beq{uuppvv}
c_kn^{1+{1\over 2k-1}}<
\ext(n,\{C_3,\dots,C_{2k}\})\le\ext(n,C_{2k})\le \tilde c_kn^{1+{1\over k}}.\eeq
\end{theorem}

\Proof of Theorems \ref{DegenCriter}, \ref{TreeUpper}, and
\ref{ClassiBX}.
If there is a bipartite $L\in\LL$, then Theorem
\ref{KovSosTurThB} implies the sharper upper bound of Theorem
\ref{DegenCriter}.
Indeed, for $v=v(L)$,
by $L\subseteq K([v/2],v)$,
we have,
$$\ext(n,\LL)\le\ext(n,L)\le \ext(n,K([v/2],v))<\half\root v\of
{2v}\cdot n^{2-(2/v(L))}=O(n^{2-c}).$$

If the minimum chromatic number $p=p(\LL)\ge 3$, then
$\Tnp$ contains no forbidden $L\in\LL$. Therefore
$$\ext(n,\LL)>e(T_{n,2})\ge e(\Tnp)=\Tur(p){n\choose2}+O(n).$$
Actually, $e(T_{n,2})=\turtwo n$.
This completes the proof of Theorem \ref{DegenCriter}.\Qed

It is easy to show that if $\Gn$ has minimum degree at least $r-1$,
then it contains every tree $T_r$ (by induction on $r$). An induction
on $n$ yields \eqref{TreeFolk}, implying half
of Theorem \ref{TreeUpper}, when $\LL$ contains a tree (or a
forest). If $\LL$ is finite and contains no trees, i.e., all the
forbidden graphs contain some cycles, then we use Theorem \ref{Girth},
or simply Corollary \ref{SuperLinX}, proved by probabilistic
methods.\footnote{There are also deterministic proofs of Corollary
 \ref{SuperLinX}, e.g., via the Margulis--Lubotzky--Phillips--Sarnak
 construction of Ramanujan graphs, see Construction
 \ref{LuboPhilSarnak}.}

\begin{remark}[Infinite families]

For infinite families the situation is different: if e.g. $\C$ is the family of all cycles, then
$\ext(n,\C)=n-1$: all graphs but the forests are excluded.
There are many further families without trees where the extremal
number is linear, see Section \ref{InfiFami}.
\end{remark}

\Proof of Theorem \ref{Girth}.
The lower bound comes from a random graph argument of Erd\H{o}s.
Concentrate on the upper bound.
If we are not interested in the value of the constant, then
we can basically use the following argument: Take a graph $\Gn$ with
$cn^{1+\alpha}$ edges. Delete its minimum degree vertex, then the
minimum degree vertex in the remaining graph, etc. At the end we get a
$\Gm$ with minimum degree at least $c_1m^{\alpha}$.
In the obtained graph $\Gm$ fix a vertex $x$ and denote by $S_j$ the set of vertices at distance $j$ from $x$.
If $\girth(\Gn)>2k$, -- as we assumed -- then
{\emR basically} $|S_j|>\mindeg(\Gm)\cdot|S_{j-1}|$. Hence
$m>|S_k|>c_1^km^{\alpha k}$. So $\alpha\le 1/k$.\qed

Assume for a second that $\Gn$ itself is asymptotically regular:
$${\mindeg(\Gn)\over \maxdeg(\Gn)}\to 1.$$
Then the previous argument asserts that $d:=\mindeg(\Gn)< n^{1/k}.$
Therefore
$$e(\Gn)\le \left(\half+o(1)\right) n d\approx \half n^{1+{1\over k}}.$$

We shall return to the case of excluded trees, namely, to the
Erd\H{o}s-S\'os conjecture on the extremal number of trees, and to the
related Koml\'os-S\'os conjecture in Section~\ref{TreeExtremal}. One
final question could be if $\ext(n,\LL)$ can be sublinear. This is
answered by the following trivial result.

\begin{theorem}\label{Bounded}

If $\LL$ is finite and $\ext(n,\LL)<[n/2]$, then $\ext(n,\LL)=O(1)$.
\end{theorem}

\proof. Consider $n/2$ independent edges: this must contain an
$L_1\in\LL$. Hence, there is an $L_1\in\LL$ contained in the
union of $t$ independent edges, for some $t$. Also, there exists an
$L_2\subseteq K(1,n-1)$. Hence an extremal graph $\Sn$ has bounded
degrees and bounded number of independent edges. This proves
\ref{Bounded}.\qed

Theorem \ref{Bounded} easily extends to hypergraphs.

\subsection{General conjectures on bipartite graphs}

We have already formulated Conjecture \ref{RacConj} on the rational
exponents. We have to remark that for hypergraphs this does not hold:
the Behrend construction  \cite{Behrend} is used to get lower bounds in the
Ruzsa--Szemer\'edi Theorem, (Thm \ref{RuzsaSzemTh}), showing
that there is no rational exponent in that case. Yet, Erd\H{o}s and
Simonovits conjectured that for ordinary graphs there is. One could also conjecture the inverse extremal problem:

\begin{conjecture}\label{RatiExpo}

For every rational $\alpha\in(0,1)$ there is a finite $\LL$ for which
$c_1n^{1+\alpha}<\ext(n,\LL)<c_2n^{1+\alpha}$, for some constants $c_1,c_2>0$.
\end{conjecture}

The third conjecture to be mentioned here is on ``compactness''  \cite{ErdSimComp}:

\begin{conjecture}

For every finite $\LL$ there is an $L\in\LL$ for which
$\ext(n,\LL) > c\cdot\ext(n,L)$, for some constants $c=c_{\LL}>0$.
\end{conjecture}

\section{Excluding complete bipartite graphs}

\subsection{Bipartite $C_4$-free graphs and the Zarankiewicz problem}\label{C4ZaraS}

Tur\'an type extremal results (and Ramsey results as well) can often
be applied in Mathematics, even outside of Combinatorics. Tur\'an
himself explained this applicability by the fact that -- in his
opinion -- the extremal graph results were generalizations of the
Pigeon Hole Principle.

Recall that
 $Z(m,n,a,b)$ denotes the maximum number of 1's in an
$m\times n$ matrix not containing an $a\times b$ minor consisting
exclusively of 1's. In 1951 Zarankiewicz~ \cite{Zarank} posed the
problem of determining $Z(n,n,3,3)$ for $n\le 6$, and the general
problem has also become known as {\emR the problem of
 Zarankiewicz}.\footnote{In Graph Theory two problems are connected
 to Zarankiewicz' name: the extremal problem for matrices that we
 shall discuss here and the Crossing Number conjecture which is not
 our topic. Actually, the crossing number problem comes from Paul
 Tur\'an, see  \cite{TurWel}.} Obviously, $Z(m,n,1,b)=m(b-1)$ (for
$n\ge b-1$). Observe that $Z(m,n,a,b)=\ext^*(m,n,\KK ab)$
(where $\ext^*(m,n,\LL)$ was defined following Remark \ref{OurProblems}.)
Considering the adjacency matrix of a $K_{a,b}$-free graph on $n$
vertices we get $2 \ext(n, K_{a,b})\le Z(n,n,a,b)$. We will use this
upper bound many times.

We will see that the easy upper bound in Theorem \ref{KovSosTurThB} is
pretty close to the truth for $a\le 2$. Actually, K\H ov\'ari,
T. S\'os and Tur\'an  \cite{KovSosTur} proved an upper bound for the
Zarankiewicz function \beq{KST_Zaran} Z(m,n,a,b) \le \root
a\of{b-1}\cdot mn^{1-(1/a)} +(a-1)n \eeq which was slightly improved by
Zn\'am~ \cite{ZnamA},  \cite{ZnamB},  (he halved the last term to
$(a-1)n/2$ in the case of $m=n$) and Guy~ \cite{GuyTihany}.

\def\cS{{\cal S}}

A bipartite graph $G[M,N]$ where $|M|=m$, $|N|=n$ is $C_4$-free if its
``bipartite'' $m\times n$ adjacency matrix contains no $2\times 2$ full 1
submatrix.\footnote{Here the ``bipartite adjacency matrix''
 $A=(a_{ij})_{m\times n}$ is defined for a bipartite graph $G[U,V]$ and
 $a_{ij}=1$ if $u_j\in U$ is joined to $v_j\in V$, otherwise $a_{ij}=0$.} In
other terminology, the hypergraph defined by the rows of this matrix is linear,
and their hyperedges pairwise meet in at most one element. There is an
important class of such hypergraphs, the Steiner $k$-systems $S(n,k,2)$. A
family $\cS$ of $k$-subsets of an $n$-set $N$ is a Steiner $k$-system if
every pair of elements is covered exactly once. For such an $\cS$, clearly,
$|\cS|=m={n\choose 2}/{k \choose 2}$. Such families are known to exists for
$(m,n,k)=(q^2+q+1, q^2+q+1, q+1)$ (called finite projective planes of order
$q$),
 and $(m,n,k)=(q^2+q, q^2, q)$ (affine planes) whenever $q$ is a power of
a prime. Also for any given $k$ there exists an $n_0(k)$ such that
$S(n,k,2)$ exists for all {\emR admissible} $n> n_0(k)$, i.e., when
$(n-1)/(k-1)$ and $n(n-1)/k(k-1)$ are integers (Wilson's existence
theorem~ \cite{WilsonExistence}).

K\H{o}v\'ari, T. S\'os and Tur\'an  \cite{KovSosTur} proved that

\begin{theorem}

$Z(n,n,2,2)=(1+o(1))n^{3/2}$, and
\beq{KSTZ} Z(n,n,2,2)<[n^{3/2}]+2n.\eeq
Further, if $p$ is a prime, then
$$Z(p^2+p,p^2,2,2)=p^3+p^2.$$
\end{theorem}

Reiman  \cite{ReimanZara} returned to this topic,
(see also  \cite{ReimanLapok}), slightly improving \eqref{KSTZ}

\begin{theorem}[Reiman  \cite{ReimanZara}]\label{th:Reiman}

\beq{ReimZa}Z(m,n,2,2)\le \half\left(m+\sqrt{m^2+4mn(n-1)}\right).\eeq
\end{theorem}

For large $m,n\to\infty$, and $m=o(n^2)$,
this yields
$$Z(m,n,2,2)\le \left(\half+o(1)\right)n\sqrt{m}.$$
Further, for $m=n$, we get
\beq{ReimanZ}Z(n,n,2,2)\le \half n\left(1+\sqrt{4n-3}\right)\approx n\sqrt{n}.\eeq
Reiman also
provides infinitely many graphs, using Finite Geometries,
 showing the sharpness of \eqref{ReimZa} and \eqref{ReimanZ}.
We have equality when $m=n(n-1)/k(k-1)$ and a Steiner system $S(n,k,2)$ exists.
Thus he determined the case
 \beq{ReimanZq}Z(n,n,2,2)= \half n\left(1+\sqrt{4n-3}\right)=(q^2+q+1)(q+1)\eeq
for $m=n=q^2+q+1$ when a projective plane of order $q$ exists.
Actually, in  \cite{ReimanLapok}, Reiman also speaks about
Zarankiewicz-extremal graphs connected to incidence-graphs of higher
dimensional finite geometries.

Since Reiman's theorem the theory of finite geometries developed tremendously.
We cite here a recent result whose proof used the most modern tools and
 stability results.

\begin{thm}[Dam\'asdi, H\'eger, and Sz\H onyi~\cite{DamasdiHegerSzonyiZaran}]
Let $q\ge 15$, and $c\le q/2$. Then
$$
Z(q^2 + q + 1 - c, q^2 + q + 1, 2,2) \le (q^2 + q + 1 - c)(q + 1).
  $$
Equality holds if and only if a projective plane of order $q$ exists.
Moreover, graphs giving equality are subgraphs of the bipartite incidence graph
 of a projective plane of order $q$ obtained by omitting $c$ rows of its incidence matrix.
  \end{thm}

They proved many more exact results when a projective plane of order $q$ exists.
The extremal configurations are submatrices of the incidence matrix of a projective plane.
\begin{eqnarray*}
Z(q^2 + c, q^2 + q, 2,2) &=& q^2(q + 1) + c q \quad (0 \le c \le q + 1),\\
Z(q^2 -q + c, q^2 + q -1, 2,2) &=& (q^2 - q)(q + 1) + c q \quad (0 \le c \le 2q),\\
Z(q^2 - 2q + 1 + c, q^2 + q - 2, 2,2) &=& (q^2 - 2q + 1)(q + 1) + cq \quad (0 \le c \le 3(q - 1)).
  \end{eqnarray*}

These refer to bipartite host graphs.
As we will see later, such exact results are rare
 for the general (non-bipartite) case.
To estimate $\ext(n,C_4)$ seems to be harder, because the
 corresponding $0$-$1$ matrices, the incidence matrix of a graph, should be symmetric.

\subsection{Finite Geometries and the $C_4$-free graphs}\label{FinGeomA}

The method of finite geometric constructions is very important and
powerful in combinatorics. In particular, it is often the best way to
obtain lower bounds. It is for this reason that we include this
section.

We give several constructions: the first two show that the K\H
ov\'ari-T. S\'os-Tur\'an theorem (Theorem \ref{KovSosTurThB}) is sharp
for both $K_{2,2}$ and $\HK33$.

\begin{remark}\label{K33SharpnessR}
 When we write that an upper bound is sharp, mostly we mean that it is
 sharp up to a multiplicative constant: it yields the correct
 exponent. There are a few exceptions, where sharpness means that the
 ratio of the upper and lower bounds tends to 1. This is the case for
 $C_4=\KK22$ and we have this also for $\KK33 $. Here, however, the
 matching upper bound for Construction \ref{BrownConstr} below is
 given not by Theorem~\ref{KovSosTurThB} but by the F\"uredi
 improvement  \cite{FureK33}.
\end{remark}

Perhaps the application of finite geometries in Extremal Graph Theory
started in the Erd\H{o}s paper, with the construction of Eszter Klein
 \cite{ErdTomsk}, to prove the sharpness of Theorem \ref{TomskC4}. The
expression ``Finite Geometry'' was not mentioned there. We skip the
description of this whole story, since it was described in several
places, e.g.,  \cite{SimErdosInflu},
 \cite{SimPragNew}.

Much later, Erd\H{o}s and R\'enyi  \cite{ErdRenyiDiam} used finite geometry for
a diameter-extremal problem.  This is a very large area, connected to our
problems, yet we have to skip it.  The interested reader is refered to
 \cite{ErdRenyiDiam}, (translated into English in  \cite{RenyiColl}).

Sharp extremal graph results were obtained by
Reiman  \cite{ReimanZara} and a Polarity Graph was used in
 \cite{ErdRenyiSos} and  \cite{BrownThom} to give asymptotically sharp
lower bound on $\ext(n,C_4)$.
This lower bound can also be found in  \cite{ErdRenyiDiam}, {\bf implicitly}:
Erd\H{o}s and R\'enyi considered the diameter-extremal problem, and do mention
the properties of this graph.

The real breakthrough came by the Erd\H{o}s-R\'enyi-T. S\'os
paper  \cite{ErdRenyiSos}, (sharp lower bound for $C_4$)
and by the Brown paper  \cite{BrownThom},
providing asymptotically sharp lower bounds for $\ext(n,C_4)$
and for $\ext(n,\HK33)$. (See Remark \ref{K33SharpnessR}.)

We know from Theorem \ref{KovSosTurThB} that
$\ext(n,C_4)\le \half n^{3/2}+o(n^{3/2}),$ but is this result sharp?
In analyzing the proof, we realize that if it is sharp (that is, if there are
infinitely many graphs $\Gn$ not containing $C_4$ and having $\approx \half n\sqrt n$ edges),
then almost all degrees are $\approx \sqrt n$ and almost every
pair of vertices must have a common neighbor (and no pair has two). This
suggests that the neighborhoods $N(x)$ behave much like the lines in a
projective plane, in that the following statement ``almost'' holds: any two
vertices lie in a common set, and any two sets intersect in one
vertex.


\begin{theorem}[Erd\H{o}s--R\'enyi--T. S\'os
 \cite{ErdRenyiSos}, and Brown  \cite{BrownThom}, see also
 \cite{KovSosTur}]\label{ErdRenyiSosTh}
\index{Erdos-Renyi-Sos} \index{Brown}
$$\ext(n,C_4)= \half n^{3/2}+O(n^{3/2-c}).$$
\end{theorem}

For the lower bound for $\ext(n,C_4)$ we use the following

\begin{construction}\label{PolarAff}

 Let $p$ be a prime, $n=p^2-1$. Construct
a graph as follows: the vertices are the $p^2-1$ non-zero pairs $(x,y)$ of
residues (modulo $p$), and $(x,y)$ is joined to $(a,b)$ by an edge if
$ax+by=1$. (This graph may contain loops, but we simply delete them.)
\end{construction}

With $n=p^2-1$, the resulting graph $\Hn$ has the necessary properties to show
the sharpness of Theorem \ref{KovSosTurThB} for $C_4$:

(a) for a given pair $(a,b)$, mostly there are $p$ solutions to $ax+by=1$, so
that, even after the loops are deleted, there are at least $\half (p^2-1)(p-1)$
edges in $\Hn$ and hence $e(\Hn)>\half n^{3/2}-n$;

(b) if $\Hn$ had a 4-cycle with vertices $(a,b),$ $(u,v),$ $(a',b')$ and
$(u',v'),$ then the two equations $ax+by=1$ and $a'x+b'y=1$ would have two
solutions, which is impossible. Since the primes are ``dense'' among the integers,
this completes the proof of the the sharpness of Theorem \ref{KovSosTurThB}
for $a=b=2$.

\begin{remark}

An alternative possibility is to use the much more symmetric polarity
graph of the projective plane (we explain this in the next section):
here we used the Affine Geometric Variant
because {\emR here} we did not wish to use anything from Projective Geometry.
\end{remark}

\subsection{Excluding $C_4$: Exact results}\label{TruncPol}

The polarity graph
\footnote{ These $C_4$-free graphs were studied earlier in finite geometry.
The bipartite point-line incidence graph appeared in Levi's book (1942)
 and polarity graphs (modulo loops) obtained from {\it Levi graph}
 had been described already by Artzy (1956).
For more details and references see Bondy~\cite{BondyHB}.},
used in  \cite{ErdRenyiDiam},
was also used in  \cite{ErdRenyiSos} and  \cite{BrownThom} to prove that
\beq{C4PolarLower}\ext(n,C_4)\ge \half q(q+1)^2,\Text{for} n=q^2+q+1.\eeq
if $q$ is a prime power.

\begin{construction}[The Polarity Graph from the finite field]\label{PolarityGraph}

Assume that $q$ is a prime power. Consider the Finite Field $GF(q)$.
The vertices of our graph are the equivalence classes of the non-zero
 triples $(a,b,c)\in GF(q)^3$
 where two of them, $(a,b,c)$ and $(a',b',c')$ are considered the same if
 $(a',b',c')=\lambda(a,b,c)$ for some
$\lambda\ne0$. There are $(q^3-1)/(q-1)=q^2+q+1$ such classes.
Further, the equivalence class of $(a,b,c)$ is connected by an edge to
 the class of $(x,y,z)$ if $ax+by+cz=0$.
Finally, we delete the $q+1$ loops, i.e. those edges, where
$a^2+b^2+c^2=0$.
This graph is $C_4$-free and it has $\frac{1}{2}(n(q+1) -(q+1))$ edges.
\end{construction}

In general, a polarity corresponds to a {\emR symmetric} incidence
 matrix of a
 finite plane of size $(q^2+q+1)\times (q^2+q+1)$.
According to a theorem of Baer  \cite{BaerPolar} such a matrix has at least
 $q+1$ non-zero elements in its diagonal.
Therefore using the polarity graph we cannot avoid losing on loops.
This way Erd\H{o}s, R\'enyi and S\'os~ \cite{ErdRenyiSos} showed that
 indeed $\ext(q^2+q+1, C_4)< \half n(q+1)$.
Yet one could hope to get a better construction.
Erd\H{o}s conjectured  \cite{ErdRome},  \cite{ErdRome1976-35} that
there are no better constructions, that is,
\eqref{C4PolarLower} is sharp if $n=q^2+q+1$, ($q$ is a prime power).

F\"uredi settled this conjecture in the following sense: First he
proved  \cite{FureC4JCT83} that if $q=2^k$, then Erd\H{o}s' conjecture
holds. Next he 
 settled the case $q\ge q_0$.
Later he found a much shorter proof of the weaker assertion that
the Polarity graphs are extremal;
however, this shorter version did not give the extremal structure.
So F\"uredi published the shorter version, while the longer version can be found on his homepage.

\begin{theorem}[F\"uredi  \cite{FureC4JCT96},  \cite{FureC4Prep}]

If $q\ne 1,7,9,11,13$ and $n=q^2+q+1$, then $\ext(n,C_4)\le \half q(q+1)^2$
 and for $q> 13$ the extremal graphs are obtained from a polarity of a finite projective plane.
Hence if $q>13$ is a prime power, then $\ext(n,C_4)=\half q(q+1)^2$.
\end{theorem}
The second part of this result probably holds for $q\in \{7,9,11,13\}$, too.

Recently a new sharp construction has been found for $n=q^2+q$.

\begin{theorem}[Firke, Kosek, Nash and Williford  \cite{FirkeWilli}]
 Suppose that $q$ is even, $q> q_0$. Then
$$\ext(q^2+q,C_4)\le \half q(q+1)^2-q.$$
Consequently, if $q> q_0$, $q=2^k$ and $n=q^2+q$ then
 $\ext(n,C_4)= q(q+1)^2-q$.
\end{theorem}

They also announced that in a forthcoming paper they show that for all but
 finitely many even $q$, any $\Sn\in\EXT(q^2+q,C_4)$ is derived from an orthogonal
 polarity graph by removing a vertex of minimum degree
 (the 1-vertex-truncated Polarity graph, see Construction \ref{PolarityGraph}).
This result shows a kind of stability of the Polarity graph.
More generally, McCuaig (private communication, 1985) {\bf conjectured}
that each extremal graph is a subgraph of some polarity graph.
So this is true for infinitely many cases, but one of the
present authors strongly disagrees and he believes just the opposite
that for e.g., $n=q^2+q+2$ maybe the extremal graphs are obtained by
{\emR adding} an extra vertex and some edges to a polarity graph.

\begin{remark}

W. McCuaig calculated $\ext(n,C_4)$ for $n\le 21$ (unpublished letter, 1985).
Clapham, Flockart and Sheehan determined the corresponding extremal graphs
 \cite{ClapFlocShee}, and Yuansheng and Rowlinson  \cite{YuangRowliC4}, -- using
computers, -- extended
these results to $n\le 31$.
(They also determined the graphs in $\EXT(n,C_6)$ for $n\le 26$,
 \cite{YuangRowliC6}.)
 Garnick,  Kwong, Lazebnik, and Nieuwejaar
  \cite{GarnLazebKwong},  \cite{GarnickNieu}
determined the values of $\ext(n, \{C_3,C_4\})$ for all $n \leq 30$.

\end{remark}

\subsection{Excluding $K(2,t+1)$, $t>1$}\label{Zara2S}

A slightly sharper form of the upper bound \eqref{KST_Zaran} was
 presented by Hylt\'en-Cavallius  \cite{Hylten}
\beq{HylCav}Z(m,n,2,k)\le \half n+\left\{(k-1)nm(m-1)+{1\over 4}n^2\right\}^{1/2}.
\eeq
Obviously, for fixed $k$ and large values of $n,m$, if $n=o(m^2)$, then
the right hand side of \eqref{HylCav}
is $\approx \sqrt{k-1}m\sqrt n.$
Using again the observation $2\ext(n, K_{2,t+1})\le Z(n,n,2,t+1)$
 one obtains the upper bound
\beq{K2tupper}
\ext(n, K_{2,t+1})\le \half n\sqrt{tn-t+1/4}+(n/4).
\eeq
The following theorem shows that the above (easy) upper bound is
 the best possible asymptotically.

\begin{theorem}[F\"uredi 
   \cite{FureK2t}]

For any fixed $t\ge 1$
\beq{MorsSymm}
 \ext (n,K_{2,t+1})= \half\sqrt{t} n^{3/2} + O(n^{4/3}). \eeq
 \end{theorem}

To prove this Theorem one needs an appropriate lower bound, a
 construction.
Let $q$ be a prime power such that $(q-1)/t$ is an integer.
We will construct a $K_{2,t+1}$-free graph $G$ on $n=(q^2-1)/t$ vertices
 such that every vertex has degree $q$ or $q-1$.
We will explain this below (Construction~\ref{constK2t}).
Then $G$ has more than $(1/2)\sqrt{t} n^{3/2}-(n/2)$ edges.
The gap between the lower and upper bounds
 is only $O(\sqrt n)$ for $n=(q^2-1)/t$.
The lower bound for the Tur\'an number for all $n$ then follows from the
 fact that such prime powers form a dense subsequence among the integers.
This means that for every sufficiently large $n$ there exists a prime $q$
 satisfying $q\equiv 1$ (mod $t$) and
 $\sqrt{nt}-n^{1/3}< q< \sqrt{nt}$ (see  \cite{HuxleyIwaniec}).

Construction~\ref{constK2t} below is inspired by constructions of Hylt\'en-Cavallius
 and M\"ors given for Zarankiewicz's problem $Z(n,n,2,t+1)$.

\begin{theorem}[Hylt\'en-Cavallius  \cite{Hylten}]

$Z(n,n,2,3)=\sqrt2n^{3/2}+o(n^{3/2}).$
Also
$$\sqrt {\lfloor k/2\rfloor}\le \liminf_{n\to\infty}
{Z(n,n,2,k)\over n^{3/2}}.$$
\end{theorem}

\begin{theorem}[M\"ors  \cite{Mors}]

For all $t\ge 1$,
$${Z(n,n,2,t+1))\over n^{3/2}}\to \sqrt{t},\Text{as}n\to\infty. $$
\end{theorem}

The topic was so short of constructions that, as a first
 step, P. Erd\H{o}s   \cite{ErdRome,ErdosKeszthely73} even proposed the problem
 whether $\lim_t (\liminf_n \ext(n,K_{2,t+1})n^{-3/2})$ goes to $\infty$
 as $t\to \infty$.

\begin{remark}
 Here we see three distinct quantities, exactly as it is
described in Problem \ref{OurProblems}.
$Z(m,n,2,t+1)=\ext^*(m,n,K_{2,t+1})$, estimated from below by M\"ors, by a construction,
and $\ext(m,n,K_{2,t+1})$ estimated by F\"uredi by the same construction.
F\"uredi showed that the matrix of M\"ors contains
neither a $(t+1)\times 2$ submatrix, nor a $2\times (t+1)$ submatrix of 1's;
finally, F\"uredi, slightly changing the definitions in M\"ors's
construction
extended this ``asymmetric matrix'' result to the symmetric case
and provided a non-bipartite graph, proving \eqref{MorsSymm}.
\end{remark}

\begin{construction}\label{constK2t}

Let $GF(q)$ be the $q$-element finite field, and let
\mbox{$h\in GF(q)$} be an element of order $t$.
This means, that $h^t=1$ and the set $H=\{ 1, h, h^2, \dots , h^{t-1}\}$
 form a $t$-element subgroup of $GF(q)\setminus \{ 0 \}$.
For $q\equiv 1\pmod t$ such an element $h\in GF(q)$ always exists.

We say that $(a,b)\in {GF(q)}\times {GF(q)}$, $(a,b)\neq (0,0)$
 is equivalent to $(a',b')$, in notation $(a,b)\sim (a',b')$, if there exists
 some $h^\alpha\in H$ such that $a'=h^\alpha a$ and $b'=h^\alpha b$.
The elements of the vertex set $V$ of $G$ are the $t$-element equivalence
 classes of ${GF(q)}\times{GF(q)}\setminus (0,0)$.
The class represented by $(a,b)$ is denoted by $\left< a,b \right>$.
Two (distinct) classes $\left< a,b \right>$ and $\left< x,y \right>$ are joined by an edge
 in $G$ if $ax+by \in H$.
This relation is symmetric, and $ax+by\in H$, $(a,b) \sim (a',b')$, and
 $(x,y)\sim(x',y')$ imply $a'x'+b'y'\in H$.
So this definition is compatible with the equivalence classes.
 \end{construction}

For any given $(a,b)\in {GF(q)}\times{GF(q)}\setminus (0,0)$ (say,
 $b\neq 0$) and for any given $x$ and $h^{\alpha}$, the equation
 $ax+by=h^\alpha$ has a unique solution in $y$.
This implies that there are exactly $tq$ solutions $(x,y)$ with
 $ax+by\in H$.
The solutions come in equivalence classes, so there are exactly
 $q$ classes $\left< x,y\right>$.
One of these classes might coincide with $\left< a,b\right>$ so the degree of the
 vertex $\left< a,b\right>$ in $G$ is either $q$ or $q-1$.

We claim that $G$ is $K_{2,t+1}$-free.
First we show, that for $(a,b), (a',b')\in {GF(q)}\times{GF(q)}\setminus
 (0,0)$, $(a,b)\not\sim (a',b')$ the equation system
\beq{tag2} ax+by =h^{\alpha} \Text{and} a'x+b'y =h^{\beta}
  \eeq
 has at most one solution $(x,y)\in {GF(q)}\times{GF(q)}\setminus (0,0)$.
Indeed, the solution is unique if the determinant
 $\det\begin{pmatrix}a & b \cr a' & b'\end{pmatrix}$ is not 0.
Otherwise, there exists a $c$ such that $a=a'c$ and $b=b'c$.
If there exists a solution of (2) at all, then multiplying the
 second equation by $c$ and subtracting it from the first one we get
 on the right hand side $h^\alpha -ch^\beta=0$.
Thus $c\in H$, contradicting the fact that $(a,b)$ and $(a',b')$
 are not equivalent.

Finally, there are $t^2$ possibilities for $0\le \alpha, \beta <t$ in \eqref{tag2}.
The set of solutions again form $t$-element equivalence classes, so
 there are at most $t$ equivalence classes $\left< x,y\right>$ joint
 simultaneously to $\left< a,b\right>$ and $\left< a',b'\right>$. \Qed

Since then, there have been two additional almost optimal constructions,
strongly related to the Construction~\ref{constK2t} above.

\begin{construction}[Lazebnik, Mubayi~ \cite{LazebMubayi}]\label{co:K2tLazebnikMubayi}

Let $GF(q)^*$ be the finite field of order $q$ without the zero element.
Suppose $q\equiv 1$ (mod $t$) and let $H$ be the $t$-element multiplicative subgroup of $GF(q)^*$.
Define the graph $G^\times$ as follows.
Let $V(G^\times)= (GF(q)^*/H)\times GF(q)$.
For $\left< a\right>, \left< b\right>\in (GF(q)^*/H)$ and $x,y\in GF(q)$,
 make $(\left< a\right>, x)$ adjacent to $(\left< b\right>, y)$ if
 $x+y\in \left< ab\right>$.
\end{construction}

This graph (after deleting the eventual loops) is $K_{2,t+1}$-free and every vertex has degree
 $q-1$ or $q-2$.
Actually, Construction~\ref{co:K2tLazebnikMubayi} differs from
 Construction~\ref{constK2t} only in that its vertex set is smaller and
 instead of using the rule that $\left<a,b\right>$ is adjacent to $\left<x,y\right>$
 if $ax+by\in H$ they use the rule $ay+bx\in H$.
This change allows them to generalize it to multipartite hypergraphs.

 The following example works only if $t$ is a power of a prime, and $t|q$.

\begin{construction}[Lenz, Mubayi~ \cite{LenzMubayi}]\label{co:K2tLenzMubayi}

Suppose that $t$ divides $q$ and let $H$
be an additive subgroup of $GF(q)$ of order $t$. Define the graph
$G^+$ as follows. Let $V(G^+)= (GF(q)/H)\times GF(q)^*$. We will
write elements of $GF(q)/H$ as $\left<a\right>$. It is the additive
coset of $H$ generated by $a$, $\left<a\right> = \{h + a : h\in H\}$.
For $\left< a\right>, \left< b\right>\in (GF(q)/H)$ and $x,y\in
GF(q)^*$, make $(\left< a\right>, x)$ adjacent to $(\left< b\right>,
y)$ if $xy\in \left< a+b\right>$. (This, in fact, means that there
exists an $h\in H$ such that $xy=a+b+h$).
\end{construction}

\subsection{Excluding $K(3,3)$, and improving the upper bound}\label{Zara33}

The main result of this section is the description of the
asymptotically sharp value of $\ext(n,\KK 33)$.

\begin{theorem}[Brown  \cite{BrownThom} and F\"uredi  \cite{FureK33}]

$$\ext(n,\HK33)=\half n^{5/3}+O(n^{(5/3)-c})\Text{for some} c>0.$$
 \end{theorem}

The lower bound can be obtained from Brown's example (discussed below
as Construction~\ref{BrownConstr}) who gave a $(p^2-p)$-regular
$K_{3,3}$-free graph on $p^3$ vertices for each prime $p$ of the form
$4k-1$.

 Improving the upper bound in Theorem \ref{KovSosTurThB} F\"uredi showed that
Brown's example is asymptotically optimal.

\begin{theorem}[F\"uredi \cite{FureK33}]\label{FureImprove}

For all $m\ge a$, $n\ge b$, $b\ge a\ge 2$ we have
\beq{eq:KSTjav}
 Z(m,n,a,b)\le (b-a+1)^{1/a}mn^{1-(1/a)}+(a-1)n^{2-(2/a)}+(a-2)m. \eeq
\end{theorem}

For fixed $a,b\ge 2$ and $n,m\to \infty$ the first term is the largest one
 for $n=O(m^{a/(a-1)})$.
This upper bound is asymptotically optimal for $a=2$ and for
 $a=b=3$ ($m=n$).
We obtain
\beq{eq:K33ZFupper} \ext(n,K_{3,3})\le\half Z(n,n,3,3) \le \half n^{5/3}+ n^{4/3}+{1\over 2}n.
 \eeq
Alon, R\'onyai and Szab\'o~ \cite{AlonRonySzab} gave an example (discussed as
 Construction~\ref{co:AlonRonySzabo}) showing that
 \beq{eq:K33ARSzabo} \ext(n,K_{3,3}) \ge \half n^{5/3}+{1\over 3}n^{4/3}-C.
 \eeq
for some absolute constant $C>0$ for every $n$ of the form $n=p^3-p^2$, $p$ is a prime.
Their example shows that the upper bound \eqref{eq:K33ZFupper} (and \eqref{eq:KSTjav})
 is so tight that that we cannot leave out the second order term.
It would be interesting to see whether \eqref{eq:KSTjav} is tight
for other values of $a$ and $b$, too.

The first step of the proof of Theorem~\ref{FureImprove} is that
 given a $K_{a,b}$-free graph $G$, we apply the
 original bound \eqref{KST_Zaran} to the bipartite subgraphs $G[N(x), V\setminus N(x)]$ generated by the
 neighborhood of a vertex $x$ and its complement.

%
When Brown gave his construction, the matching upper bound of Theorem~\ref{FureImprove}
 was not known yet. He wrote that even the existence of
 $\lim_{n\to\infty} \ext(n,\HK33)/n^{5/3}$ was unknown.

\begin{construction}\label{BrownConstr}

 Let $p$ be an odd prime $n=p^3$ and $d\in GF(p)$, $d\neq 0$ a quadratic residue if
 $p$ is of the form $4k-1$ and $d$ be a non-residue otherwise.
Construct a graph $B_n$ whose vertices are the triples $(x,y,z)$ of residue
classes (modulo $p$) and whose edges join vertices $(x,y,z)$ and
$(x',y',z')$ if
\beq{Br3}(x-x')^2+(y-y')^2+(z-z')^2=d.\eeq
\end{construction}

It is easy to see that the graph $B_n$ has $\half n^{5/3}+O(n^{4/3})$ edges.
Given a vertex $(x',y',z')$, the equation \eqref{Br3} has $p^2-p$ solutions
 by a theorem of Lebesgue.
Thus $(x',y',z')$ has this many neighbors.

We claim that $B_n$ does not contain $\HK33$.
The geometric idea behind Construction~\ref{PolarAff} (concerning $C_4$-free graphs)
 was to join a point of the finite plane to the points of its ``polar''
 (with respect to the unit circle), and then to use the fact that two lines
 intersect in at most one point.
In contrast, the Brown construction
uses the fact that, if points
of the Euclidean space $\EE^3$ at distance $1$ are joined, then the resulting infinite
graph $G$ does not contain $\HK33$. This is easily seen as follows:
suppose $G$ does contain $\HK33$. Then the three points of one
color class cannot be collinear since no point is equidistant from
three collinear points. On the other hand, only two points are
equidistant from three points on a circle, and so $\HK33$ cannot
occur. There is one problem with this ``proof'': in finite fields
$\sum_i x_i^2=0$ can occur even if not all $x_i$'s are 0's. Therefore in
finite geometries, in some cases, a sphere can contain a whole line.
So here the geometric language must be translated into the language of analytic
geometry,
 and the right hand side of \eqref{Br3} (that is $d$) must be chosen
 appropriately.

\begin{theorem}[Nikiforov,  \cite{NikifZaran}]\label{NikifZaraTh}

For $b\ge a\ge 2$ let $k\in[0,a-2]$ be an integer. Then
$$Z(m,n,a,b)\le (b-k-1)^{1/a}mn^{1-(1/a)}+(a-1)n^{1+(k/a)}+km.$$
\end{theorem}

For $k=0$ we get back Theorem \ref{KovSosTurThB}, and 
substituting $k=a-2$ we obtain \eqref{eq:KSTjav}.
Nikiforov remarks that letting $k$ run from $0$ to $a-2$,
we may get the best results for various values of $k$ as the relation
of $m$ and $n$ varies, but we still have no constructions to substantiate this.
Nikiforov also proves results on the spectral radius. 

\subsection{Further applications of Algebraic Methods}\label{NormGraphS}

Most of the constructions providing sufficiently good lower bounds for
Bipartite Extremal Graph Problems are coming either from Geometry or
from Algebra\footnote{The Random Graph Methods are very nice but
 mostly they are too weak to provide sufficiently sharp lower
 bounds.}. In all these cases the vertices of the graph-construction
are ``coordinatized'' and two vertices are joined if some
(usually polynomial) equations are satisfied.\footnote{Some of the
 constructions may seem number theoretic.}  

Actually, this motivated Conjecture \ref{RacConj} or its weakening:
If we use a typical finite geometric construction, then there is a
$d$-dimensional space, where each vertex is joined to a
$t$-dimensional subspace. Hence $n=p^d$, the degrees are around $n^{t/d}$,
so the construction has around $n^{1+(t/d)}$ edges. The conjecture
 suggests that there are always such almost extremal
 constructions.\footnote{Here we have to make some remarks about our
  ``Conjectures'': Many of them have the feature that it is not that
  interesting if they are true or false: in proving any alternative,
  we get new, important knowledge about our topics. The first such
  ``Conjecture'' was that of Tur\'an on ``Diagonal'' Ramsey Numbers, that lead
  to the Erd\H{o}s Random Graph Approach, see Remark \ref{ProbabStarts}.}

The most important question in this part is if one can find
constructions\footnote{Or ``random constructions''.} to provide lower
bounds where the exponents match the exponents in the upper bounds.
Here we shall discuss when do we know the sharpness of the
K\H{o}v\'ari-T. S\'os-Tur\'an upper bound, $\ext(n, K_{a,b})=O(n^{2-(1/a)})$.

As we have mentioned in Section \ref{CentralExamplesS}, Koll\'ar, R\'onyai and
T. Szab\'o  \cite{KollRonySzab} gave a construction which was
improved by Alon, R\'onyai and Szab\'o  \cite{AlonRonySzab}
(Constructions \ref{co:KollRonySzabo} and \ref{co:AlonRonySzabo} below).
The basic idea of their proofs was -- at least in our interpretation
-- the same as that of William G. Brown; however, much more advanced.
In the three dimensional Euclidean space $\EE^3$ the Unit Distance Graph
 contains no $\HK33$. If we change the underlying field to a
 finite field $GF(q)$ (as Brown did in Construction~\ref{BrownConstr})
 then we obtain a finite 
 graph having $n=q^3$ vertices.
The neighborhood of each vertex will have
 $\approx q^2$ neighbors, and therefore $\approx\half n^{5/3}$ edges.
Now comes the crucial part: despite the fact, that this is highly nontrivial,
 we could say, that -- because of the geometric reason, -- this graph contains
 no $\KK33 $ proving the sharpness of \eqref{KSTshf}.\footnote{Actually, as we have
 already discussed this in Subsection \ref{BrownConstr}, the Will
 Brown's lower bound also proves this sharpness, only, the lower
 bound of  \cite{AlonRonySzab} is a little better.}

If we wish to extend the above
construction to get lower bounds for $\ext(n,K_{a,a})$ and we
 mechanically try to use unit balls in the $a$-dimensional
 space $GF(q)^a$ then several problems occur.
We would need that any $a$ of them intersect in at most $a-1$ points. Then we would be home.

In $\EE^4$ we can choose two orthogonal circles of radii $1\over \sqrt 2$,
e.g.,
$$\left\{(x_1,x_2,0,0):~
x_1^2+x_2^2=1/2\right\}
\text{and}\left\{(0,0,x_3,x_4):~
x_3^2+x_4^2=1/2\right\},$$ then
 each point on the first one has distance 1 from each point in the second one.
Hence the ``Unit Distance Graph'' contains $K(\infty,\infty)$.
 (Similarly, the ``Unit Distance Graph'' of $GF(q)^4$ contains a
 $K_{q,q}$.)
So everything seems (!) to break down?
Not quite, by the Koll\'ar-R\'onyai-Szab\'o construction.
Instead of the `Euclidean metric' they use a so-called {\emR norm}
 in the space $GF(q^a)$.
Two vectors $\xx$ and $\yy$ are connected if
 the norm of their sum is 1; $N(\xx+\yy)=1$.
(In this context there is not much difference
between connecting them this way or take a bipartite graph and
connecting the vertices in it if $N(\xx-\yy)=1$).

\begin{theorem}[Koll\'ar, R\'onyai, and T. Szab\'o  \cite{KollRonySzab} for $b> a!$, Alon, R\'onyai, and Szab\'o  \cite{AlonRonySzab} for $b> (a-1)!$]

There exists a $c_a>0$ such that for $b>(a-1)!$ we have $$\ext(n,\KK ab )>c_an^{2-(1/a)}.$$
\end{theorem}

Below we provide the Koll\'ar-R\'onyai-Szab\'o construction and a short verification.
The norm of an element $\xx\in GF(q^a)$ is defined as
$$N(\xx):=\xx\cdot \xx^q\dots\cdot \xx^{q^{a-1}}.$$

\begin{construction}[Koll\'ar-R\'onyai-T. Szab\'o  \cite{KollRonySzab}, the Norm Graph]\label{co:KollRonySzabo}

The vertices of $G(q,a)$ are the elements $\xx\in GF(q^a)$.
The elements $\xx$ and $\yy$ are joined if
$N(\xx+\yy)=1$.
\end{construction}

We claim that $G(q,a)$ is $K_{a,b}$-free where $b=a!+1$.
If we have a
$\KK ba\subseteq G(q,a)$, then fixing -- as parameters -- the $a$ vertices
$\yy_1,\dots,\yy_a$, we get $a$ equations of the form $N(\xx+\yy_i)=1$ with $b$ solutions
 $\xx\in \{ \xx_1,\dots,\xx_b\}$.
Then we can use the following result from Algebraic Geometry with $t=a$.

\begin{lemma}\label{le:KRSz}

Let $K$ be a field and $\alpha_{i,j},\beta_i\in K$ for $1\le i, j\le t$ such that
$\alpha_{i_1,j}\ne \alpha_{i_2,j}$ if $i_1\ne i_2$. Then the system of equation
\begin{eqnarray*}
(x_1-\alpha_{1,1})(x_2-\alpha_{1,2})\dots(x_t-\alpha_{1,t})&=&\beta_1\\
(x_2-\alpha_{2,1})(x_2-\alpha_{2,2})\dots(x_t-\alpha_{2,t})&=&\beta_2\\
\vdots ~\hskip2cm~ ~& &\vdots\\
(x_t-\alpha_{t,1})(x_2-\alpha_{t,2})\dots(x_t-\alpha_{t,t})&=&\beta_t
\end{eqnarray*}
has at most $t!$ solutions $(x_1,x_2,\dots,x_t)\in K^t$. \Qed
\end{lemma}

\begin{construction}[Alon-R\'onyai-T. Szab\'o  \cite{AlonRonySzab}]\label{co:AlonRonySzabo}
The vertices of the graph $H(q,a)$ are the elements
 $(x,X)\in GF(q)^*\times GF(q^{a-1})$ and $(x,X)$ and $(y,Y)$ are joined if $N(X+Y)=xy$.
\end{construction}

Here the norm $N(X)$ is defined in $GF(q^{a-1})$ and so it is $X\cdot
X^q\dots\cdot X^{q^{a-2}}$. The graph $H(q,a)$ has $(q-1)q^{a-1}$
vertices, it is $q^{a-1}-1$ regular, and contains no $K_{a,b}$ with
$b=(a-1)!+1$. To show this we use Lemma~\ref{le:KRSz} with $t=a-1$
only.

\begin{theorem}[Ball and Pepe  \cite{BallPepe}]\label{co:BallPepe}

The Alon-R\'onyai-T. Szab\'o graph $H(q,4)$ does not contain $K_{5,5}$. Hence
$\ext(n,K_{5,5})\ge (\half+o(1))n^{7/4}$.
\end{theorem}

This is better than the earlier lower bounds of $\ext(n,K_{a,b})$
for $a=5$, $5\le b\le 12$, and $a=6$, $6\le b\le 8$.

Recently, Blagojevi\'c, Bukh, and Karasev ~ \cite{BBK} gave a new
algebraic construction to provide lower bounds on $Z(m,n,a,b)$
matching the \eqref{KST_Zaran} upper bound. Their example is weaker than the
Koll\'ar-R\'onyai-Szab\'o in the sense that it only works for
$b>(a^2(a+1))^{a}$. On the other hand, they give new insight about
the limits of the Algebraic Geometric method on which constructions
may and which may not work.

We close this section mentioning that
Noga Alon has a survey paper in the Handbook of Combinatorics
 \cite{AlonHB}
providing ample information on the topics treated here
 (i.e., applications of algebra in combinatorics).

\subsection{The coefficient in the K\H{o}v\'ari-T. S\'os-Tur\'an bound}

Alon, R\'onyai and Szab\'o  \cite{AlonRonySzab} observed that their
Construction \ref{co:AlonRonySzabo} can be factored with a $t$-element
subgroup $H\subset GF(q)^*$ (when $t$ divides $q-1$) in the same way
as it was done in Construction~\ref{constK2t}. Namely, the vertex set
of the new graph $H^t(q,a)$ are the elements $(x,X)\in GF(q)^*/H\times
GF(q^{a-1})$ and $(x,X)$ and $(y,Y)$ are joined if
$N(X+Y)x^{-1}y^{-1}\in H$. Then the graph $H^t(q,a)$ has
$n=(q-1)q^{a-1}/t$ vertices, its degrees are about $q^{a-1}$, and
it contains no $K_{a,b}$ for $b=(a-1)!t^{a-1}+1$.
Let $q\to\infty$. Then also $n\to \infty$,
and we get
that for
these fixed values of $a$ and $b$
one gets
$$
  \ext(n, K_{a,b})\geq (1-o(1)){{\root a\of{b-1}} \over {2\root a \of {(a-1)!}}} n^{2-(1/a)}.
 $$
This  shows that the order of magnitude of the coefficient
 in the KST bound \eqref{KST_Zaran} should be indeed $\root a\of{b-1}$.

Mont\'agh~ \cite{Montagh} found a clever factorization of the Brown
graph (using the spherical symmetry of the balls) thus proving the
same result with even a slightly better bound than the bound of Alon,
R\'onyai and Szab\'o, for $a=3$.

\subsection{Excluding large complete subgraphs}

The following theorem was discovered many times because its
connections with Computer Science problems: given a graph $G$ with $n$
vertices, there exists a decomposition of its edges into complete
balanced bipartite graphs $K_{a_i,a_i}$ having altogether $O(n^2/\log
n)$ vertices, $\sum_i a_i =O(n^2/\log n)$. Lately, Mubayi and
Gy.~Tur\'an~ \cite{MubayiTuranGy} gave a a polynomial algorithms
finding such a subgraph partition efficiently. Strictly speaking,
this is not a Tur\'an type problem but their result implies, e.g.,
that there exists a polynomial algorithm to find a $K_{a,a}$ in a
graph of $n^2/4$ edges of size $a=\Theta(\log n)$. The bound $O(\log
n)$ is the best possible (shown by the random graph).

It is also not very difficult to show that usually the random graph gives
 the correct order of the Tur\'an number $\ext(n,K_{a,a})$ for $n,a\to \infty$ simultaneously.

The case when $a,b$ are very large i.e. $a+b=\Omega(n)$ was
considered by Griggs, Quyang, and Ho  \cite{GriggsQuyang},
 \cite{GriggsHo}. In this case $Z(m,n,a,b)$ is almost $mn$ so they
considered the dual question.
Let us mention only one result of this type by Balbuena,
Garc\'{\i}a-V\'azquez, Marcote, and Valenzuela, who have more papers on
this topic.

\begin{theorem}[see  \cite{BGMV2007},  \cite{BGMV2008} and the references there]

$Z(m,n;a,b)=mn-(m+n-a-b+1)$ if $\max\{m,n\}\le a+b-1$.
\end{theorem}

There is another direction of research, when the ratio of $m$ and $n$
is extreme. Here we only mention a classical result, that it is easy
to solve the case when $n$ is very large compared to $m$.

\begin{theorem}[\v{C}ul\'{\i}k  \cite{CulikZara}]\label{CulikTh}

$$Z(m,n,a,b)=(a-1)n+(b-1){m\choose a}\Text{for} n\ge(b-1){m\choose a}.$$
\end{theorem}

\section{Excluding Cycles : $C_{2k}$}

To start with, Bondy wrote a long chapter in the Handbook of
Combinatorics \cite{BondyHB} and also a very nice survey on Erd\H{o}s
and the cycles of graphs \cite{BondyOnErdos}.

Let $\C$ be a (finite or infinite) set of cycles. The study of
$\ext(n, \C )$ is especially interesting if $\C$ has a member of even
length. However, constructions of dense graphs without some given even
cycles is usually very difficult; the examples use polarities of
finite geometries (generalized polygons \cite{LazebUstimWoldPolar}),
or Ramanujan graphs \cite{MargulisExpander}, \cite{LuboPhilSarA} or
some other families of polynomials \cite{LazebUstimWoldarGirth}.

An odd cycle, $C_{2k+1}$ is chromatically critical. Hence
 a theorem of Simonovits  \cite{SimSymm77} implies that $\ext(n,C_{2k+1})=[\reci 4 n^2]$ for $n>n_k$
 and the only extremal graph is $K_{{\lfloor n/2\rfloor},{\lceil n/2 \rceil}}$.

In this Section we concentrate on even cycles $C_{2k}$.

\subsection{Girth and Tur\'an numbers, upper bounds}

What is $\ext(n, \{ C_3, C_4, \dots, C_{g-1}\})$, the maximum number of edges
 in a graph with $n$ vertices and girth $g$?
This problem can be considered in a dual form, what is the least number of
 vertices $n=n(d,g)$ in a graph of girth $g$ and an average degree at least $d$?
If we replace `average' with `minimum' $\delta$ then a
 simple argument gives the so-called {\emR Moore bound} for odd girth:
\beq{eq:Moore}
 |V(G)|=n\ge n_0(\delta, 2k+1):= 1+ \delta \sum_{ 0\le i \le k-1} (\delta -1)^i.
 \eeq
Alon, Hoory and Linial~ \cite{AHL} showed that \eqref{eq:Moore} holds for
 the average degree, too.
Rearranging we have $d_{\rm ave}<n^{1/k} +1$, in other words

\begin{theorem}[Upper bound when the girth is odd]
\beq{eq:girth}
 \ext(n, \{ C_3, C_4, \dots, C_{2k}\})< \half n^{1+(1/k)} + \half n.
 \eeq
\end{theorem}

To prove an upper bound $n^{1+(1/k)}$ is trivial by induction on $n$.
Then \eqref{eq:girth} was improved but with a larger linear additive term.

\begin{theorem}[Lam and Verstr\"aete \cite{LamVer}, Excluding only even cycles]
\beq{eq:evengirth}
 \ext(n, \{ C_4, C_6, \dots, C_{2k}\})< \half n^{1+(1/k)} + 2^{k^2} n.
 \eeq
\end{theorem}
They also note that for $k=2,3,5$ the $n$-vertex polarity graphs of generalized $(k+1)$-gons
 (defined by Lazebnik, Ustimenko and Woldar  \cite{LazebUstimWoldPolar} described below as Construction
 \ref{co:LazPolarity})
have $\half n^{1+(1/k)} + O(n)$ edges and have no even cycles of length at most $2k$.

\begin{corollary}[\cite{LamVer} and  \cite{LazebUstimWoldPolar}
Even girth is 6, 8 or 12]\label{co:EvenGirth}

~
In case of $2k\in \{4,6,10\}$ we have
\beq{eq:evengirth:corollary}
 \ext(n, \{ C_4, C_6, \dots, C_{2k}\})= (1+o(1))\half n^{1+(1/k)}.
 \eeq
\end{corollary}

On the other hand, F\"uredi, Naor and Versta\"ete  \cite{FureNaorVerstra}
showed that if we exclude only $C_{2k}$, then
 $\ext(n,C_6)> 0.53 n^{4/3}$
 (see below as Construction~\ref{co:randomC6})
and Lazebnik, Ustimenko, and Woldar~ \cite{LazUstWolProperties}, showed that
$\ext(n, C_{10})> 0.57 n^{6/5}$
(see below as Construction \ref{co:LazMultiplying}).

Concerning the Moore bound for even girth we have
 \beq{eq:Moore2}
 |V(G)|=n\ge n_0(\delta, 2k+2):= 2 \sum_{ 0\le i \le k} (\delta -1)^i. \eeq

Alon, Hoory and Linial~ \cite{AHL} showed that \eqref{eq:Moore2} holds
for the average degree, too.
Rearranging, we have $d_{\rm ave}<(n/2)^{1/k} +1$, in other words

\begin{theorem}[Upper bound when the girth is even]
\beq{eq:girth2} \ext(n, \{ C_3, C_4, \dots, C_{2k+1}\})< {1\over
 2^{1+(1/k)}} n^{1+(1/k)} + \half n. \eeq
\end{theorem}

 This upper bound with a weaker error term was also proved earlier by Erd\H{o}s
 and Simonovits~\cite{ErdSimComp}.

Note that because of Theorems \ref{th:Reiman}, \ref{th:Benson}, and
\ref{th:Benson12} one can easily show that asymptotic bound holds in
\eqref{eq:girth2} for $2k=4,6,10$. The other cases are unsolved.

\begin{theorem}

For $2k=4,6$ and $10$ as $n\to \infty$ we have

\beq{eq:knownGirth} \ext(n, \{ C_3, C_4, \dots, C_{2k+1}\})=(1+o(1)){1\over 2^{1+(1/k)}} n^{1+(1/k)}.
  \eeq
Moreover, infinitely many exact values are obtained for $2k=4,6,10$: for $n=2(q^k+q^{k-1}+\dots +q+1)$,
 \beq{eq:knownExactGirth}
 \ext(n,\{ C_3, C_4, \dots, C_{2k+1}\})= (q+1)(q^k+q^{k-1}+\dots +q+1) \eeq
whenever $q$ is a power of a prime.
\end{theorem}

\subsection{Excluding a single $C_{2k}$, upper bounds}
Concerning our central problem,
Erd\H{o}s showed that excluding just one even cycle has essentially the same effect as
 excluding all smaller cycles as well.
This is far from trivial! Erd\H{o}s never published a proof of his result.

\begin{theorem}[Erd\H{o}s, The Even Cycle Theorem]

\beq{eq:Erd}\ext(n,C_{2k})=O(n^{1+(1/k)}).
  \eeq
\end{theorem}

The first proof was published by Bondy and Simonovits in the following
 stronger form.

\begin{theorem}[Bondy and Simonovits  \cite{BondySim}]\label{BSthm}

 Let $\Gn$ be a graph with $e$ edges, and let $t$ satisfy
$2\le t\le e/(100n)$ and $tn^{1/t}\le e/(10n)$. Then $\Gn$ contains a
$C_{2t}$. \end{theorem}

In some sense, this is a ``pancyclic theorem'': there is a
meta-principle, that if some reasonable conditions ensure the
existence of a Hamiltonian cycle, then they ensure the existence of
all shorter cycles. Here we go the other direction: if we ensure the existence of a $C_{2k}$, then we ensure the existence of all longer cycles, up to a natural limit, with the natural parity.

\begin{corollary}\label{co:BSthm}
 If $\Gn$ has at least $100kn^{1+(1/k)}$ edges, then it
contains a $C_{2t}$, for every $t\in[k,kn^{1/k}]$.
\end{corollary}

The Erd\H{o}s-Bondy-Simonovits upper bound together with
 earlier known constructions 
imply that the exponent $1+(1/k)$ is sharp for $C_4$
 (see, e.g., Theorem~\ref{ErdRenyiSosTh}), $C_6,$ and $C_{10}$
 (Theorems \ref{co:Benson} and \ref{co:Benson12} below).

\begin{corollary}[The only known exact exponents for single cycles]\label{co:8}
 ${}$\newline
 $\ext(n,C_4)=\Theta(n^{3/2})$, \quad $\ext(n,C_6)=\Theta(n^{4/3})$, \quad $\ext(n,C_{10})=\Theta(n^{6/5})$.
\end{corollary}

The upper end of the interval in Corollary~\ref{co:BSthm} is also sharp,
 apart from the constant 100
take the disjoint union of complete graphs of order $200kn^{1/k}$.
We made the following conjecture:

\begin{conjecture}[Erd\H{o}s--Simonovits]
 $\ext(n,C_{2k})\ge c_kn^{1+(1/k)}$. Moreover,
$$\ext(n,C_{2k})\over n^{1+(1/k)}$$ converges to a positive limit.
\end{conjecture}

It is only known for $C_4$.
A weakening of this conjecture would be the following:
Let $\Theta_{k,\ell}$ denote the graph of order $2+(k-1)\ell$ in which
two vertices are joined by $\ell$ paths of length $k$.

\begin{conjecture}[Simonovits]\label{ThetaCon}

For each $k$ there is an $\ell=\ell(k)$ for which
$\ext(n,\Theta_{k,\ell})\ge c_kn^{1+(1/k)}$.
\end{conjecture}

Perhaps the first very annoying unsolved problem on this area is 

\begin{conjecture}

$\ext(n,C_8)\ge c_4n^{5/4}$.
\end{conjecture}

Returning to the Tur\'an number of $C_{2k}$, the multiplicative constant
 of the upper bound in the Bondy-Simonovits theorem was improved by
 Verstra\"ete  \cite{VertreAPCyc} from 100 to 8.
The best known upper bound today is that of Oleg Pikhurko:

\begin{theorem}[Pikhurko,  \cite{PikhurkoC2k}]
$$\ext(n,C_{2k})\le (k-1)n^{1+(1/k)} + 16(k-1)n.$$
\end{theorem}

\begin{histrem}

\dori (a) Pikhurko, in his very nice paper  \cite{PikhurkoC2k}
gives a short description of the whole story.

(b) Pikhurko mentions that the Bondy-Simonovits proof gives a
constant 20: originally it was stated as 100.
It would be extremely interesting 
 if the upper bound $k-1+o(1)$ for $\ext(n,C_{2k})/n^{1+(1/k)}$ could be improved to $o(k)$.
\end{histrem}

\subsection{Eliminating short cycles, a promising attempt}

It was relatively easy to prove the upper bound
\eqref{eq:girth} 
for the number of edges for a graph $\Gn$ with girth exceeding $2k$,
$e(\Gn)=O(n^{1+(1/k)})$.  Suppose that $G$ has no $C_{2k}$.  Erd\H{o}s
bipartite subgraph lemma \ref{ErdLemHalf} states that there is a bipartite
subgraph $H$ with $e(H)\ge \half e(G)$.  This way we have eliminated all the
odd cycles $C_3$, $C_5, \dots , C_{2k-1}$ from $G$.  It is a natural to ask
whether one can eliminate other short cycles, thus obtaining an easy proof for
the Erd\H{o}s-Bondy-Simonovits upper bound, \eqref{eq:Erd}.

\begin{problem}\label{prob:girth}

Is it true that there exists a constant $\alpha_{2k}>0$ such that each $C_{2k}$-free
$\Gn$ contains an $\Hn$ with $\girth(\Hn)>2k$ and $e(\Hn)>\alpha_{2k} e(\Gn)?$
\end{problem}

The answer is still unknown.
The first step was done by E. Gy\H ori.
The following lemma implies that $\alpha_6$ exists and it is at least $1/4$.

\begin{lemma}[Gy\H{o}ri  \cite{GyoriC6}]
  If $\Gn$ is bipartite and it does not contain any $C_6$, then it
 contains an $\Hn$ with $$e(\Hn)\ge\half e(\Gn)+1, $$
not containing $C_4$'s either (for $e(G)\geq 2$). This is sharp only for
$\Gn=\KK 2{n-2}$.
\end{lemma}

We mention two generalizations.

\begin{theorem}[F\"uredi, Naor and Verstraete  \cite{FureNaorVerstra}]
\label{delete}

Let $G$ be a hexagon-free graph. Then there exists a subgraph of
$G$ of girth at least five, containing at least half the edges
of $G$.
\end{theorem}

Furthermore, equality holds if and only if $G$ is
a union of edge-disjoint complete graphs of order four or five.
We got $\alpha_6=1/2$.

\begin{theorem}[Getting rid of $C_4$'s, K\"uhn and Osthus  \cite{KuehnOstGyori}]

Every bipartite $C_{2k}$-free
graph $G$
contains a $C_4$-free subgraph
$H$
with $e(H)\ge e(G)/(k-1)$.
\end{theorem}

The factor $1/(k-1)$ is best possible, as the example $K_{k-1, n-k+1}$ shows.

These theorems settle some special cases
(namely $\LL=\{ C_4, C_{2k}\}$) of the following compactness conjecture of
Erd\H{o}s and Simonovits.

\begin{conjecture}[Compactness. Erd\H{o}s-Simonovits  \cite{ErdSimComp}]\label{conj:Compactness}

For every finite family of graphs $\LL$ (containing bipartite members as well)
there exists an $L_0\in\LL$ for which $\ext(n,\LL) = O (\ext (n, L_0))$.
\end{conjecture}

The following result of K\"uhn and Osthus makes a little step toward solving
Problem~\ref{prob:girth} and Conjecture~\ref{conj:Compactness}.

\begin{theorem}[\cite{KuehnOstGyori}]

Let $g\geq 4$ be an even integer and let $\ell(g)=: \Pi_{1\leq i\leq g/2}\, i$.
Suppose that $k-1$ is divisible by $\ell(g)$ and $\Gn$ is a $C_{2k}$-free graph.
Then $\Gn$ contains an $\Hn$ with $\girth(\Hn)>g$ such that
 $e(\Hn)\ge e(\Gn)/2(4k)^{(g-2)/2}$.
\end{theorem}

In other words, for some very special values of $k$'s a $C_{2k}$-free
graph contains a subgraph having a positive fraction of the edges and
of girth at least $\Omega(\log k / \log \log k)$.

\subsection{A lower bound for $C_6$ : The Benson Construction}\label{BensonS}

In the preceding section, we asserted that the Erd\H{o}s theorem on even
circuits is sharp for $C_4$, $C_6$ and $C_{10}$ (and is conjectured to
be sharp in all cases). For $C_4$, the sharpness follows from
Construction \ref{PolarAff}. For $C_6,$ it can be deduced from the
Benson construction  \cite{Benson} which we explain below.
Note that (about the same time) Singleton  \cite{Singleton} described the same graph
but his definition was much more complicated.

The points of the $d$-dimensional finite projective geometry $PG(d,q)$
 are the equivalence classes of the nonzero vectors of $GF(q)^{d+1}$
 where $\xx$ and $\yy$ are equivalent if there is a $\gamma\in GF(q)^*$
 such that $\xx = \gamma \yy$.
There are $(q^{d+1}-1)/(q-1)$ such classes.
Then the $i$-dimensional subplanes are generated by the $(i+1)$-dimensional
 subspaces of the vector space $GF(q)^{d+1}$.

Let $$A=\begin{pmatrix}0&0&0&0&1\cr
   0&0&0&1&0\\
   0&0&1&0&0\\
   0&1&0&0&0\\
   1&0&0&0&0\\
\end{pmatrix}.$$
Clearly, $A$ is non-singular.
Define the surface $S$ by the equation ${\bf x}A{\bf x}^T=0$,
 $$
  S:=\{ \xx \in PG(4,q): {\bf x}A{\bf x}^T=0 \}$$

\begin{construction}[Benson's $C_6$-free bipartite graph]
\label{Benson6}

Let $\LL$ be the set of lines of $PG(4,q)$ contained entirely in $S$.
The vertex set of the bipartite graph $B_q$
is $S\cup\LL$,
and $\xx\in S$ is joined to $L\in \LL$ if $\xx\in L$.
\end{construction}

\begin{theorem}[Benson  \cite{Benson}, Singleton  \cite{Singleton}]
\label{th:Benson}

$B_q$ is a $(q+1)$-regular, bipartite, girth 8 graph with $2(q^3+q^2+q+1)$ vertices.
\end{theorem}

\begin{corollary}
\label{co:Benson}
$\ext(n, C_6)\geq (1+o(1)) (n/2)^{4/3}$.
\end{corollary}

First, we can see that $S$ does not contain a full 2-dimensional projective plane.
We can use the fact that for ${\bf x}$ and
${\bf y}$ on
$S$, the line ${\bf xy}$ consists of the points ${\bf z}=a{\bf x}+(1-a){\bf
y}$, and lies entirely in
$S$ if both ${\bf y}A{\bf y}^T=0$ and ${\bf x}A{\bf y}^T=0$.

Second, the number of lines from $\LL$ containing a given point ${\bf x}\in S$ is $q+1$.
Since the number of points on a line is $q+1$ we immediately get that $|S|=|\LL|$.

 Furthermore, $B_q$ contains no cycles of length 3, 4, 5 or 7.
(For the odd cases this is because it is bipartite, and the existence of a 4-cycle
would imply that two points of $S$ are on two distinct lines.)
Now suppose that $B_q$ contains a
6-cycle $v_1w_1v_2w_2v_3w_3v_1$.
Then $S$ must contain the three lines $v_1v_2,$
$v_2v_3,$ and $v_3v_1$, and so it must contain the plane $<v_1v_2v_3>$.
But this is impossible.
If we apply a coordinate transformation $T$ with $v_1$,
$v_2$ and $v_3$ as the first three base vectors, we get the matrix
$$\begin{pmatrix}0&0&0&?&?\\
   0&0&0&?&?\\
   0&0&0&?&?\\
   ?&?&?&?&?\\
   ?&?&?&?&?\\
\end{pmatrix}$$
since $v_iAv_j^T=0.$ But then $A'$ cannot be regular, contradicting the
regularity of $A$. Hence $B_q$ cannot contain $C_6$ either.

All these imply that $|S|=q^3+q^2+q+1$ and that every $\xx\in S$, $L\in \LL$
 if $\xx\not\in L$ then there exists a unique line $L'\in \LL$
 such that $\xx\in L'$ and $L\cap L'\neq\emptyset$.

In concluding this section, we note that finite geometry constructions can
also be used in hypergraph extremal problems (see  \cite{BrownErdSosA},  \cite{BrownErdSosB} and  \cite{SimColumbus}).

\subsection{Girth 12 graphs by Benson and by Wenger}\label{SWengerC12}

\begin{theorem}[Benson  \cite{Benson}]
\label{th:Benson12}

Let $q$ be an odd prime power.
There is a $(q+1)$-regular, bipartite, girth 12 graph $B^*_q$ with $2(q^5+q^4+q^3+q^2+q+1)$ vertices.
\end{theorem}

\begin{corollary}
\label{co:Benson12}

$\ext(n, C_{10})\geq (1+o(1)) (n/2)^{6/5}$.
\end{corollary}

One half of the vertex set of $B^*_q$ are the points of the quadric
$Q_6$ in $PG(6,q)$ defined by $x_0^2+x_1x_{-1}+ x_2x_{-2}+
x_3x_{-3}=0$. Its size is exactly $(q^6-1)/(q-1)$. Then we select a
set of lines $\LL$ contained entirely in $Q_6$ and covering each point
of $Q_6$ exactly $q+1$ times. The family $\LL$ is selected as
follows: If $\xx\in Q_6$ and $\xx, \yy\in L\in \LL$ then $\xx$ and
$\yy$ must satisfy the following six bilinear equations:
$$x_0y_i-x_iy_0+ x_{-j}y_{-k}- x_{-k}y_{-j}=0$$ where $(i,j,k)$ is a
cyclic permutation of $(1,2,3)$ or $(-1,-2,-3)$.

\begin{construction}\label{Benson12}
The bipartite graph $B_q^*$ is defined, as before, by the incidences $\xx\in L$.
\end{construction}

Now consider the much simpler example of Wenger.

\begin{construction}[Wenger  \cite{Wenger}]\label{co:Wenger}

Let $p$ be a prime, $k=2,3$ or $5$.
$H_k(p)$ is defined as a bipartite graph with two vertex classes $\bA$ and $\bB$, where
$|\bA|=|\bB|=p^k$ and
 the vertices of $\bA$ are $k$-tuples $\aa=(a_0,a_1,\dots,a_{k-1})\in GF(p)^k$
 and same for $\bb=(b_0,b_1,\dots,b_{k-1})\in \bB$.
 The vertices $\aa$ and $\bb$ are joined if
$$b_j\equiv a_j+a_{j+1}\cdot b_{k-1}~(\mod p)\text{for}j=0,1,\dots,k-2.$$
\end{construction}

One can see that for every $\aa\in\bA$ each $b_{k-1}$ determines
exactly one $\bb\in\bB$ joined to it. This easily implies that
$G[\bA,\bB]$ is $p$-regular, with $n=2p^k$ vertices and
$p^{k+1}=(n/2)^{1+(1/k)}$ edges.

Wenger gives an elegant proof of that $H_2(p)$ has no $C_4$, $H_3(p)$ has no
$C_4$, nor $C_6$.
Finally, $H_5(p)$ contains no $C_4,C_6$ or $C_{10}$, however, it has many $C_8$'s.

\subsection{Short cycles, $C_6$ and $C_{10}$}\label{Sshort}

The densest constructions of $2k$-cycle-free graphs for certain small values of $k$ arise from the existence of rank two geometries called {\emR generalized $d$-gons}.
These may be defined as
rank two geometries whose bipartite incidence graphs are
regular graphs of diameter $d$ and girth $2d$.
These are known to exist when $d$ is three, four or six.
This is the background of the above Constructions~\ref{Benson6} and \ref{Benson12}.

\begin{construction}[Lazebnik, Ustimenko and Woldar~ \cite{LazebUstimWoldPolar}]\label{co:LazPolarity}

One can use the existence of polarities of the generalized $(k+1)$-gons
to obtain dense $2k$-cycle-free graphs when $k \in \{2,3,5\}$.
In particular, for these $k$'s
\beq{eq:polarityval:1per2} \ext(n,C_{2k}) \; \ge \; \frac{1}{2}n^{1+(1/k)} + O(n)
 \eeq
for infinitely many $n$.
 \end{construction}

In  \cite{ErdSimComp}, Erd\H{o}s and Simonovits formulated the following
conjecture. For fixed $k$ and $n\to\infty$,
$\ext(n,C_{2k})=\frac{1}{2}n^{1 + (1/k)}+o(n^{1+(1/k)})$.
This holds for $C_4$ (Theorem \ref{ErdRenyiSosTh}), but was disproved
first for $C_{10}$, then for $C_6$ by the following two examples.

\begin{construction}[Lazebnik, Ustimenko and Woldar~ \cite{LazebUstimWoldPolar}]\label{co:LazMultiplying}

Consider a bipartite graph $G[A,B]$ of girth exceeding $2k$.
Replace each vertex of $A$ by $k-1$ new vertices with the same neighborhood.
Then the new graph \mbox{$G[(k-1)A,B]$} is still $C_{2k}$-free.
In particular, starting with the girth 12 bipartite graph of
Theorem~\ref{th:Benson12} (here $k=5$)
 one gets a graph of about $5q^5$ vertices and about $4q^6$ edges, implying
\beq{eq:C10bestLower} \ext(n,C_{10}) \; \ge \; 4 (n/5)^{6/5} \; > \; 0.5798 n^{6/5}
 \eeq
for infinitely many $n$.
 \end{construction}

Since the $C_6$-free graph of Construction \ref{co:LazPolarity}
 does not have $C_3$ and $C_4$ either, doubling a random subset appropriately,
 one obtains a denser $C_6$-free graph:

\begin{theorem}[F\"uredi, Naor and Verstr\"aete  \cite{FureNaorVerstra}]\label{co:randomC6} 

For infinitely many $n$,
$$ \ext(n,C_6) \; \; > \; \; \frac{3(\sqrt{5}-2)}{(\sqrt{5}-1)^{4/3}} n^{4/3} + O(n) > 0.5338 n^{4/3}.$$
\end{theorem}

They also showed that
\begin{theorem}[F\"uredi, Naor and Verstr\"aete  \cite{FureNaorVerstra}]\label{C6_felso}

$\ext(n,C_6) \le \lambda n^{4/3} + O(n)$, where
$\lambda \approx 0.6271$ is the real root of $16\lambda^3-4\lambda^2+\lambda-3=0$.
\end{theorem}

These theorems give the best known lower and upper bounds for $\ext(n,C_6)$.
The proof of Theorem \ref{C6_felso} requires a statement about
hexagon-free bipartite graphs, which is interesting in its
own right (see de Caen and Sz\'{e}kely~ \cite{CaenSzek}).
Let $\ext(m,n,C_6)$ be the maximum number
of edges amongst all $m$ by $n$ bipartite hexagon-free graphs.
Then

\begin{theorem}[F\"uredi, Naor and Verstra\"ete  \cite{FureNaorVerstra}]
\label{bipC6}

Let $m,n$ be positive integers with $n \ge m$. Then
$$ \ext(m,n,C_6) < 2^{1/3}(mn)^{2/3} + 10n.$$
Furthermore, if $n = 2m$ then as $n$ tends to infinity,
$$ \ext(m,n,C_6) = \left\{\begin{array}{ll}
2^{1/3}(mn)^{2/3} + O(n) & \mbox{ for infinitely many }m\\
2^{1/3}(mn)^{2/3} - o(n^{4/3}) & \mbox{ for all }m.
\end{array}\right.
$$
\end{theorem}

The lower bound is given by the graph defined in
Construction~\ref{co:LazMultiplying} starting with the Benson graph
(Theorem~\ref{th:Benson}, $k=3$).

\subsection{Bipartite hosts with extreme sides}\label{SextremSides}

We have already seen two such results concerning the Zarankiewicz number, by
Reiman (Theorem \ref{th:Reiman}) and \v{C}ul\'{\i}k (Theorem \ref{CulikTh}).
Andr\'as S\'ark\"ozy and Vera S\'os formulated the following
conjecture
\footnote{A weaker version of this conjecture was formulated by Erd\H{o}s several years earlier.}

\begin{conjecture}\label{SarkoSosConj}
$$\ext(m,n,C_6)<2n+c(nm)^{2/3}.$$
\end{conjecture}

A weaker version of this was proved by G\'abor N. S\'ark\"ozy,
 \cite{SarkoGNC6} and later Gy\H{o}ri  \cite{GyoriC6} proved a stronger

\begin{theorem}\label{th:Gyori}

There exists a constant $c_k>0$ for which if $G[A,B]$ is a
bipartite graph with color classes $A,B$, and $|A|=m, |B|=n\ge m^2,$ and
$$e(G[A,B])\ge (k-1)n + c_km^2,$$
then $G[A,B]\supset C_{2k}$.
 \end{theorem}

This means that for $n>m^2$ the extremal number becomes linear.
For more recent results see, e.g., Balbuena, Garc\'{\i}a-V\'azquez,
Marcote, and Valenzuela  \cite{BGMV2007Gyori}.
Later Gy\H{o}ri
\cite{Gyori2006CPC}  showed that $c_3=1/8$, proving
$$\ext(m,n, C_6)\leq 2n + {1\over 8}m^2,$$ for $n,m > 100$, $n\ge m^2/16$
 and here equality holds if $m$ is a multiple of $4$.

\subsection{The effect of odd cycles}\label{C4C5}

Let $\LL$ be a set of graphs and let $\extbip(n,\LL)$ denote
 the {\emR bipartite Tur\'an number} of $\LL$, the size of the largest
 $\LL$-free bipartite graph on $n$ vertices.

\begin{theorem}
[Erd\H{o}s and Simonovits~ \cite{ErdSimComp}]
\label{th:C4C5}

$$ \ext(n, \{C_4, C_5\})=(1+o(1)) \extbip(n, C_4)=(1+o(1))(n/2)^{3/2}.$$
 \end{theorem}

They also {\bf conjecture} that the same holds for
 $\{ C_3, C_4\}$ (i.e., for the girth problem) but this is still unsolved.
Then, they make the following bold conjecture.

\begin{conjecture}[Erd\H{o}s and Simonovits~ \cite{ErdSimComp}
 on the effect of odd cycles]

${}$\newline
Let ${\mathcal C}^{\rm odd}_{2\ell+1}$ denote the set of odd cycles $\{ C_3, C_5, \dots, C_{2\ell+1} \}$.
For any family $\LL$ consisting of bipartite graphs
 there exists an odd integer $2\ell+1$ such that
 $\ext(n, \LL\cup {\mathcal C}^{\rm odd}_{2\ell+1})\approx \extbip(n,\LL)$.
\end{conjecture}

This conjecture was verified in a few cases by extending and sharpening
 Theorem~\ref{th:C4C5} as follows.

\begin{theorem}[Keevash, Sudakov and Verstra\"ete~ \cite{KeevSudVersES}]\label{th:KSV_ErdosSimC5Conj}
 Let ${\mathcal C}^{\rm even}_{2k}$ denote the set of even
cycles $\{ C_4, C_6, \dots, C_{2k} \}$. Suppose that $2k\in
\{4,6,10\}$ and suppose that $2\ell+1>2k$. Then
$$ \ext(n, {\mathcal C}^{\rm even}_{2k}, C_{2\ell+1})=(1+o(1)) \extbip(n, {\mathcal C}^{\rm even}_{2k})
  \sim (n/2)^{1+(1/k)}.$$
 \end{theorem}

They even proved a stability result (when $n\to \infty$) and, using
it, an exact version: 
If $2k\in \{4, 6,10\}$ and $2\ell+1\geq 5$, $15$, or $23$, respectively,
and $n=2(q^{k}+ q^{k-1}+ \dots +q+1)$ then for $n> n_{2\ell+1}$
we have $$ \ext(n, {\mathcal C}^{\rm even}_{2k}\cup C_{2\ell+1})\leq (q+1)n$$
 and here equality holds only if there is a generalized $(k+1)$-gon of order $q$.

In a more recent work Allen, Keevash, Sudakov and Verstra\"ete
 \cite{AllenKeevSudVers}
verified the stronger form of the Erd\H{o}s-Simonovits conjecture proving
 that for any fixed $2\ell+1\geq 5$ one has
 $\ext(n, \{ K_{2,t}, C_{2\ell+1}\})\sim \extbip(n,K_{2,t})$
 and $\ext(n, \{K_{3,3}, C_{2\ell+1}\})\sim \extbip(n,K_{3,3})$.
They also show $$
 \ext(n, \{ K_{2,t}, B_t, C_{2\ell+1}\})\sim \extbip(n,\{ K_{2,t}, B_t\})\sim (n/2)^{3/2}
 $$ for any fixed $t\geq 2$ and $2\ell+1\geq 9$, where $B_t$ is a
 ``book'' of $t$ $C_4$'s sharing and edge: it has $2t+2$
 vertices and $3t+1$ edges. Their main tool is the smoothness of the
 corresponding Tur\'an number's and the sparse regularity lemma of
 A. Scott~ \cite{ScottRegu}.

On the other hand, for any $t\geq 1$ and prime $q> 2^{t^4}$, they construct $(t + 2)$-partite graphs
 $G_{q,t}$ with no triangle or $K_{2,2t+1}$ having $n=(t+2)q^2$ vertices and
 ${t+2\choose 2}q^2(q-1)$ edges. This implies
\beq{KeevA}\ext(n, \{ K_{2,2t+1}, C_3\})\geq (1+o(1)){t+1\over\sqrt{t+2}}n^{3/2}.\eeq
 So, using $\extbip(n, K_{2t+1})\sim \sqrt{t}n^{3/2}$, which follows easily from
 \eqref{HylCav} and \eqref{MorsSymm}, they obtain
 \beq{KeevB}
\liminf_{n\to \infty} {\ext(n, \{ K_{2,2t+1}, C_3\}) \over
 \extbip(n, K_{2,2t+1})}\geq {t+1\over\sqrt{t(t+2})} >1.
 \eeq
In particular the ratio is $2/\sqrt{3} +o(1)$ for $K_{2,3}$.
We explain their
construction yielding \eqref{KeevA} only for $t=1$.

\begin{construction}[Allen, Keevash, Sudakov and Verstra\"ete  \cite{AllenKeevSudVers}]

Let $q\equiv 2 \pmod 3$ be a prime.
 Let $G^q$ be a
three-partite graph with parts $A_1$, $A_2$ and $A_3$ which are copies of $GF(q)\times GF(q)$.
 Join $(x_1,x_2)\in A_i$ to $(y_1,y_2)\in A_{i+1}$ if
 $$
  (y_1,y_2)=(x_1,x_2)+(a,a^2)
  $$
  for some $a\in GF(q)$, $a\neq 0$.
 \end{construction}

The obtained graph is $K_{2,3}$ and $C_3$-free, and has $n=3q^2$ vertices and almost $n^{3/2}/\sqrt 3$ edges.
This yields the ratio $2/\sqrt{3}+o(1)$ for $K_{2,3}$  in \eqref{KeevB}.
They believe that Erd\H{o}s' Conjecture \ref{ErdosC3C4} is false:

\begin{conjecture}[\cite{AllenKeevSudVers}]
 $$
\liminf_{n\to \infty} {\ext(n, \{ C_3, C_4\}) \over
 \extbip(n, C_4)} >1.
 $$
\end{conjecture}

\subsection{Large girth: Ramanujan graphs}\label{LuboPhilSarnakS}

Until this point we were fixing the excluded subgraphs. However, there
is a subcategory of extremal graph problems, which we could also call
``Parametrized Extremal Graph Problems''. Instead of defining them we
give an almost trivial but important example:
Horst Sachs and Erd\H{o}s  \cite{ErdSachs} reformulated the Moore bounds
 \eqref{eq:Moore} \eqref{eq:Moore2}, in a slightly different form.

\begin{theorem}

If the minimum degree of $\Gn$, $d:=\mindeg(\Gn)>2$ then $\Gn$
contains a $C_\ell$ with
\beq{TutteErd}\ell<{2\log n\over \log (d-1)}.\eeq
\end{theorem}


Here we arrived at an area where some constructions (for lower bounds)
were needed, and the lower bounds were easily obtained by
probabilistic arguments; however they were very difficult to obtain
them in a constructive way.  Instead of going into details, we mention
a result of Margulis \cite{MargulisCCA} that
  \eqref{TutteErd} is sharp up to a constant: there are -- not too
  complicated -- Cayley graphs of constant (even) degrees $d$ and
  girth at least $c\log_{d-1} n$. Here -- surprisingly, Margulis'
  construction is better than the random graph and a construction of
  Imrich yields an even better constant $c$:

\begin{theorem}[Imrich \cite{ImrichGirth}]
For every integer $d>2$ one can (effectively) construct
infinitely many $d$-regular Cayley graphs $X_n$ with $$\girth(X_n)>0.4801{\log n
  \over \log(d-1)}-2.$$
\end{theorem}

The next step in this area was a much
deeper and more important results of Margulis
\cite{MargulisExpander,MarguTaskent}, Lubotzky, Phillips and Sarnak,
\cite{LuboPhilSarA} on the Expander graphs, that are
eigenvalue-extremal. In this sense the
Margulis-Lubotzky-Phillips-Sarnak graph is very nice.  There is only
one problem with it. While defining these graphs is non-trivial, but
not extremely complicated, to verify their extremal properties
requires deep mathematical tools. Below we give a very compressed
description of it.

\begin{definition}
   Given a connected $k$-regular graph $X$, we denote by $\lambda(X)$
  the largest of the absolute values of eigenvalues of the adjacency matrix of
  $X$, different from $k$.  An $n$-vertex $k$-regular graph $X_{n,k}$ is a
  {\emR Ramanujan graph} if $\lambda(X_{n,k})\le 2\sqrt{k-1}.$
\end{definition}

\begin{remark}

 In case of $k$-regular graphs, the largest absolute values of the
 eigenvalues is $k$. The bipartite graphs have the property that if
 $\lambda_i$ is eigenvalue, then $-\lambda_i$ is also an eigenvalue.
By the Alon-Boppana inequality, (see Proposition 4.2 of  \cite{LuboPhilSarA})
$$\liminf_{n\to\infty} \lambda(X_{n,k})= 2\sqrt{k-1}.$$
\end{remark}

Ramanujan graphs are important because they are expander graphs,
which are extremely important in Theoretical Computer Science.

There are quite a few cases, where -- instead of using ``random graph
constructions'' one tries to use Cayley Graphs.
{\emR Cayley graphs} are graphs whose vertices are the elements of some group
$\GG$ and the edges are the pairs $(g,\alpha_ig)$, where $g\in\GG$ and
$\alpha_1,\dots,\alpha_k$   are some elements of 
                                            $\GG$. If
we look for a digraph, then this is a correct definition.
However, if
we are looking for ordinary graphs, then we have to assume that
$S:=\{\alpha_1,\dots,\alpha_k\}$ is closed under taking the inverse:
if $\alpha\in S$ then $\alpha^{-1}\in S$ as well. If we
choose $\GG$ and $S$ appropriately, then the obtained graph will
provide us with nice constructions; mainly, because it behaves as if
it were a random graph, or, occasionally, even better.

\begin{construction}[\cite{LuboPhilSarA}]\label{LuboPhilSarnak}
Let $p$ and $q$ be unequal primes congruent to 1 mod 4. The Ramanujan
graphs $X^{p,q}$
of  \cite{LuboPhilSarA}
are $p+1$-regular Cayley graphs\footnote{There are two of them, a bipartite and a
 non-bipartite, we forget the bipartite one.}
 of the group ${\bf PSL}(2,\ZZ/q\ZZ)$: $p+1$ generators
of the group are fixed, which are obtained from the solutions of
\beq{Ramanu}p=a^2+b^2+c^2+d^2, \text{where}a>0 \text{is odd and} b, c, d\text{ are even.}\eeq
 \end{construction}

The number of solutions of \eqref{Ramanu} is connected to the famous
Ramanujan conjecture, which is still open.
However, good approximations are known, by Eichler and Igusa, enough for the purposes of  \cite{LuboPhilSarA}.
Originally most of the authors were interested in the eigenvalue properties
(spectral gap) of these graphs, that are also strongly connected to them being expander graphs
(see Alon,  \cite{AlonExpand}, Alon--Milman  \cite{AlonMilman}).

From here on, $X_{n,k}=X^{p,q}$ is a special sequence of Ramanujan graphs,
which is non-bipartite if the Jacobi symbol $({q\over p})=1$; then it has
$n=(q^3-q)/2$ vertices.

\begin{theorem}
 For $k=p+1$, $X^{p,q}$ is $k$-regular, its eigenvalues are
$\lambda=\pm k$ or $|\lambda |\le 2 \sqrt{k-1}$.
\end{theorem}

This property is optimal and leads to the best known explicit expander graphs.
Alon turned the attention of the authors to that
 these graphs satisfy a number of extremal combinatorial properties.

\begin{theorem}[Observation of Alon]

The girth of $X_{n,k}$ is asymptotically $\ge {4\over3} {\log n\over \log (k-1)}$.
\end{theorem}

This gives larger girth than what was previously known by explicit or
non-explicit constructions. Also, it is one of the ``cleanest'' way to define graphs with large girth and high chromatic number:

\begin{theorem}[\cite{LuboPhilSarA}]

If $X_{n,k}$ is a non-bipartite Ramanujan graph, then
its independence number and chromatic number satisfy
$$\alpha(X_{n,k})\le {2\sqrt{k-1}\over k}n\Text{and}
\chi(X_{n,k})\ge {k\over 2\sqrt{k-1}}.$$
\end{theorem}

For a more informative description of these and many other related areas see
the survey of Alon in the Handbook  \cite{AlonHB}.

\subsection{The girth problem: the Lazebnik-Ustimenko approach}

After 20 years  Theorem \ref{th:bestGirth} still yields
  the best known lower bound for the girth problem:
Lazebnik, Ustimenko and Woldar's work  \cite{LazebUstimWoldarGirth}
gives a slight improvement (an $O(1)$ in the denominator of the exponent)
to what we can get from the Ramanujan' graphs.

\begin{theorem}[
  \cite{LazebUstimWoldarGirth}]\label{th:bestGirth}
 $\ext(n, \{ C_3, C_4, \dots, C_{2k+1}\})=\Omega(n\cdot n^{2/(3k-3+\varepsilon)})$
 where $k\ge 2$ is fixed, $\varepsilon=0$ if $k$ is odd, $\varepsilon=1$ if $k$ is even
 and $n\to \infty$.
 \end{theorem}

We have seen basically two approaches on how to {\emR construct}
graphs with high girth. One was the use of Finite Geometries, and the
other the use of Cayley Graphs of some matrix groups (Ramanujan graphs).
There is (at least) one further important approach to this question
which we find in the works of Lazebnik and Ustimenko and later Lazebnik,
Ustimenko and Woldar.

\begin{remark}[History]  
 In this
 survey many important areas had to be skipped.
One of them is the family of Lazebnik-Ustimenko type algebraic
constructions. This family of constructions is much more flexible
than many earlier ones, and provides a lot of new constructions in
extremal graph theory, in Ramsey type problems, for graphs and
hypergraphs as well. The first results were achieved by Lazebnik and
Ustimenko  \cite{LazUstiExa}. Lazebnik and his
coworkers created a school in this area. The reader is referred here
to  \cite{LazebMubayi}. 
\end{remark}

The main feature of this approach can be described (perhaps slightly cheating) as follows.
We take a set $R$ (a finite or infinite ring or field), its $d\th$ power,
 and a sequence of polynomials $f_2, \dots, f_d$.
Define a bipartite graph,
where the colour classes $A$ and $B$ consist of vectors
$(a_1,\dots,a_d)$ and $(b_1,\dots,b_d)$ that are joined if
\begin{eqnarray*}
a_2+b_2&=&f_2(a_1,b_1)\\
a_3+b_3&=&f_3(a_1,b_1,a_2,b_2)\\
\dots\\
a_d+b_d&=&f_d(a_1,b_1,\dots,a_d,b_d).
\end{eqnarray*}
We may also identify $A$ and $B$ to get non-bipartite graphs as well.
In general,
 either we get digraphs, or
some symmetry conditions are assumed on the functions $f_i$, ensuring that if
$(a_1,\dots,a_d)$ is joined to $(b_1,\dots,b_d)$, then
$(b_1,\dots,b_d)$
and
$(a_1,\dots,a_d)$
are joined as well.
Yet, it is not an easy area to describe it on a few pages: this is why we basically skip it.
Perhaps the more interested reader should look at  \cite{LazebWoldarEqu}.

\subsection{Cycle length distribution}

As a measure of the density of the cycle lengths in a graph $G$, Erd\H{o}s
introduced the number $L(G)$, the sum of the reciprocals of the distinct cycle
lengths of $G$.  The following beautiful theorem, due to Gy\'arf\'as, Koml\'os
and Szemer\'edi, proves a conjecture of Erd\H{o}s and Hajnal, asserting that
in some sense the complete graph or the complete bipartite graph are the
densest concerning cycle lengths:

\begin{theorem}[\cite{GyarfasKomSzem}]

There exists a positive constant $c>0$ such that if $\mindeg(G)\ge k$,
 then for the sum of the reciprocals of the cycle lengths $\ell_i$ of $G$ we have
$$L(G)=\sum {1\over \ell_i}>c\log k.$$
\end{theorem}

The union of complete graphs $K_{k+1}$ or bipartite graphs $K_{k,m}$
(where $m \ge k$) show that this lower bound is sharp.

Generalizing a theorem of Bondy and Vince  \cite{BondyVince},
Gengua Fan proved several nice results on the distribution of cycle lengths.
We mention only one of them.

\begin{theorem}[G. Fan  \cite{FanGDistri}]
 Let $xy$ be an edge in a 2-connected graph $G$, $k$ be a positive
integer and suppose that all the vertices of $G$ but $x$ and $y$ have degrees
 at least $3k$.
Then $xy$ is contained in $k+1$ cycles $C^0,C^1,\dots,C^k$, such that
$k+1<|E(C^0)|<|E(C^1)|<\dots<|E(C^k)|$, $|E(C^i)|-|E(C^{i-1})|=2$ for
 $i=1, \dots, k-1$ and $1\le |E(C^k)|-|E(C^{k-1})|\le 2$.
\end{theorem}

A related result concerning $k$ odd cycle lengths can be found in Gy\'arf\'as~ \cite{GyarfasK}.

Next we recall a conjecture of Burr and Erd\H{o}s.

\begin{conjecture}[Burr and Erd\H{o}s]\label{BurrErdTop}

For every odd integer $k>0$, and every integer $\ell$, there exists a $c_k$ such that if
$e(\Gn)>c_kn$, then some $m\equiv \ell \pmod k$, we have
$C_m\subseteq\Gn$.
\end{conjecture}

This was proved by Bollob\'as  \cite{BolloCycleMod}
 with $c_k\leq ((k+1)^k-1)/k$.
H\"aggkvist and Scott (\cite{HaggkScottArith},  \cite{HaggkvScottCyc}) decreased $c_k$ and
 extended the Bollob\'as result, proving that every graph $\Gn$
 with minimum degree at least $300k^2$ contains $k$ cycles of consecutive even lengths.
Soon after, the right order of magnitude of $c_k$ was established.

 \begin{theorem}[Verstra\"ete  \cite{VertreAPCyc}]

Let $\Gn$ be a graph with $e(\Gn)\geq 4kn$.
Then there are cycles of $k$ consecutive even lengths in $\Gn$.
\end{theorem}

We close this part with the following theorem:

 \begin{theorem}[Sudakov, Verstra\"ete  \cite{SudakVerstra}]
 Let $\girth(\Gn)=g$ be fixed and $d=2e(\Gn)/n$. Let $\C(G)$
denote the set of cycle-lengths in $G$. Then $\C(\Gn)$ contains at
least $\Omega(d^{\lfloor(g-1)/2\rfloor)})$ consecutive even integers,
as $d\to\infty$.
\end{theorem}

\section{Paths and long cycles}

In this section we shall describe
results connected with $\ext(n,\Pk)$, $\ext(n,\Cg k)$, (where the
cycles of at least $k$ vertices are excluded).
This problem was proposed by Tur\'an
and the (asymptotic) answer were given by Erd\H{o}s--Gallai.

\subsection{Excluding long cycles}

\begin{theorem}[Erd\H{o}s and Gallai  \cite{ErdGallai}]\label{ErdGallaiCyc}

Let $\Gn$ be a graph with more than $\half(k-1)(n-1)$ edges, $k\ge 3$.
Then $\Gn$ contains a cycle of length at least $k$.
This bound is the best possible if $n-1$ is divisible by $k-2$.
\end{theorem}

A matching lower bound $\half(k-1)n-O(k^2)$ can be obtained
gluing together complete graphs of sizes at most $k-1$.
If $k$ is odd, then there are nearly extremal graphs having a
 completely different structure.
Namely, one can take a complete bipartite
 graph with partite sets $A$ and $B$ of sizes $|A|={{k-1}\over 2}$ and
 $|B|=n- {{k-1}\over 2}$ and add all edges in $A$, too.

The exact value was determined by Woodall  \cite{WoodallA}
 and independently and at the same time by Kopylov  \cite{Kopylov}.

\begin{theorem} [\cite{Kopylov},  \cite{WoodallA}]
\label{th:Kopylov}

Let $n=m(k-2)+r$, where $1\le r\le k-2$, $k\geq 3$, $m\ge 1$ integers.
If
$$e(\Gn) > m{k-1\choose 2} + {r\choose 2},$$
then $\Gn$ contains a cycle of length at least $k$,
and this bound is the best possible:
\beq{eq:Kop_cycle}
 \ext(n , \Cg k)= \half (k-1)n - \half r(k-r).
 \eeq
 \end{theorem}

Caccetta and Vijayan~ \cite{CV1991} gave an alternative proof of the result.
We need a definition.

\begin{construction}

Let $H_{n,k,s}$ be an $n$-vertex graph consisting of a complete graph $K_{k-s}$
 on the set $A\cup B$, $|A|=k-2s$, $|B|=s$ and a complete bipartite graph
 $K_{s, n-(k-s)}$ with parts $B$ and $C$
 where $A, B$ and $C$ form a partition of $V(H)$ (hence $|C|=n-(k-s)$ and
 $n\ge k$, $(k-1)/2\ge s\ge1$).
\end{construction}

The graph $H$ contains no cycle of size $k$ or larger and for $s\geq 2$ it is 2-connected.
Denote its size by $h(n,k,s)$.

They all (\cite{Kopylov},  \cite{WoodallA},  \cite{CV1991}) characterized the
structure of
 the extremal graphs in Theorem~\ref{th:Kopylov}. Namely either \dori
--- the blocks of $\Gn$ are $m$ complete graphs $K_{k-1}$ and a $K_r$, or
\dori
--- $k$ is odd, $r=(k+1)/2$ or $(k-1)/2$ and $q$ of the blocks of $\Gn$ are
 $K_{k-1}$'s and a copy of a $H_{n-q(k-2), k, (k-1)/2}$.

The strongest result on the field is due to Kopylov who also investigated 2-connected graphs.

\begin{theorem} [Kopylov  \cite{Kopylov}]
\label{th:Kopylov2}
 Suppose that $n\geq k\geq 5$ and the 2-connected graph $\Gn$
contains no cycles of length of $k$ or larger. Then
$$ e(\Gn)\leq \max\{ h(n,k,2), h(n, k, \lfloor \half(k-1)\rfloor)\}$$
and this bound is the best possible.

Moreover, only the graphs $H_{n,k,s}$ could be extremal, $s\in \{ 2, \lfloor (k-1)/2\rfloor\}$.
 \end{theorem}

This theorem was also conjectured by Woodall~ \cite{WoodallA} and he also proved it
for $n\geq (3k-5)/2$. It was also reproved much later in
 \cite{FanXWang}.

\subsection{Excluding $P_k$}

One of the oldest problems is the question of determining $\ext(n ,P_k)$.

\begin{theorem}[Erd\H{o}s and Gallai~ \cite{ErdGallai}]
\label{ErdGalPathTh}

If $\Gn$ is a graph containing no $\Pk$,
 ($k\ge 2$),
then
$$
 e(\Gn)\le {{k-2}\over 2}n
 $$
with equality if and only if $k-1$ divides $n$ and all connected
 components of $G$ are complete graphs on $k-1$ vertices.
\end{theorem}

Consider
 the $n$-vertex graph $G_n$ which is the union of $\lfloor n/(k-1)\rfloor$ vertex-disjoint $K_{k-1}$
 and a $K_r$ ($0\le r\le k-2$). 
If $\Tk$ is any connected $k$-vertex graph, then $\Tk\not\subseteq\Gn$.
Hence
\beq{TreeLower}\ext(n,\Tk)\ge {k-2\over 2}n -\reci8 k^2 .\eeq
In particular,
\begin{equation} 
 \ext(n,\Pk)\ge {{k-2}\over 2}n -\reci8 k^2.
 \end{equation} 

\begin{figure}[ht]
\begin{center}
\epsfig{file=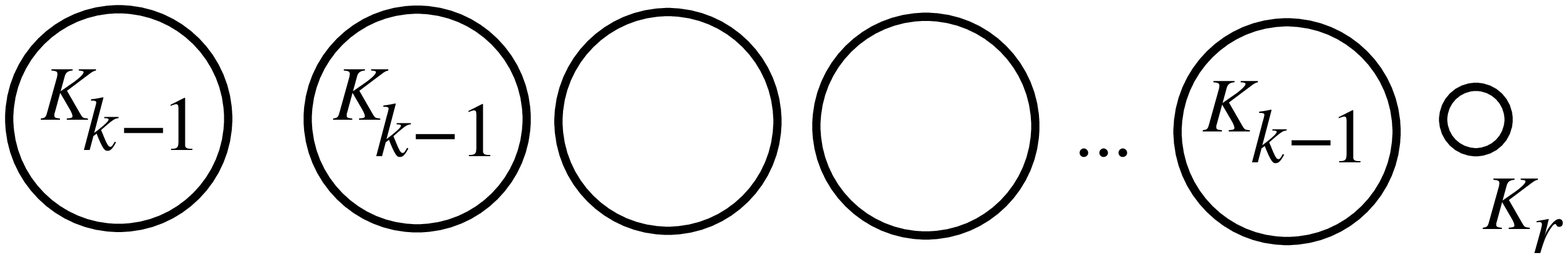,height=1cm}
\qquad\epsfig{file=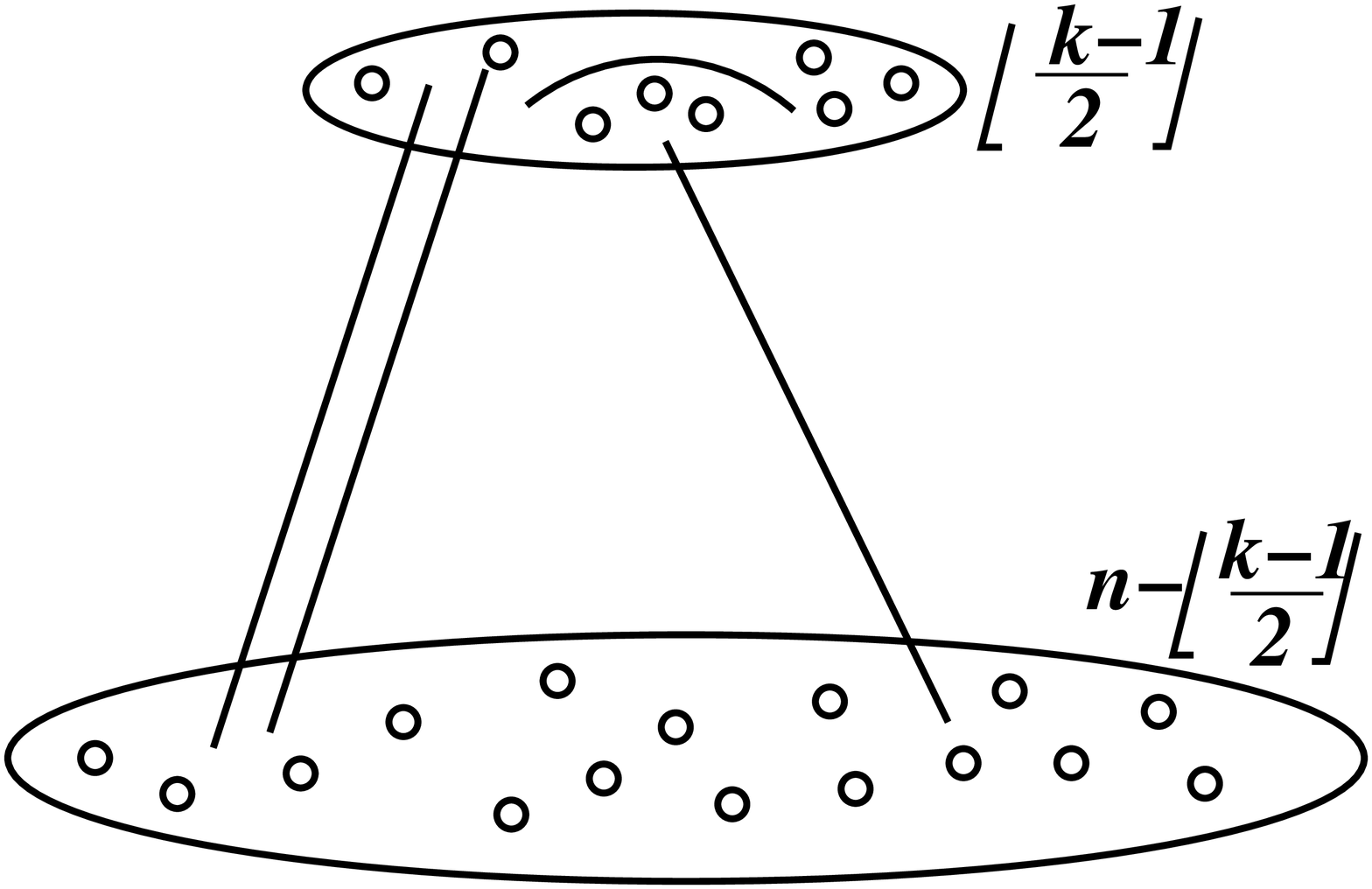,height=2cm}
\end{center}
\caption{Potential extremal graphs}
\label{PotExtrPk}
\end{figure}

If $k$ is even, then there are nearly extremal graphs having a
 completely different structure.
Namely, one can take a complete bipartite
 graph with partite sets $A$ and $B$ of sizes $|A|={{k-2}\over 2}$ and
 $|B|=n- {{k-2}\over 2}$ and add all edges in $A$, too (Fig.~\ref{PotExtrPk}).
Faudree and Schelp  \cite{FaudreeSchelpJCT75} proved that
 the extremal graph for $\Pk$ can indeed be obtained in this way for all $n$ and $k$.
They needed this to prove some Ramsey theorems on paths.
The variety of extremal graphs makes the solution difficult.

 \begin{theorem} [Faudree and Schelp  \cite{FaudreeSchelpJCT75} and independently Kopylov  \cite{Kopylov}]
\label{th:Kopylov_FaudreeSchelp} ${}$
Let $n\equiv r$ {\rm (mod $k-1$)}, $0\le r< k-1$, $k\geq 2$. Then
\beq{eq:Kop_path}
 \ext(n ,P_k)= \half (k-2)n - \half r(k-1-r).
 \eeq
 \end{theorem}

Faudree and Schelp also described the extremal graphs which are either
 \dori --- vertex disjoint unions of $m$ complete graphs $K_{k-1}$ and a $K_r$, or
 \dori --- $k$ is even and $r=k/2$ or $k/2-1$ and another extremal graphs can be obtained by
  taking a vertex disjoint union of $t$ copies of $K_{k-1}$ ($0\le t <m$) and a
  copy of $H(n-t(k-1), k-1, k/2-1)$.

\begin{theorem}[Kopylov~ \cite{Kopylov}]
\label{KopyThB}
Let $\Gn$ be a connected graph containing no $\Pk$,
 $(k\ge 4)$ and $n\ge k$.
Then
$$
 e(G)\le \max\{ h(n,k-1,1), h(n, k-1, \lfloor \half(k-2)\rfloor)\}$$
and this bound is the best possible.

Moreover, only the graphs $H_{n,k-1,s}$ could be extremal, $s\in \{ 1, \lfloor (k-2)/2\rfloor\}$.
  \end{theorem}

Balister, Gy\H{o}ri, Lehel and Schelp  \cite{BalisGyoriLehSchelp}
 also provided the extremal structures.

\subsection{Proof ideas}

Let $\extcon(n,\Pk)$ be the maximum number of edges in connected,
 $n$-vertex, $\Pk$-free graphs, and let $\extconn(n,{\mathcal C}_{\ge k})$ denote the
 maximum number of edges in 2-connected, $n$-vertex, ${\mathcal C}_{\ge k}$-free graphs.
Determining these functions give upper bounds for
 $
 \ext(n,\Pk)$
and $\ext(n,{\mathcal C}_{\ge k})$.

Indeed, every $\Pk$-free graph is a vertex disjoint union of $\Pk$-free
 components, we have
$$
 \ext(n,\Pk)=\max_{\sum n_i =n, \, \, n_i\ge 1} \sum \extcon(n_i,\Pk).
 $$
Similarly, a maximal ${\mathcal C}_{\ge k}$-free graph is connected
and every connected graph is a cactus-like union of 2-connected
blocks (and edges)
 so we have
\begin{equation}\label{eq:coonectedcycle}
 \ext(n,{\mathcal C}_{\ge k})=\max_{\sum (n_i-1) =n-1,\, \, n_i\ge 2}
      \sum \extconn(n_i, {\mathcal C}_{\ge k}),
 \end{equation}
where we define $\extconn(2,\Cg k)=1$.

Let $G$ be a connected, $n$-vertex, $P_k$-free graph.
Add a new vertex to it and join to all other vertices.
We obtain $G_{n+1}$ with $e(G_{n+1})=e(G)+n$.
This new graph has no cycle of length exceeding $k$
 and its connectivity is one larger than that of $\Gn$.
We obtain
\begin{equation}\label{eq:pathANDcycle}
 \ext(n,P_k)+n \le \ext(n+1,{\mathcal C}_{\ge k+1})
 \end{equation}
and
\begin{equation}\label{eq:2conn}
 \extcon(n,P_k)+n \le \extconn(n+1, {\mathcal C}_{\ge k+1}).
 \end{equation}

So Theorem~\ref{ErdGallaiCyc} and \eqref{eq:pathANDcycle} imply Theorem~\ref{ErdGalPathTh}.
Similarly, Theorem~\ref{th:Kopylov} and \eqref{eq:pathANDcycle} imply Theorem~\ref{th:Kopylov_FaudreeSchelp}.

The upper bounds for $\extconn(n,{\mathcal C}_{\ge k+1})$ yield upper bounds
 for $\extcon(n,\Pk)$.
 (Actually, \eqref{eq:2conn} and Theorem~\ref{th:Kopylov2} lead to the solution
 of $\extcon(n,\Pk)$, Theorem~\ref{KopyThB}).

Again Theorem~\ref{th:Kopylov2} and \eqref{eq:coonectedcycle} lead to
 Theorem~\ref{th:Kopylov} which is obviously stronger than Theorem~\ref{ErdGallaiCyc}.

Finally, the proof of Theorem~\ref{th:Kopylov2} uses induction on $n$ and $k$, by deleting small degree vertices, contracting edges,
and finally applying P\'osa's theorem on Hamiltonian graphs.

\subsection{Generalizations}

In a recent work Lidick\'y, Hong Liu and Cory Palmer~ \cite{LLP:forests}
 determined the exact Tur\'an number (and the unique extremal graph)
 when the forbidden graph $L$ is a {\it linear forest}, each component
 is a path.
They also considered {\it star-forests}.

Gy\'arf\'as, Rousseau, and Schelp  \cite{GyRS}
determined $\ext(K(m,n), P_k)$ for all $m,n,k$.
Their formula and proof are rather involved, they distinguish 10 subcases.

\section{Excluding trees}\label{TreeExtremal}

Here we shall discuss two extremal problems on trees: the Erd\H{o}s-S\'os conjecture
 and the Loebl-Koml\'os-S\'os conjecture.

\subsection{Erd\H{o}s-S\'os conjecture}

We have already discussed the Erd\H{o}s-Gallai theorems. Since the
extremal numbers for $\Pk$ and for the star $K_{1,k-1}$ are roughly the
same, this led Erd\H{o}s and T. S\'os to the following famous
conjecture.

\begin{conjecture}[Erd\H{o}s-S\'os  \cite{ErdSmole}] \label{ErdSosConj}
Let $\Tk$ be an arbitrarily fixed $k$-vertex tree. If a graph $\Gn$
contains no $\Tk$, then
\begin{equation}\label{sharpconj} e(\Gn)\le \half(k-2)n.
\end{equation}
\end{conjecture}

As we have seen -- by \eqref{TreeLower} -- the disjoint union of
complete graphs $K_{k-1}$ shows that $\ext(n,T_k)\ge \half
(k-2)n-\reci8 k^2$. Though several partial cases were settled, the
upper bound was unknown until Ajtai, Koml\'os, Simonovits,  and
Szemer\'edi proved:

\begin{theorem}[Main Theorem, Sharp
 \cite{AKSSzApprox,AKSSzDenseBlock,AKSSzSharp}]\label{SharpTh}
There exists an integer $k_0$ such that if $k>k_0$ and
 $\Tk$ is an arbitrarily fixed
 $k$-vertex tree, and the graph $\Gn$ contains no $\Tk$, then
\beq{sharpfo}e(\Gn)\le \half(k-2)n.
\eeq
\end{theorem}

Below we list a few subcases where this conjecture is verified, but we
do not try to give a complete list.

\begin{theorem}[Sidorenko  \cite{SidorTree}]
If $\Tk$ has a vertex $x$ connected to at least $k/2$ vertices of degree $1$ (i.e., leaves)
 then the Erd\H{o}s--S\'os conjecture holds for this $\Tk$.
\end{theorem}

\begin{theorem}[McLennan  \cite{LennanTree}]
If the diameter of $\Tk$ is at most 4, then
the Erd\H{o}s--S\'os conjecture holds for this $\Tk$.
\end{theorem}

Dobson (and coauthors) have several results in this area, under some
strong condition of sparsity. We mention only the Brandt-Dobson
theorem  \cite{BrandtDobson}, or Sacle and Wozniak,  \cite{Wozn96},  \cite{SacWoz}.

\subsection{Sketch of the proof of Theorem \ref{SharpTh}}

We are given a $\Tk$, and a $\Gn$ violating \eqref{sharpfo}.
We wish to embed $\Tk$ into $\Gn$ ($\TembG$).
The proof is very involved and will be given in three rather long
papers. The following weakening plays a central role.

\begin{theorem}[$\eta$-weakening  \cite{AKSSzApprox}]\label{AKSSzApproxThm}
 For any (small) constant $\eta>0$ there exists a $k_0(\eta)$
such that for $n\ge k>k_0(\eta)$ , if
\beq{ESeta} e(\Gn)>\half(k-2)n+\eta kn,\eeq
then each $k$-vertex tree $\Tk$ is contained in $\Gn$.
\end{theorem}

(a) First, in  \cite{AKSSzApprox} we prove this theorem. If, in
addition, we assume that $\Gn$ is dense: for some $c>0$, $k>cn$, then
we can apply the Szemer\'edi Regularity Lemma  \cite{SzemRegu}.
The proof of this theorem follows basically the line which was later used to
prove the Loebl Conjecture, by Ajtai, Koml\'os and Szemer\'edi
 \cite{AKSLoebl95} and later by Yi Zhao  \cite{YiZhao11}.
Also it was used in  the Koml\'os-S\'os conjecture by Diana Piguet and Maya Stein  \cite{PiguetStein08},
 Cooley  \cite{Cooley09}, Hladk\'y and Piguet  \cite{HladkyPiguet}, in stronger and stronger form, and
now the publication of that proof is almost finished by
Hladk\'y, Koml\'os, Piguet, 
  Simonovits, Stein, and Szemer\'edi  \cite{HladKomPiguSimSteinSzem}.

(b) In the second part,  \cite{AKSSzDenseBlock} we prove several theorems
asserting that under some very special conditions
$\Tk\subseteq\Gn$. Some of these steps are ``stability arguments''.

Analyzing the proof of Theorem \ref{AKSSzApproxThm}, shows that
either we can gain at some points, in some of the estimates $\eta kn$
edges, and therefore Theorem \ref{AKSSzApproxThm} (more precisely, its
slightly modified proof) implies the sharp version,
 Theorem~\ref{SharpTh}, or else $\Gn$ must have a very special
 structure: it contains a smaller copy of the conjectured extremal
 graphs: for some $m\approx k$,

($b_1$) either it contains a $\Gm$ which is almost a $K_m$;

($b_2$) or a $\Gm$ which is almost a $K(m/2-\eps m,m+\de m)$ .

(c) In both cases, if many edges connect $\Gn-\Gm$ to $\Gm$,
then we can embed $\Tk$ into $\Gn$, embedding a smaller part of $\Tk$
outside of $\Gm$, a larger part in the dense $\Gm$, concluding that
$\TembG$.

(d) If, on the other hand, we have found such a
``mini-almost-extremal'' $\Gm\subseteq\Gn$, but $e(\Gm,\Gn-\Gn)$ is
``small'', then we prove that
$$e(\Gn-\Gm)>\half(n-m)(k-2).$$ Hence we may forget the larger $\Gn$:
replace it by the smaller $\Gn-\Gm$. (In other words,
we can apply induction on $n$.)

(e) The real difficulty comes when we have sparse graphs:
$e(\Gn)=o(n^2)$. Then we partition $V(\Gn)$ into three parts: $\CC$
contains the vertices of high degrees, $\BB$ contains a part of
$V(\Gn)$ not containing dense subgraphs, and therefore behaving in a
pseudo-random way, and $\AA$ behaves very similarly to the graphs we
have in the dense cases.

\subsub How do we handle the dense case?//

(i) Applying the Regularity Lemma to $\Gn$, we get a so called Cluster
Graph $\Hnu$. If this cluster graph has an (almost)-1-factor, then we can
relatively easily embed $\Tk$ into $\Gn$, using the extra $\eta kn$
edges of \eqref{ESeta}.

(ii) Next we extend this case to a more general situation, when $\Gn$
contains a so called Generalized 1-factor. We can prove the
$\eta$-weakening in this case as well.

(iii) If the Cluster Graph $\Hnu$ does not contain an almost-1-factor,
then we apply the Gallai-Edmonds structure-theorem (on graphs without
1-factors) to $\Hnu$. In this case we can either embed $\Tk$ into
$\Gn$ directly, or reduce this case to Case (ii) above. Case (iii) is a very
important subcase, with 3-4 subsubcases (depending on, how do we count
them). Some of them go back to Case (ii) and in some others we directly
(pseudo-greedily) embed $\Tk$ into $\Gn$.

\subsection{Koml\'os-S\'os conjecture on median degree}

The Koml\'os-S\'os conjecture was already formulated in Section
\ref{OtherTypes}. This is a generalization of the Loebl conjecture:

\begin{conjecture}[Loebl--Koml\'os-S\'os Conjecture  \cite{EFLS95}]

If $\Gn$ has at least $n/2$ vertices of degree at least $k-1$, then $\Gn$ contains all the $k$-vertex trees $\Tk$.
\end{conjecture}

The authors of  \cite{HladKomPiguSimSteinSzem}
plan to write up the sharp version as well, which asserts the following.

\begin{theorem}

If $k$ is sufficiently large, then the Loebl--Koml\'os--S\'os conjecture is
true.
\end{theorem}

\begin{remarks}

(a) The Loebl conjecture originates from a paper of Erd\H{o}s, F\"uredi,
Loebl, and T. S\'os, on the discrepancy of trees  \cite{EFLS95}.

(b) P\'osa's theorem on the existence of Hamiltonian cycles also is -- in some
sense -- a theorem asserting that if $G$ has many vertices of sufficiently
high degree, then it is Hamiltonian.  There were earlier cases, when Woodall
 \cite{WoodallB}, proved an Erd\H{o}s--Gallai type theorem on cycles, using the
condition that there are many vertices of high degree. Also, Erd\H{o}s,
Faudree, Schelp, and Simonovits -- trying to prove some Ramsey type theorems,
-- found a similar statement  \cite{ErdFaudSchelpSim}, but not for all the
trees, only for the paths, and they proved there an almost sharp
theorem. Their sharp 
conjecture was later proved by Hao Li.

(c) There were many important steps to reach the theorem above. We
 should mention here Ajtai--Koml\'os--Szemer\'edi,  \cite{AKSLoebl95},
 then Yi Zhao  \cite{YiZhao11}, next Piguet and Stein
  \cite{PiguetStein08},  \cite{PiguetStein12}, Cooley  \cite{Cooley09},
 Hladk\'y and Piguet  \cite{HladkyPiguet}, and many others.
\end{remarks}

\section{More complex excluded subgraphs}\label{MoreComplex}

In this section we present three theorems, each leading to a
 reduction method to prove new results from old estimates.
Still there is no general theory to determine the bipartite Tur\'an numbers.

The three results we pick are the Erd\H{o}s-Simonovits cube theorem
(Theorem~\ref{CubeThB}), $\ext(n, Q_8)=O(n^{8/5})$ which led to the
Erd\H{o}s-Simonovits reduction, the Faudree-Simonovits theorem
(Theorem~\ref{FauSimTheta}) concerning Theta graphs, a generalization
of the Erd\H{o}s-Bondy-Simonovits theorem, $\ext(n,
C_{2k})=O(n^{1+(1/k)})$, and F\"uredi's theorem
(Theorem~\ref{FureM11Th}) on two levels of the Boolean lattice which
implies a general upper bound $\ext(n, L)=O(n^{2-(1/r)})$ for any
graph $L$ with vertices of degrees at most $r$ on one side of $L$.

\subsection{The Erd\H{o}s-Simonovits Reduction and the Cube theorem
}\label{CubeReduS}

We have already mentioned Theorem \ref{CubeTh},
on the extremal number of the cube. Here we formulate a sharpening of it.

\begin{theorem}[\cite{ErdSimCube}] \label{CubeThB}

Let $Q_8$ denote the graph determined by the
8 vertices and 12 edges of a cube, and $Q_8^+$ denote the graph obtained by joining two
opposite vertices of this cube. Then
$$ \ext(n,Q_8)\le \ext(n,Q_8^+)= O(n^{8/5}).$$
\end{theorem}

\BBalAbra{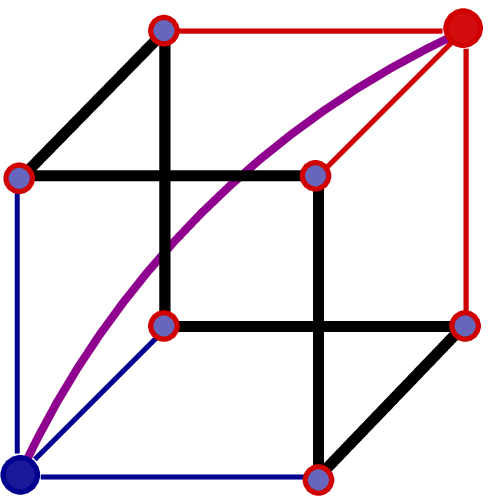}{24} One reason why Erd\H{o}s and Simonovits considered the extremal
problem of the Cube graph was that this was one of Tur\'an's originally
posed problems. The reason that $Q_8^+$ was also considered was that
Erd\H{o}s and Simonovits got it for free: their proof of Theorem
\ref{CubeTh} gave the same upper bound for $Q_8^+$.

Let $L$ be a
bipartite graph with partite sets $X$ and $Y,$ and let $K_{t,t}* L$ 
denote the graph obtained by completely joining one partite set of $K_{t,t}$
to $X$ and the other to $Y.$

\begin{theorem}[Erd\H{o}s and Simonovits Reduction Theorem  \cite{ErdSimCube}]\label{CubeRedu}
 ${}$ If $L$ is a bipartite graph with
$\ext(n,L)=O(n^{2 -a})$, $a\leq 1$,
and $b$ is defined by ${1\over b}={1\over a}+t$,
then $\ext(n, K_{t,t}*L) 
=O(n^{2-b})$.
\end{theorem}

The proof can go by induction on $t$ and by counting the number of $C_4$'s.

Let $H$ be the graph obtained by deleting just three independent edges from $K_{4,4}$.
Since $H=K_{1,1}*C_6$, Theorem \ref{CubeRedu} and $\ext(n, C_6)=O(n^{4/3})$
 (Corollary~\ref{co:8}) imply Theorem \ref{CubeThB}.

Since $\ext(n,L)=O(n)$ if $L$ is a tree, so we have the following result:

\begin{corollary}\label{co:treeplus}

For any tree $L$, $\ext(n, K_{t,t}* L) 
 =O(n^{2-(1/(t+1))})$.
\end{corollary}

\BBalAbra{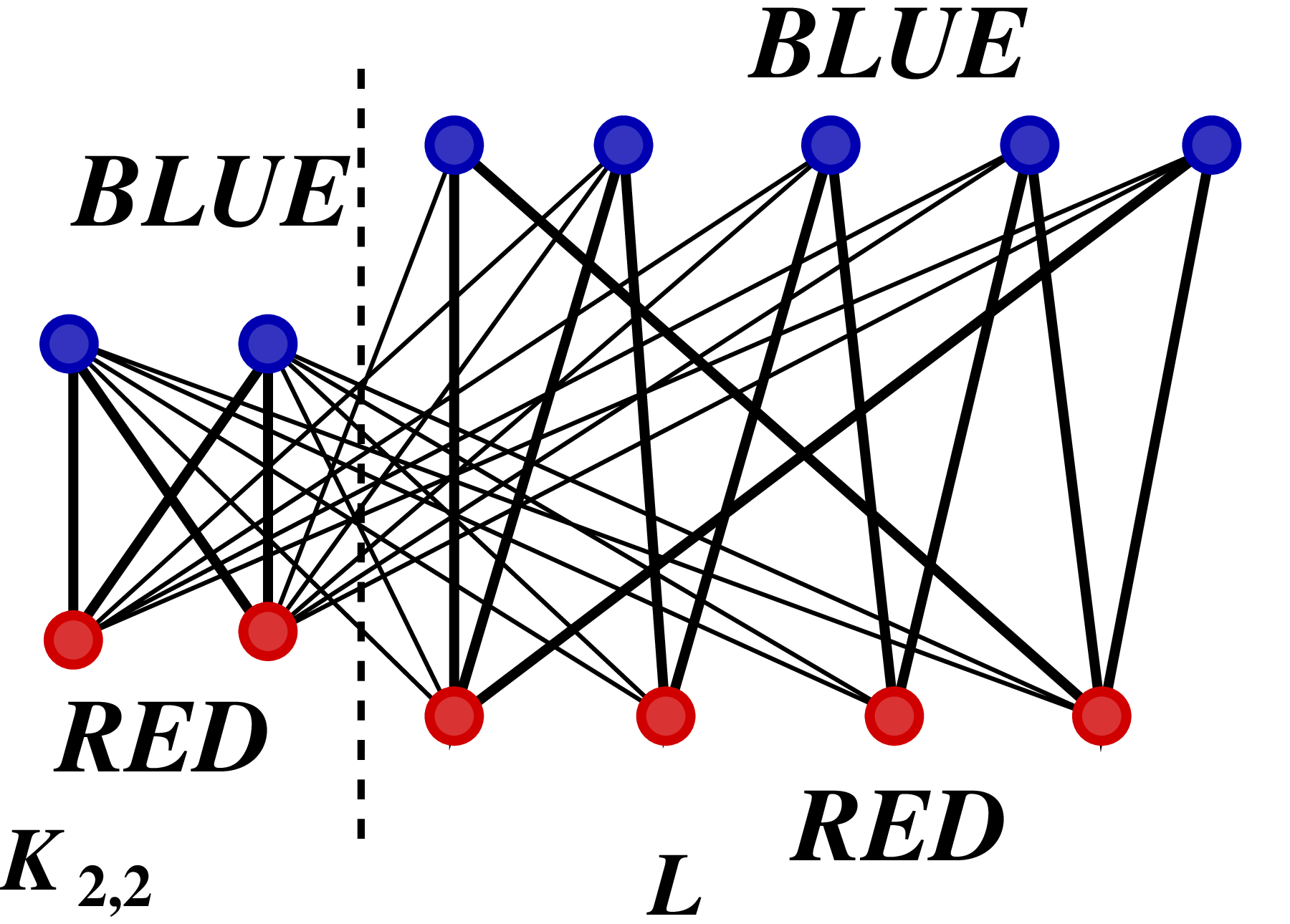}{40} This can be considered as a generalization of
 the K\H ov\'ari-T. S\'os-Tur\'an
theorem, since for $L=K_2$ we have 
$K_{t+1,t+1}=K_{t,t} * K_{1,1}$.
Since $Q_8-e$ is a subgraph of $K_{1,1}* P_6$,
Corollary \ref{co:treeplus} implies

\begin{corollary}[Erd\H{o}s and Simonovits]
 $\ext(n,Q_8-e)=O(n^{3/2})$.
\end{corollary}

Further,

\begin{theorem}[Erd\H{o}s]

Delete an edge from $\KK aa$. For the resulting $L=K_{a,a}-e$ we have
$\ext(n,L)=O(n^{2-{1\over a-1}}).$
\end{theorem}


Indeed, for $a\geq 3$ the graph $K_{a,a}-e$ is a subgraph of
$K_{a-2,a-2}* P_4$.

Since $K_{b,b}-K_{a,a}$ (for $b-2\ge a\geq 1$) can be written as
 $K_{b-a-1,b-a-1}*T$ where $T$ is a double star, Corollary~\ref{co:treeplus} also implies that
$\ext(n, K_{b,b}-K_{a,a})=O(n^{2-(1/(b-a))})$.
For this important case F\"uredi and West gave a sharper upper bound.

\begin{theorem}[\cite{FureWest}]

For every $n\geq b> a$ we have
$$\ext(n, K_{b,b}-K_{a,a})\leq
  \half(b+a-1)^{1/(b-a)}n^{2-(1/(b-a))} + \half(b-a-1)n.$$
\end{theorem}

In particular, it gives $\ext(n,
K_{3,3}-e)\le \half \sqrt{3} n^{3/2}+O(n)$. This was further improved
by J. Shen  \cite{Shen_K33MinusEdge} to
$$
  \ext(n, K_{3,3}-e)\le { \sqrt{15}\over 5} n^{3/2}+O(n).
 $$
He also showed that $\ext(n, n, K_{3,3}-e)\le (4/ \sqrt{7}) n^{3/2}+(n/2).$

Pinchasi and Sharir extended the cube theorem, using a somewhat different proof: 

\begin{theorem}[Pinchasi and Sharir  \cite{PinchasiSharir}]\label{th:PinchasiSharir}

A bipartite graph $G[A,B]$ with $|A|=m$ and $|B|=n$, not containing the cube $Q$ has $$O(n^{4/5}m^{4/5}+mn^{1/2}+nm^{1/2})$$ edges.
\end{theorem}

Another, more explicit proof for Theorem~\ref{th:PinchasiSharir} was presented
in  \cite{FureCube}.

\subsub Historical remarks://

(a) Erd\H{o}s and Simonovits first proved the Cube
theorem, using the $\Delta$-almost-regularization.

(b) It seems that Theorem \ref{CubeRedu}
covered all the cases known until that point.

(c) This (i.e. the Cube Recursion Theorem) was the first case, where one got an exponent, different from
$2-(1/a)$ and $1+(1/a)$. Actually, Erd\H{o}s thought earlier, that all
the exponents must be of this form, (see  \cite{Erd1975-42}).
This was disproved in their paper  \cite{ErdSimCube}: not by the cube, since there is no good lower
bound for the cube: even $\ext(n, Q_8)/ n^{3/2}\to\infty$ is not
known. However, a more complicated example, for which the lower bound
-- using random graphs -- was good enough, disproved Erd\H{o}s' conjecture.
Actually, one thinks that each rational $\alpha\in(0,1)$ is extremal
exponent for some finite $\LL_\alpha$, see Conjecture~\ref{RatiExpo}.

To disprove the Erd\H{o}s conjecture
  concerning the exponents are of the form $1+(1/a)$ or $2-(1/a)$
it is enough to notice that we have graphs $H$ with
$$c_Hn^{(8/5)-\eps(H)}<\ext(n,H)=O(n^{8/5}),\qquad(c_H>0,\,\eps(H)<{1\over 10} ).$$

\BBalAbra{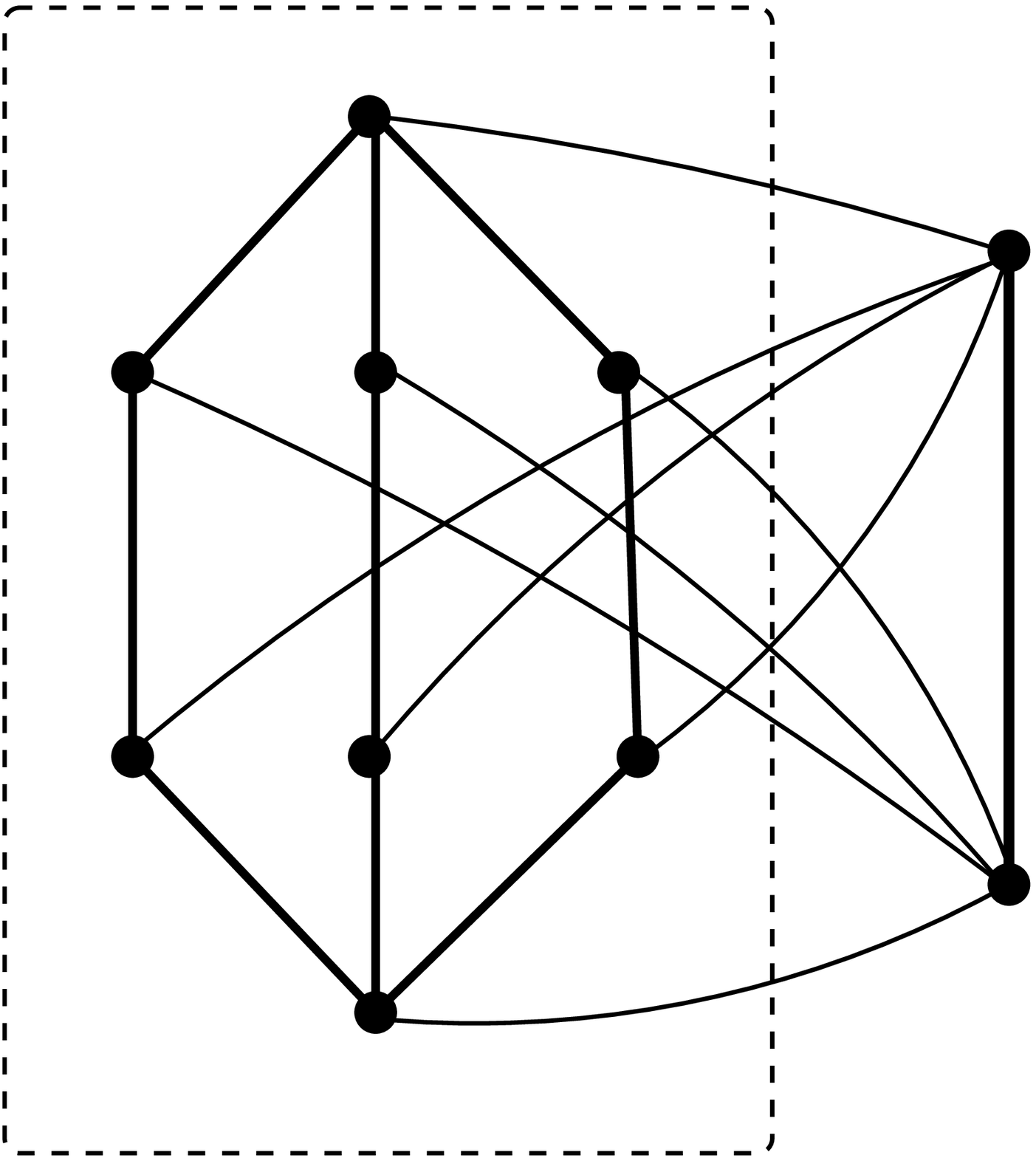}{30} More generally, consider the graph $H(t,\ell)$
obtained by connecting a $\Theta(3,\ell)$ to $K(t,t)$, as described in
Theorem~\ref{CubeRedu}.  By the Theorem~\ref{ErdRenyEv} (lower bound) and
 Theorem~\ref{CubeRedu} and Theorem~\ref{FauSimTheta} (upper bound) we obtain
$$c_{\ell,t}n^{2-{{2\ell+2t \over 3\ell+t^2+2t(\ell+1)-1}}} <\ext(n,H_{t,\ell})
\le \ti c_{\ell,t}n^{2-{2\over 2t+3}}. $$
So all the numbers $2-{2\over2t+3}$ are points of accumulations of exponents,
in this sense.
Actually, applying this argument with $t=1$, $\ell=3$,
  we get a simple counterexample, with the upper bound $O(n^{8/5})$ and
a lower bound  $cn^{2-(8/17)}$ ($c>0$).

(d) In  \cite{ErdSimCubeWat}, Erd\H{o}s and Simonovits proved the
Supersaturated graph theorem (see Section \ref{SupersatS})
 corresponding to the cube, thus providing a second proof of the
Cube Theorem, that needed ``less regularization''. 

\subsection{Theta graphs and the Faudree-Simonovits reduction}

\index{Faudree--Simonovits recursion}
There is an alternative proof for the Bondy--Simonovits Theorem in  \cite{FaudreeSimCCA}.
This proof enabled a generalization 
to $\Theta$-graphs.
Recall that $\Theta_{k,\ell}$ denotes the graph consisting of $\ell$ paths of
length $k$ with the same endpoints but no inner intersections.
We have $v(\Theta_{k,\ell})= 2+(k-1)\ell$ and $e(\Theta_{k,\ell})= k\ell$.

\begin{theorem}[Theta-graph, Faudree--Simonovits  \cite{FaudreeSimCCA}]
\label{FauSimTheta}
 For fixed $k$ and $\ell\geq 2$ one has
$\ext(n,\Theta_{k,\ell})
=O(n^{1+(1/k)})$.
\end{theorem}

 \BalAbra{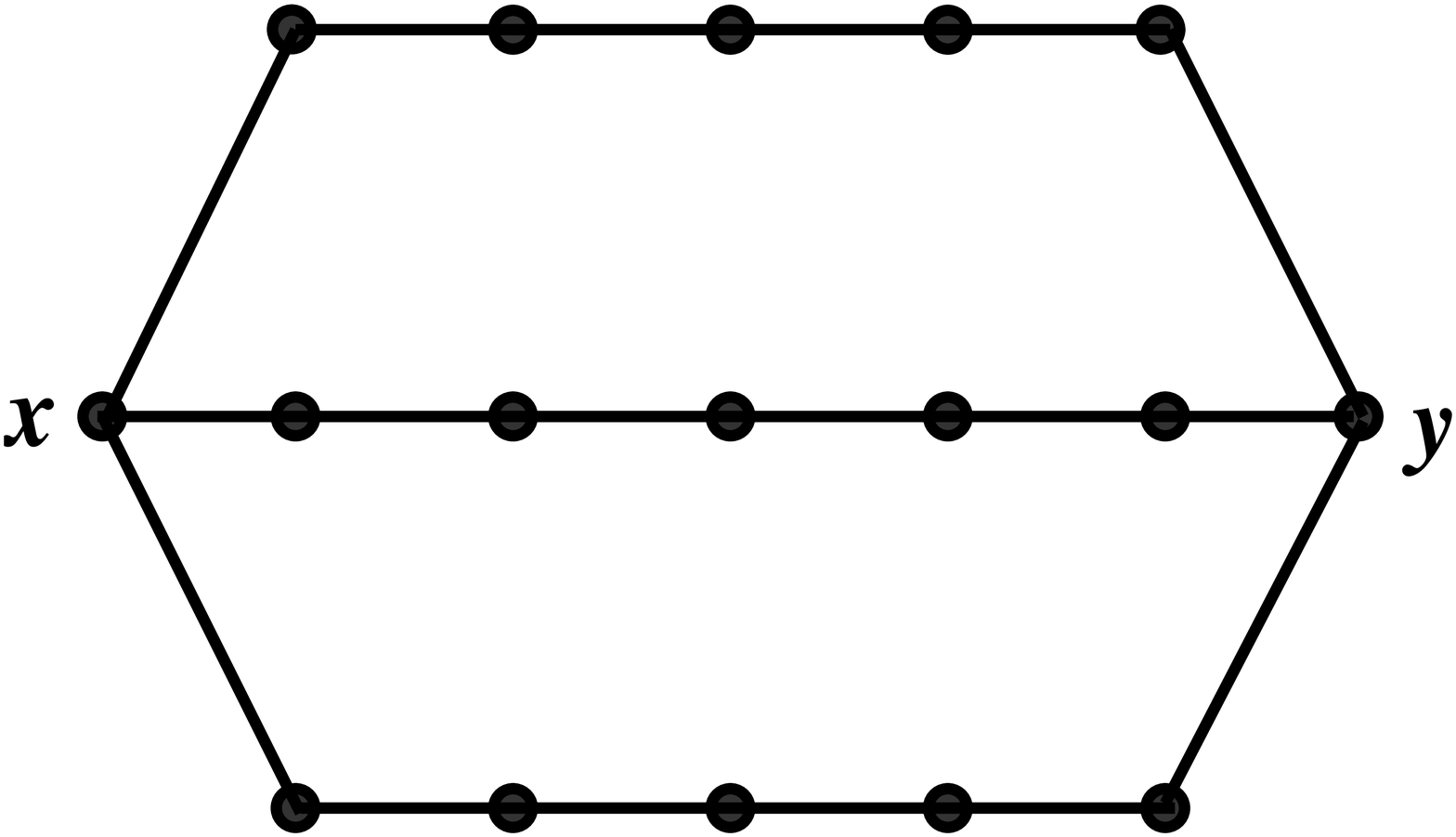}
This exponent is conjectured to be the best possible, see Conjecture~\ref{ThetaCon}.

Applying the Erd\H{o}s-R\'enyi Random Lower bound (Theorem \ref{ErdRenyEv}) in its simpler form
 to $\Theta_{k,\ell}$ 
  we get
$$\ext(n,\Theta_{k,\ell})> 
c_{k,\ell}n^{1+{1\over k}-{2\over k\ell}},$$
asymptotically matching the upper bound's exponent.

\BBalAbra{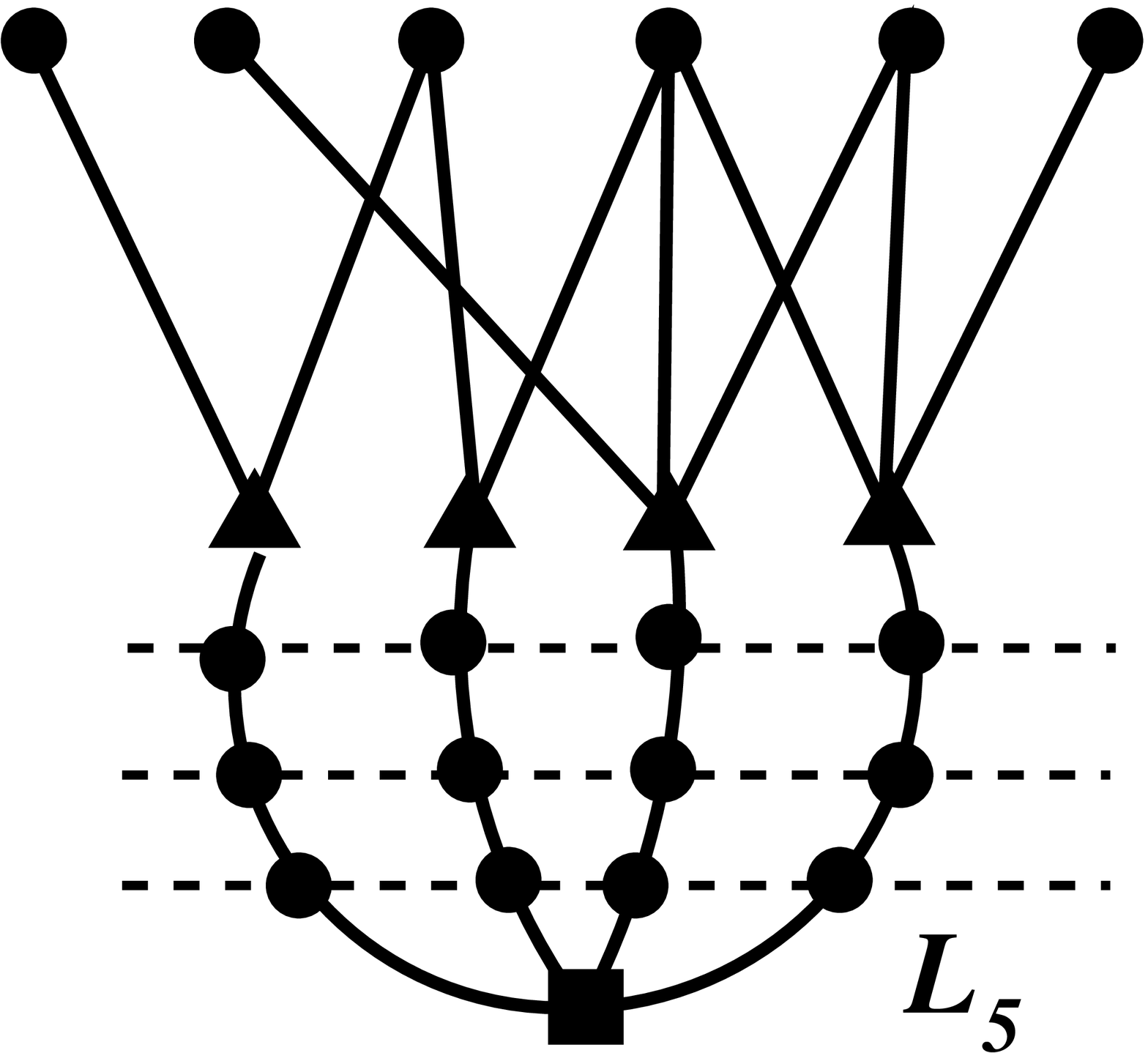}{30}
The proof of 
Theorem \ref{FauSimTheta} came from a ``Recursion'' theorem, asserting that if
one knows good upper bounds for an $L$, and $L^*$ is built from $L$ in
a simple way, then one has a good 
upper bound on $\ext(n,L^*)$ as well.

\begin{definition}\label{FauSimD}

Let $L$ be a bipartite graph, with a fixed 2-colouring $\psi$ in
RED-BLUE with $h$ RED vertices.
Let \mbox{$x\not\in V(L)$} be a vertex from which $h$ independent
paths of $k-1$ edges go the RED vertices of $L$, 
(these paths intersect only in $x$).
Denote the obtained graph by $\Lt k$.
\end{definition}

\begin{theorem}[Faudree-Simonovits Reduction, Trees  \cite{FaudreeSimCCA}]

If $L$ is
a tree, then $$\ext(n,\Lt k)=O(n^{1+(1/k)}).$$
\end{theorem}

The Theta graph $\Theta_{k,\ell}$ is obtained from a star of $\ell$ edges.
One has to be cautious with the next theorem, see Remark \ref{caution}.

\begin{theorem}[Faudree-Simonovits Reduction, General
  Case  \cite{FaudreeSimCCA,FaudreeSimDegII}]
\label{ThTWO}
 Let $L$ be an arbitrary bipartite graph with a fixed coloring $\psi$
and assume that
\begin{equation}
\label{JJJ}\ext^*(n,L) = O(n^{2-\alpha}).
\end{equation}
Then for
\begin{equation}\label{KKK}
\beta = {\alpha+ \alpha^2 + \dots + \alpha^{k-2} \over
1+\alpha+\alpha^2+\dots+\alpha^{k-2} } $$
we have
$$
\ext(n,L_k(L,\psi)) \le \ext^*(n,L_k(L,\psi)) = O(n^{2-\beta} ).
\end{equation}
\end{theorem}

\begin{remark}\label{caution}
 Most probably, this recursion is never sharp but for trees.
In its proof one has to apply standard arguments to subgraphs of
$K(m,n)$ where $n\gg m$. We very seldom have matching lower and upper
bounds in such cases.
\end{remark}

\subsection{A universal graph and dependent random choice  
}\label{M11Excluded}

\leftskip=56mm

Erd\H{o}s asked the following question: what are the extremal numbers
for the two graphs on the left: The left one will be called $M_{10}$, the
right one $M_{11}$ and they are described as special cases of the
following

\leftskip=0mm

\vskip-2cm
\noindent
\epsfig{file=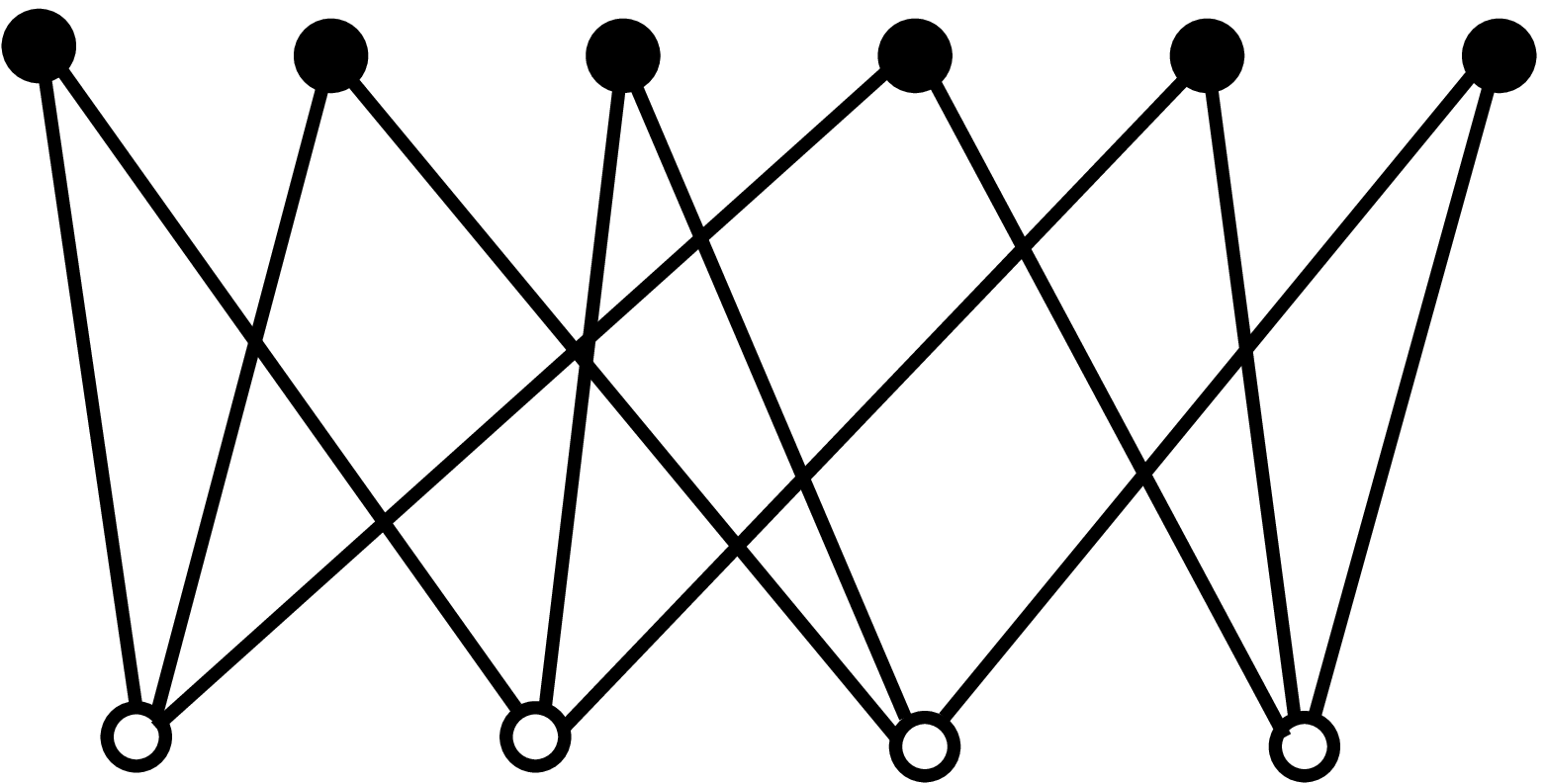,width=25mm,height=15mm}
\epsfig{file=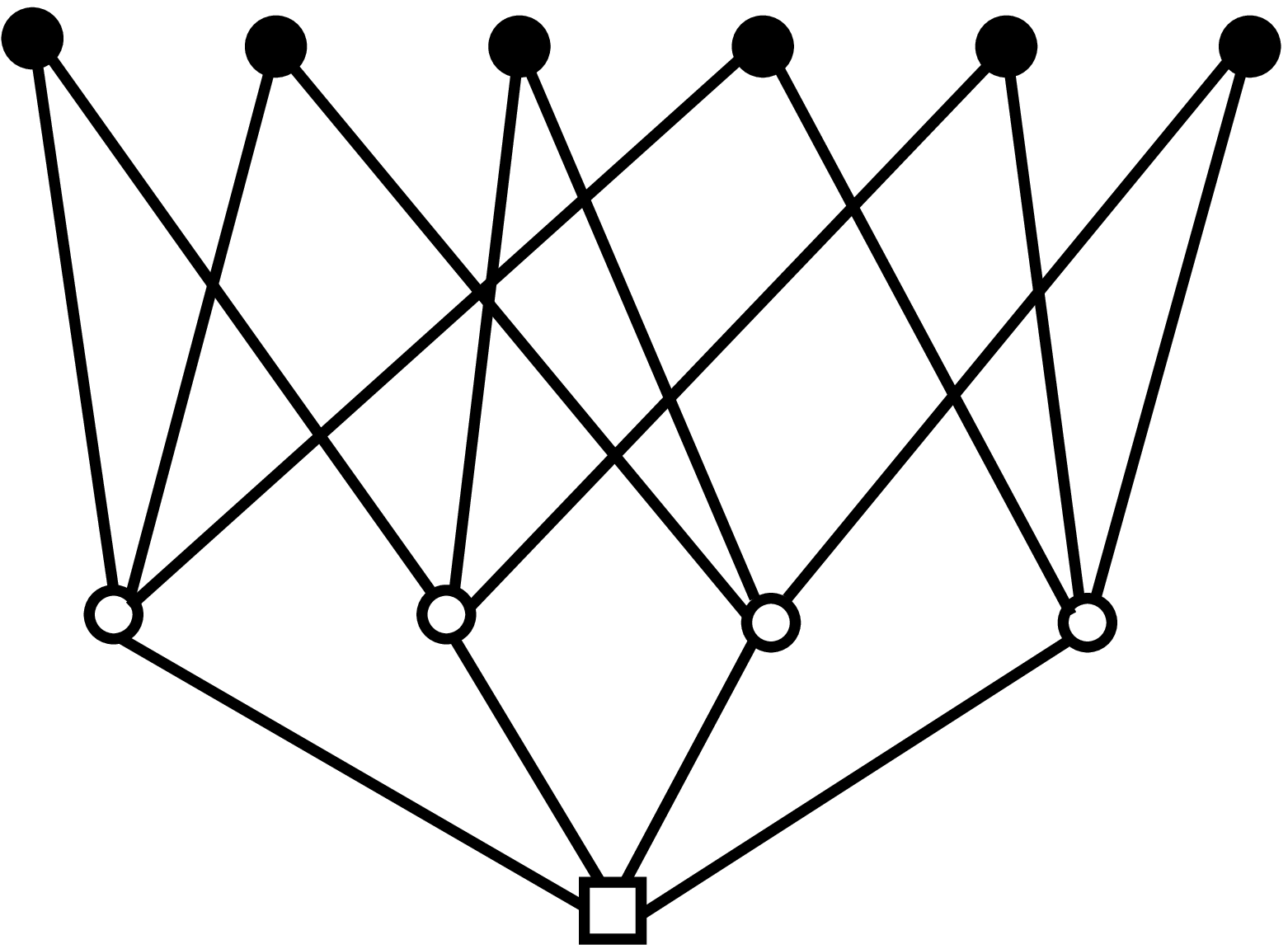,width=25mm,height=17mm}

\begin{definition}\label{univgr}

Let $k,r$ and $t$ be given positive integers. $U(k,r,t)$ is
obtained from the $k$ vertices $x_1,\dots,x_k$ by joining to each of
the $r$-element subsets of $\{x_1,\dots,x_k\}$ $t$ distinct vertices
$y^i_{i_1,\dots,i_r}$. $U^+(k,r,t)$ is obtained from $U(k,r,t)$ by joining a new vertex $w$ to all $x_h$,
$h=1,\dots,k$.
\end{definition}

\begin{problem}[Erd\H{o}s]

Determine (or estimate) $\ext(n,L)$ for $L:=M_{10}=U(4,2,1)$
and $L:=M_{11}=U^+(4,2,1)$.
\end{problem}

One could ask for the motivation: why
these graphs? Perhaps having obtained the cube theorem, we had good
upper and lower bounds only in very special cases, when $L$ contained
some sample graphs -- say a $C_4$ for which we have already provided
sharp lower bounds. 
$U(4,2,1)$ clearly needed a new approach, and e.g. $U^+(k,2,1)$ contains
many $C_4$'s but the earlier methods did not yield appropriate upper
bounds.
F\"uredi  \cite{FurediL11} answered this question proving that
$\ext(n,U^+(k,2,1)))<k^{3/2}n^{3/2}$. More generally,

\begin{theorem}[F\"uredi  \cite{FurediL11}]\label{FureM11Th}

Let $U^+(k,r,t)$ be the universal bipartite graph from Definition \ref{univgr}.
Then there exists a $c= c_{r}^{k,t}>0$ such that
\beq{eq:L11_UpperGeneral}
  \ext(n,U^+(k,r,t))< c n^{2-(1/r)}.
 \eeq
Concerning the Erd\H{o}s question we have $\ext(n,U^+(k,2,1))<k^{3/2}n^{3/2}$ and
more generally
\beq{eq:L11_Upper_r=2}
  \ext(n,U^+(k,2,t))<
n^{3/2}\cdot\sqrt{ tk(k-1)^2 + 2(k-2)(k-1) \over 8}+
n{k-1\over 4}.
 \eeq
\end{theorem}
Multiplying each vertex $(k-1)$ times in a $C_4$-free graph we get a $U^+(k,2,1)$-free graph
 which yields $\ext(n,U^+(k,2,1))\ge \Omega(k^{1/2}n^{3/2})$.

Erd\H{o}s had the more general conjecture

\begin{conjecture}
[Erd\H{o}s,
  \cite{ErdRome}, see also  \cite{ErdSimCubeWat},  \cite{SimWat}]
\label{co:716}

If every subgraph of the bipartite graph $L$ has a vertex of degree at most $r$, then
$$\ext(n,L)=O(n^{2-(1/r)}).$$
\end{conjecture}

The upper bound
\eqref{eq:L11_UpperGeneral} for the universal graph immediately gives

\begin{corollary}\label{co:AhasDegrees=r}

If $L$ is bipartite and has a 2-coloring where in the first color class
all but one vertex is of degree at most $r$, then
$$\ext(n,L)=O(n^{2-(1/r)}).$$
\end{corollary}

Indeed, all such graphs can trivially be embedded into an appropriate
$U^+(k,r,t)$.

Alon, Krivelevich, and Sudakov  \cite{AlonKriveSudaDep}
gave a new probabilistic proof (for graphs where on one side all vertices
 are of degree at most $r$) with a better constant $c_r^{k,t}$.
Their proof method became known as ``dependent random choice''; for a
survey see  \cite{FoxSudDRCSurv}.

\begin{lemma}[Dependent random choice, see, e.g.,  \cite{FoxSudDRCSurv}]
\label{DependentRandom}

 Let $k, t, r$ be positive integers. Let $\Gn$ be
a graph with $n$ vertices and average degree $d$, $d$ be an integer.
If there is a positive integer $a$ such that
$${d^a \over n^{a-1}}-{n\choose r} \left({t\over n} \right)^a\ge k,$$
then $\Gn$ contains a subset $U$ of at least $k$ vertices such that every
$r$
vertices in $U$ have at least $t$ common neighbors.
\end{lemma}

Note that in this lemma they do not claim that $U(k,r,t)$ is a subgraph.
Nevertheless, using this lemma they improve the constant $c=c_r^{k,t}$ in
 \eqref{eq:L11_UpperGeneral} from $O((t+1)^{1/r} k^{2-(2/r)})$ to $c\le
2^{-1+(2/t)}(t+1)^{1/r}k$.

As they mention at the end of  \cite{AlonKriveSudaDep},
 both proofs of Theorem~\ref{FureM11Th} give a bit more
  (and thus Corollary \ref{co:AhasDegrees=r} can be sharpened accordingly):
\beq{eq:L11_UpperGeneralPlus_r}
  \ext(n,U^{+r}(k,r,t))< c n^{2-(1/r)},
 \eeq
 where the graph $U^{+r}$ is obtained from $U^{+}(k,r,t)$
by replacing the vertex $w$ in Definition~\ref{univgr} by an independent
set of $r$ vertices
with the same neighbors, $x_1\dots, x_k$.

However, the method of Dependent random choice gives more.
Call a graph $L_h$ on $h$ vertices {\it $r$-degenerate} if it satisfies
the condition of Conjecture \ref{co:716}.
In other words, there is an ordering of its vertices
 $x_1, \dots, x_h$ such that for every $1\le i\le h$ the vertex $x_i$
  has at most $r$ neighbors $x_j$ with $j<i$.

\begin{theorem}[Alon, Krivelevich, and Sudakov
 \cite{AlonKriveSudaDep}]\label{th:AlonUpper_r_degenarate}

If $L$ is bipartite $r$-degenerate graph on $h$ vertices, then for every
$n\geq h$
$$\ext(n,L)\le h^{1/(2r)}n^{2-(1/4r)}.$$
\end{theorem}

Applying the above results with $r=2$, $t=1$ and $k=c\sqrt{n}$ to find a
$U(k,2,1)$ one immediately obtains the following.  Any graph on $n$ vertices
with $c_1n^2$ edges contains a $1$-subdivision of $K_k$ with $k=c_2\sqrt{n}$
for some positive $c_2$ depending on $c_1$.  This answers a question of Erd\H
os  \cite{Erd1979-17}.  The theorems of Bollob\'as and Thomason~ \cite{BolloThom} and Koml\'os and
Szemer\'edi  \cite{KomlosSzemB} also imply the existence of such a large
topological clique but their subgraph is not necessarily a $1$-subdivision.

Given any graph $L$, let $\overline{d}$ denote the
 $\max_{X\subseteq V(L)} \{2e_L(X)/ |X|\}$,
 the maximum {\it local average degree}.
Then $L$ is $\lfloor \overline{d}\rfloor$-degenerate.
Hence the upper bound of Theorem~\ref{th:AlonUpper_r_degenarate} and
the random method lower bound in \eqref{LowerExt} yield that

\begin{corollary}
 For every bipartite graph $L$,
\beq{eq:ApproxErdosConj}
  \Omega(n^{2-c})\leq  \ext(n,L) \le O(n^{2- (c/8)}),
  \eeq
where $c=2/{\overline {d}}$, is the same as in \eqref{LowerExt}.
\end{corollary}

\section{Eigenvalues and extremal problems}\label{EigenValuesS}

Let $A=A(\Gn)$ be the adjacency matrix of $\Gn$, and $\jj$ be the vector each
entry of which is 1. Since \beq{EdgeCount} e(\Gn)=\half \jj A \jj^T \eeq and,
more generally, \beq{WalkCount} w_k(\Gn)=\half \jj A^k \jj^T \eeq counts the
number of $k$-edge walks in $\Gn$, therefore it is not so surprising that
eigenvalues can be used in extremal graph problems. An easy to read source on
spectra of graphs is Cvetkovi\v{c}--Doob--Sachs  \cite{CvetkoDoobSachs}.

\begin{theorem}[Babai--Guiduli  \cite{BabaiGuidulli}]
\label{BabaiGuid}

Let $\Lambda(G)=\max |\lambda_i|$, where $\lambda_1,\dots,\lambda_n$ are the eigenvalues of $A(\Gn)$.
If $K_{a,b}\not\subseteq \Gn$, (and $2\le a\le b$) then
\beq{BabaiGui}
\Lambda\le \root a\of {b-1}\cdot n^{1-(1/a)}+o(n^{1-(1/a)}).
\eeq
\end{theorem}

Since trivially
\beq{eq:sajatertek}
  2e(\Gn)\le \Lambda n\eeq
the inequality \eqref{BabaiGui} implies Theorem \ref{KovSosTurThB} apart from the $o()$ term.

\begin{remark}

For regular or almost regular graphs $\Lambda(\Gn)\approx
{2e(\Gn)\over n}$, and then the two estimates are basically
equivalent. The constant in the above theorem is not sharp since -- as
we know from Theorem \ref{FureImprove},-- the constant can be
improved.
\end{remark}

We have already mentioned Nikiforov's result (Theorem \ref{NikifZaraTh}) on the Zarankiewicz problem.
In fact, he proved  \cite{NikifZaran} that for all $n\ge b\ge a\ge 2$ and a $K_{a,b}$-free graph $\Gn$
 we have
\beq{eq:KSTjavNikiforov}
 \Lambda(G_n)\le (b-a+1)^{1/a}n^{1-(1/a)}+(a-1)n^{1-(2/a)}+(a-2). \eeq
This improves the coefficient in Theorem \ref{BabaiGuid}.
It also implies F\"uredi's bound \eqref{eq:KSTjav} for the
 $\ext(n, K_{a,b})$ according to \eqref{eq:sajatertek}.
For $C_4$-free graphs he has $\Lambda^2-\Lambda+1\leq n$.

Recall that $T_{n,k}$ denotes the Tur\'an graph, the $k$-partite
 complete graph of maximum size.
Given a $K_{k+1}$-free graph $\Gn$ Nikiforov  \cite{NikifEigenII}  showed that
 $\Lambda(\Gn)<\lambda(T_{n,k})$ unless $G=T_{n,k}$.
For a recent reference of a generalization see Z. L. Nagy  \cite{NagyZL}.

\section{Excluding topological subdivisions}\label{InfiFami}

\subsection{Large topological subgraphs}

We have already mentioned that our classification does not hold for
infinite families of excluded subgraphs. One important phenomenon
is that $\ext(n,\LL)$ can be linear for infinite $\LL$ even
if $\LL$ contains only cycles.
\footnote{Here the simplest case is Theorems \ref{ErdGallaiCyc}.}
Here we consider a very central graph theoretical problem strongly
connected to the 4-colour conjecture.

\begin{definition}

Given a graph $H$, its {\emR subdivision} (or a
 topological $H$) is obtained from
it by replacing each edge $e$ of $H$ by some paths $P_e$ so that these
paths do not have their inner (new) vertices in common.

\end{definition}

Wagner asked if for any integer $\ell$ there exists a $k=k_\ell$ such that any
$G$ with chromatic number $\chi(G)>k_\ell$ must contain a topological subdivision of $K_\ell$.
This was proved by Gabor Dirac and H. Jung (independently).
Answering a question of Dirac, Mader proved the following important result.

\begin{theorem}[Mader,  \cite{Mader}]\label{MaderATh}

If $\Gn$ is an $n$-vertex graph, and
$$e(\Gn)\ge n(\ell-1)2^{{\ell-1\choose 2}-1},$$
then
$\Gn$ contains a subdivision of the complete $\ell$-graph.
\end{theorem}

This statement is stronger than the original Wagner conjecture , since a graph
with large chromatic number contains a subgraph with large minimum degree.
Mader, and independently, Erd\H{o}s and Hajnal conjectured that

\begin{conjecture}[Mader, Erd\H{o}s-Hajnal]\label{MadErdHaj}
There exists a constant $c>0$ such that if $e(\Gn)>c\ell^2n$, then $\Gn$ contains
a topological $K_\ell$.
\end{conjecture}

A slightly weaker form of this conjecture was proved by Koml\'os and
Szemer\'edi,  \cite{KomlosSzem}, then -- by a different method -- Bollob\'as
and Thomason  \cite{BolloThom} proved this conjecture and, almost
immediately after that, Koml\'os and Szemer\'edi  \cite{KomlosSzemB} proved
Conjecture \ref{MadErdHaj} as well.

\begin{theorem}[Bollob\'as--Thomason]

Every graph $\Gn$ of size at least $256\ell^2n$
edges contains a topological complete subgraph
of order $\ell$.
\end{theorem}

As to the small values of $\ell$, Dirac conjectured that for $n\ge
3$ every $\Gn$ with $e(\Gn)\ge 3n-5$ contains a topological $K_5$. This
improvement of the famous Kuratowski theorem was proved by Mader in
 \cite{MaderCCA1998} and the corresponding extremal graphs were characterized
in  \cite{MaderCCA2005}.
The reader is recommended the excellent ``featured review'' of Carsten
Thomassen on the paper of Mader  \cite{MaderCCA1998}, on the MathSciNet.

An excellent survey of Mader on this topic is  \cite{MaderStirin}.

\subsection{Tur\'an numbers of subdivided graphs}

Let $\varepsilon$ be a positive real, $0<\varepsilon<1$.
Kostochka and Pyber  \cite{KostoPyber} proved that every $n$-vertex graph
$\Gn$ with at least
 $4^{t^2} n^{1+\varepsilon}$ edges contains a subdivision of $K_t$
  on at most $(7t^2\ln t)/\varepsilon$ vertices, where  $0<\varepsilon<1$.
This (for $t=5$) answers a question of Erd\H os
 about finding a non-planar subgraph of size $c(\varepsilon)$ in a graph
 with $n^{1+\varepsilon}$ edges.

Recently, T. Jiang  \cite{Jiang1} improved the Kostochka-Pyber
  upper bound to $O(t^2/\varepsilon)$.
On the other hand, for each  $0<\varepsilon<1$ and $n> n_0(\varepsilon)$
  there are $n$-vertex graphs of girth at least $1/\varepsilon$
(see Corollary~\ref{CycleRandCons}).
In such a graph any subdivision of $K_t$ must contain
$\Omega(t^2/\varepsilon)$
 vertices, so Jiang's result is sharp.


\begin{theorem}[Jiang and Seiver  \cite{JiangSeiver}]\label{JiangSeiverTh}
Let $L$ be a subdivision of another graph $H$.
For each edge $xy\in E(H)$ let $\ell(x,y)$ denote the
  length of the path in $L$ replacing the edge $xy$.
Suppose that $\ell(x,y)$ is even for each edge of $H$, and let
$\min\{\ell(x,y) : xy\in E(H)\} = 2m$.
Then $\ext(n,L) = O(n^{1+(8/m)})$.
  \end{theorem}

The main tools in the proof are the Dependent Random Choice,
Lemma~\ref{DependentRandom},
 and the Erd\H os-Simonovits $\Delta$-almost-regularization,
Theorem~\ref{DeltaReg}.

\section{Hypergraph Extremal Problems}\label{DegenHyperSec}

\subsection{Positive Density problems}

This is a short detour into Hypergraph Extremal Problems. Now our
``Universe'' is the class of $r$-uniform hypergraphs. Katona, Nemetz
and Simonovits  \cite{KatNemSim} showed (using a simple averaging) that

\begin{theorem}[Katona, Nemetz and Simonovits  \cite{KatNemSim}]

$\ext_r(n,\LL)/{n\choose r}$
is monotone decreasing, and therefore convergent.
\end{theorem}

\balAbra{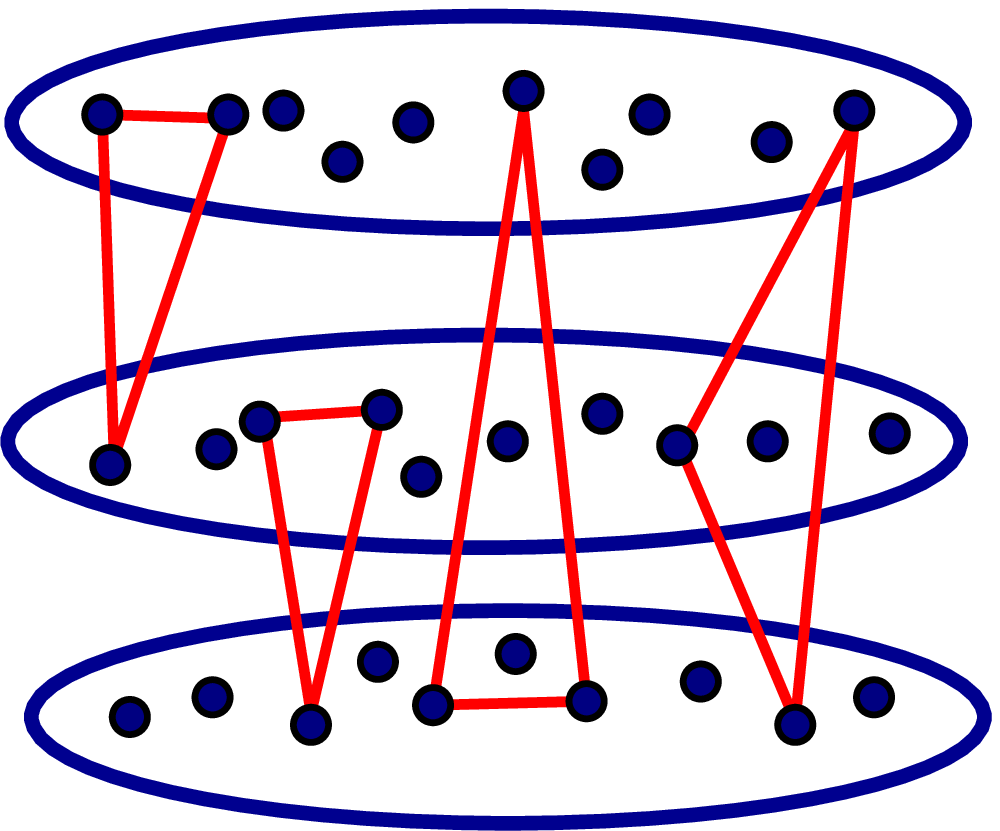}
The hypergraph extremal problems are extremely hard. Even the simplest
extension of Tur\'an's theorem is unsolved: Let $K_4\hyp 3$ be the
three-uniform hypergraph with 4 vertices and 4 triples.

\begin{construction}[Tur\'an, the simplest case]
\label{TurHyp4}
 We partition
$n$ vertices into three classes $C_1,C_2,C_3$
and we take all the triples of the form
$(x,y,z)$, where \dori (a) $x\in C_i,y\in C_i, z\in C_{i+1}$ (where the
indices are taken {\rm mod} 3);
 \dori (b) the three vertices are in three different groups.
\end{construction}

One can easily see that his construction contains no $K_4\hyp3$.

\begin{conjecture}[Tur\'an]
 Construction \ref{TurHyp4}
is asymptotically extremal for
$K_4\hyp3$. (Perhaps it is extremal, not only asymptotically
extremal, at least for $n>n_0$.)
\end{conjecture}

Here we cut it short and recommend the reader (among others) the survey
of F\"uredi on hypergraph extremal problems  \cite{FureLondon}, and also
the papers of F\"uredi and Simonovits \cite{FureSimFano},
Keevash and Sudakov  \cite{KeevSud}, F\"uredi-Pikhurko-Simonovits  \cite{FurePikhSim},
and the survey of Keevash  \cite{KeevashHyperSurv}.

\subsection{Degenerate hypergraph problems}

For $r$-uniform hypergraphs the $r$-partite graphs generalize the bipartite
graphs. An important illustration of this is the one below, extending
Theorem \ref{DegenCriter}.

\begin{theorem}
[Degenerated hypergraph problems]\label{th:HyperDegen}  For an
$r$-uniform extremal hypergraph problem of $\LL\hyp{r}$,
$\ext(n,\LL\hyp r)=o(n^r),$ if and only if there is an
$L\in\LL\hyp{r}$ which can be $r$-vertex-colored so that each
hyperedge of $L$ gets $r$ distinct colors.
 \end{theorem}

Theorem \ref{th:HyperDegen} is an easy corollary of the following
theorem of Erd\H{o}s, (which generalizes Theorem \ref{KovSosTurThB}).

\begin{theorem}[Erd\H{o}s  \cite{ErdHyper34}]
\label{th:ErdosHyper}

Let $K\hyp r(a_1, ..., a_r)$ be the $r$-uniform hypergraph with $r$
vertex-classes $C_1,\dots,C_r$, where $|C_i|=a_i$, and $a_1=t$.
Then $$\ext\hyp r(n,K\hyp r(a_1, ..., a_r))=O(n^{r-(1/t^{r-1})}).$$
 \end{theorem}

Extending some problems and results for ordinary graphs,
Brown, Erd\H{o}s and S\'os started investigating the following

\begin{problem}[Brown, Erd\H{o}s, and S\'os  \cite{BrownErdSosA},  \cite{BrownErdSosB}] \label{BrownErdSosPrb}
 Consider $r$-uniform
 hypergraphs for some fixed $r$, and denote by $\HH_{k,\ell}^r$ the
 family of $r$-uniform $k$-vertex hypergraphs with $\ell$ hyperedges.
Determine or estimate
$f_r(n,k,\ell):=\ext(n,\HH_{k,\ell}^r).$
\end{problem}

Brown, Erd\H{o}s, and S\'os proved many upper and lower bounds for
special cases of Problem \ref{BrownErdSosPrb}. We have already
mentioned one of them: the $f_3(n,6,3)$-problem.\footnote{If $r=3$, then
 we delete the subscript in $f_3$.} It is easy to see that
$f_3(n,6,3)<\reci 6 n^2$. The real question was if $f_3(n,6,3)=o(n^2)$ or not.
Ruzsa and Szemer\'edi  \cite{RuzsaSzemer} proved that the answer is YES.
We formulated this in Theorem~\ref{RuzsaSzemTh}. This theorem became a
very important one. We originate, among others, the ``Removal Lemma''
from here.

We shall return to this problem in the section on applications.

\section{Supersaturated graphs}\label{SupersatS}

The theory of Supersaturated extremal problems is a very popular area today.
Here we shall restrict ourselves to the supersaturated extremal graph problems related to bipartite excluded graphs, just mention a few further references, like
Lov\'asz and Simonovits \cite{LovSimBirk}, Razborov \cite{RazborFlag}, Lov\'asz \cite{LovaszLimitBook}, Reiher \cite{Reiher}.

Given a graph $G$, denote by $N(G,F)$ the number of subgraphs of $G$
isomorphic to $F$. Here we have to be slightly cautious: if $F$ has
non-trivial automorphisms, then we can count isomorphisms or
isomorphic subgraphs, and the ratio of these two numbers equal to the
automorphism number.

A theorem which asserts that a graph $\Gn$ contains very many graphs $L$
from a family $\LL$ is called a {\bf theorem on supersaturated graphs}.
Such theorems are not only interesting in themselves, but also are often
useful in establishing other extremal results. At this point it is
worthwhile mentioning such a result for complete bipartite graphs, obtained
by Erd\H{o}s and Simonovits  \cite{ErdSimSuper}:

\begin{theorem}[Number of complete bipartite graphs]
For any integers $a$ and
$b$ there exists a constant $c_{a,b}>0$ such that
if $\Gn$ is a graph with $e$ edges, then
 $\Gn$ contains at least
$[c_{a,b}e^{ab}/n^{2ab-a-b}]$ copies of $K_{a,b}$ .
\end{theorem}

\begin{corollary}\label{CycleSu}

Let $c>0$.
If $e(\Gn)=e>(1+c)\ext(n,C_4)$, then $\Gn$ contains at least
$\gamma e^4/n^4$ copies of $C_4$, for some $\gamma(c)>0$.
The random graph with $e$ edges shows that this is sharp.
 \end{corollary}

\Proof of the Cube Theorem (Sketch). Apply Theorem \ref{DeltaReg}
obtaining a $\Delta$-almost-regular (bipartite) $\tilde\Gn\subseteq\Gn$.
Apply the corollary to this $\ti\Gn$. It contains $\ga{e^4\over n^4}$
$C_4$'s. On the average, an edge of $\Gn$ is contained in $\ga
e^3/n^4$ copies of $C_4$. Take a typical edge $xy$: the bipartite
graph $G[U,V]$ spanned by the neighbors $U:=N(x)$ and $V:=N(y)$ -- by
$\ext(m,C_6)=O(m^{4/3})$,-- will contain a $C_6$. Now, $xy$ and this
$C_6$ will provide a $Q_8^+$: a cube with a diagonal.\Qed

Basically the same argument proves Theorem \ref{CubeRedu}.

\subsection{Erd\H{o}s-Simonovits-Sidorenko conjecture}

In this part $\chi(L)=2$. Erd\H{o}s and Simonovits  \cite{SimWat}
formulated three conjectures and also that the main idea behind these
conjectures is that the number of copies of subgraphs $L\in \Gn$ is
minimized by the random graph if $E=e(\Gn)$ is fixed and is not too
small.

To formulate these conjectures, first we calculate the ``expected
number of copies'' of $L\subseteq R_n$ if $R_n$ is a random graph with
edge probability $p=E/{n\choose 2}$. Let $v=v(L)$, and $e=e(L)$.
Clearly, if the edges are selected independently, with probability
$p$, then $K_n$ contains ${n\choose v}$ possible $v$-tuples, each
containing the same number $a_L$ of copies of $L$, and therefore
\beq{expectedL} \EE(\#(L\subseteq R_n))= (a_L+o(1)){n^v\over v!}p^e=
(a_L+o(1)){n^v\over v!}\left({2E\over n^2}\right)^e=a^*_L{E^e\over
 n^{2e-v}} \eeq

\def\ER{{\mathbb E\mathbb R}}

\begin{conjecture}[Erd\H{o}s--Simonovits,  \cite{SimWat}]

 There are two constants, $\Omega=\Omega_L>0$ and $c=c_L>0$ such that if
 $E>\Omega\cdot \ext(n,L)$, then any graph $\Gn$ with $E$ edges contains at
 least
$$c_L{E^e\over n^{2e-v}}$$
copies of $L$.
\end{conjecture}

This was the weakest form. The strongest form of this conjecture was

\begin{conjecture}[Erd\H{o}s and Simonovits]

For every $\eps>0$, if
 $E>(1+\eps)\cdot \ext(n,L)$, then any graph $\Gn$ with $E$ edges contains at
 least
$(1+\eps)\ER(n,L,E)$
copies of $L$, if $n>n_0(\eps)$, where $\ER(n,L,E)$ denotes the expected number
of edges of a random Erd\H{o}s-R\'enyi graph with $n$ vertices and $E$ edges.
\end{conjecture}

Obviously, one has to assume that $e(\Gn)>\ext(n,L)$.

\begin{remark}[Relation to Sidorenko's Conjecture]\label{ErdSimSidRel}
  At first sight Sido\-renko's conjecture  \cite{Sidorenko1993} seems to be
  sharper than the above one.  This is not the case. In fact, Sidorenko's
  Conjecture applies only to dense host graphs. There, as Sidorenko points
  out in his papers, the two versions are equivalent.
\end{remark}

Also, it is obvious that there is not much difference if we consider above
the hypergeometric model of random graphs, where the number of vertices and
edges are given, or if we fix only $n$ but the edges are taken independently,
and therefore $e(\Gn)$ follows a binomial distribution.

Sidorenko, working on
applications of extremal graph theorem to probability distribution
 translated the above conjecture to integrals and arrived at a conjecture
 \cite{Sidorenko1993}, where the error terms disappeared.
The meaning of his version was that if
one considers dense graphs and defines $L\subseteq R$ for the case when
$G$ is a function, generalizing the notion of graphs, then
the Random Continuous graph will have the least number of copies of
$L$, more precisely, that will minimize the corresponding integral.

We skip the formulation of this problem, just refer to some papers of
Lov\'asz, and Hatami  \cite{HatamiSidor},
and to the book of Lov\'asz  \cite{LovaszLimitBook}.

Jagger, \v{S}\'tov\'{\i}\v{c}ek, and Thomason
\cite{JaggerStovicekThomasonCCA1996} investigated the following
problem originating from a conjecture of Erd\H{o}s, disproved by
Thomason.


\begin{problem}
  Given a sample graph $L$, denote by $\rho_L(\Gn)$ the sum of copies of $L$
  in $\Gn$ and in its complementary graph. What is the minimum $\Gamma_n(L)$
  of this, taken over all $n$-vertex graphs?
\end{problem}

Erd\H{o}s conjectured that the random graph yields the minimum, for $K_4$. This
was disproved by Thomason  \cite{ThomasonDisproof}. Investigating the case of
general $L$, Jagger, \v{S}\'tov\'{\i}\v{c}ek, and Thomason proved some
interesting results in connection with Sidorenko's conjecture.

Here we should emphasize that there is a slight difference between
looking for copies of an $L$ in $\Gn$ or for copies of homomorphic images:
In the second case we allow vertices to map into $\Gn$ with some coincidences.

As to the Sidorenko Conjecture, the first unknown case (as Sidorenko mentions)
is when we delete the edges of a Hamiltonian cycle from $\KK5 5$.


\begin{theorem}[Conlon, Fox, Sudakov  \cite{ConFoxSudSid}
 \cite{ConFoxSudSidShort}]
The Sidorenko Conjecture holds if $L=L[A,B]$ is a bipartite graph with a vertex $x\in A$ completely joined to $B$.
\end{theorem}

\begin{remark}
 Unfortunately, we do not have sufficiently good lower bounds
for the extremal problem of the cube. The Erd\H{o}s--Simonovits
Conjecture was proved for $Q_8$ in  \cite{ErdSimCubeWat}.

Hatami proved the Sidorenko conjecture for any cube (i.e. of any
dimension), yet, that was not really enough to provide a reasonable
upper bound for the 4-dimensional cube. This reflects some difference
between extremal problems and the corresponding Supersaturated Graph
Problems (at least, for dense host graphs).
\end{remark}

\section{Ordered structures}

\def\Lol{{\overline{L}}}
\def\extdir{{\rm ext_{\rm dir}}}

\subsection{Directed graphs, ordered graphs}

\def\Lol{\overrightarrow{L}}
\def\onyil#1{\overrightarrow{#1}}
\def\extdir{\overrightarrow{{\bf ext}}}

There is an extensive literature on Digraph extremal problems, see e.g.,
the survey of Brown and Simonovits  \cite{BrownSim99}, or  \cite{BrownSimDM}.
We skip here the general theory.

Denote $\extdir(n,\Lol)$ the maximum number of edges in a directed graph
 not containing the oriented subgraph $\Lol$.
For every $\Lol$ containing a directed path of length 2 one has
  $\extdir(\Lol)\geq \lfloor n^2/4\rfloor$.
Indeed, orient the edges of $K_{n/2, n/2}$
 simultaneously into one direction.
For bipartite $\Lol$ it is more interesting to consider the minimum outdegree.

Consider the following directed graph $\Lol_{1,a,b}$ on $1+a+b$ vertices
 $w$, $x_1,\dots x_a$, $y_1, \dots , y_b$.
The oriented edges are $w$ to $x_i$ and $x_i$ to $y_j$ ($1\le i\le a$, $1\le j\le b$).


\begin{theorem}[Erd\H{o}s, Harcos and Pach  \cite{ErdHarcPach}]
  Given integers $a$ and $b$, there exists a $c=c_{a,b}>0$ such that the
  following holds.  Any oriented graph with minimum out-degree $\delta^+\ge c
  n^{1-(1/a)}$ contains a copy of $\Lol_{1,a,b}$.
 \end{theorem}

This result opened up a new interesting field with many open problems.

Another ordered Tur\'an function was defined by  Timmons  \cite{Timmons}.
He showed that if a graph with vertex set $\{ 1,2,\dots, n\}$
 has at least
 $$(1+o(1))(2/3)n^{3/2}$$
 edges, then it contains a $C_4$ with vertices $a_1b_1a_2b_2$
 such that $a_1,a_2< b_1,b_2$. 
He extended other ordinary Tur\'an problems to these zig-zag type questions.
Many problems remain unsolved.

\subsection{Erd\H os-Moser conjecture on unit distances}\label{subs:ErdosMoser}

Erd\H{o}s and Leo Moser  \cite{ErdMoser} conjectured that

\begin{conjecture}

If $n$ points of the plane are in convex position, then
the number of unit distances among them is $O(n)$.
\end{conjecture}

F\"uredi proved a slightly weaker result:

\begin{theorem}[F\"uredi  \cite{FurConvGeom}]\label{FureConvUTh}
 If $n$ points are in convex position in the plane, then
there are at most $O(n\log n)$ unit distances among them.
\end{theorem}

To prove this, F\"uredi directly
 formulated the excluded Ordered Matrix Property and
 solved a matrix-containment problem.
The crucial point of his proof was Theorem \ref{le:FurediNlogN} below.

The best known lower bound in the Erd\H os-Moser problem, $2n-7$,
 is due to Edelsbrunner and P. Hajnal~ \cite{EdelPHajnal}.

\subsection{Ordered submatrices}

The ordered matrix problems partly came from geometric problems
 (see Bienstock and Gy\H{o}ri,  \cite{BienGyori6}, F\"uredi~ \cite{FurConvGeom}),
 but they are interesting on their own, too.
A geometric application, called Erd\H os-Moser conjecture, is
 discussed above in Subsection \ref{subs:ErdosMoser}.

We have already indicated that most
 extremal graph problems have matrix forms, too:
To determine $\ext^*(m,n,L)$ we considered all $m\times n$ $0$-$1$ matrices not
containing any permutation of the bipartite adjacency matrix of $L$.

In the ordered
case here we exclud only those submatrices where the indexing of the rows
and columns of $\bM$ is fixed.
This way we exclude fewer subconfigurations.

\begin{definition}[Matrix containment]

Let $\bM$ and $\bP$ be two 0-1 matrices. We say that $\bM$
contains $\bP$ if we can delete some rows and columns of $\bM$ and
then perhaps switch some 1's into 0 so that the resulting matrix be
$\bP$. Otherwise we say that $\bM$ avoids $\bP$.
\end{definition}

So, we can delete rows and columns of $\bM$
but can not permute them.
Now we can define the Matrix Extremal Problems:

\begin{problem}[Ordered Matrix Problem]
\label{OrdMatProb}
 Given an $a\times b$ 0-1 (sample) matrix $\bP$, and a (huge)
$m\times n$ 0-1 matrix $\bM$, how many 1's can occur in $\bM$ under the
condition that $\bM$ does not contain $\bP$ in the ``ordered'' way.
Denote by $\exto(m,n,\bP)$ the maximum.
\end{problem}

One of the first nontrivial results was

\begin{theorem}[F\"uredi  \cite{FurConvGeom}]\label{le:FurediNlogN}

Let $$\bP=\begin{pmatrix}1&1&0\\1&0&1\end{pmatrix}.$$
If\, the $n\times n$ $0$-$1$ matrix $\bM$ does not contain $\bP$, then it has
at most $O(n\log n)$ $1$'s. In fact,
$\exto(n,n,\bP)=\Theta(n\log n)$.
\end{theorem}

Tardos  \cite{Tardos2005} proved that
$\exto(n,\bP)=n \log_2 n +O(n).$

Completing earlier works of F\"uredi and P\'eter Hajnal  \cite{FureHajn} Tardos
 \cite{Tardos2005} classified the ordered matrix Tur\'an numbers for all small
submatrices.  The extremely slow growing inverse Ackermann function is denoted
by $\alpha(n)$.

\begin{theorem}[\cite{FureHajn},  \cite{Tardos2005}]
If $\bP$ is a $0$-$1$ matrix with at most four $1$'s, then
$$
\exto(n,n,\bP)=
\begin{cases} 0 & or\\
\Theta(n),& or\\ \Theta(n\alpha(n)),& or\\ \Theta(n\log n),& or \\
\Theta(n^{3/2}).
\end{cases}$$
  \end{theorem}

\subsection{Ordered matrices and the Stanley-Wilf conjecture on subpermutations}\label{MatrixTardos}

Trying to prove the Erd\H{o}s-Moser Conjecture, F\"uredi and P\'eter
Hajnal  \cite{FureHajn} arrived at the following conjecture, proved by
 Marcus and Tardos.

\begin{theorem}[F\"uredi--Hajnal conj.  \cite{FureHajn}/Marcus-Tardos theorem
 \cite{MarcusTardos}]\label{FureHajnConj}
 ${}$
For all permutation matrices $\bP$ we have $\exto(n,n,\bP)=O(n)$.
\end{theorem}

This time there was a famous Stanley-Wilf conjecture around, on the
number of permutations ``avoiding'' a fixed permutation.
To formulate it, we need to define the Permutation containment:

\begin{definition}[Permutation containment]
We say that a permutation $\sigma: [1,n]\to [1,n]$ contains a permutation $\pi:[1,k]\to[1,k]$, if there exist $1\le x_1<x_2<\dots<x_k\le n$ for which
$$\sigma(x_i)<\sigma(x_j)\Text{if and only if} \pi(i)<\pi(j) .$$
\end{definition}

The famous Stanley--Wilf conjecture\footnote{Marcus and Tardos
  \cite{MarcusTardos} write that it is difficult to locate the
 corresponding reference.} states that

\begin{conjecture}[Stanley--Wilf]

For any permutation pattern $q$,
if $S_n(q)$ is the number of permutations of
 length $n$ avoiding the pattern $q$, then
there is a constant $c_q$ so
 that $S_n(q)\le c_q^n$.
\end{conjecture}

\begin{theorem}[Klazar  \cite{KlazarStW}]
The F\"uredi--Hajnal conjecture implies the Stanley-Wilf conjecture.
\end{theorem}

So Marcus and G. Tardos settled this conjecture as well.

\begin{remark}
The permutation containment is just a subcase of the more general question.
In some other cases there are
definite differences between ordinary Tur\'an type extremal problems
and the ordered matrix problems. For a special matrix, where the
corresponding graph is a tree, hence it has linear Tur\'an function,
our threshold function turns out to be $\Theta(n \log n)$.
\end{remark}

\section{Applications in Geometry}\label{AppliGeom}

\subsection{Applicability of the K\H{o}v\'ari-T. S\'os-Tur\'an bound}\label{GeomApplica}

We have mentioned that Theorem \ref{KovSosTurThB} is applicable in
several cases. Here we mention only two.

(A) The Unit Distance Graph of the Plane contains no $K(2,3)$.
Erd\H{o}s  used this to estimate the number of
unit distances by $O(n^{3/2})$.

(B) G. Megyesi and Endre Szab\'o\footnote{We use the longer versions of the
 names whenever we see chances to mix up authors of similar names.}
answered a question of F. E. P. Hirzebruch using this theorem.

Assume that we are given $k$ smooth curves in the the Complex
Projective Plane and assume that their union has only nodes and
tacnodes\footnote{Tacnode means roughly that the curve is touching
 itself.} as
singularities. Let $t(k)$ denote the maximum number of tacnodes in
such cases. Hirzebruch proved that $t(k)\le{4\over 9}k^2+{4\over 4}k$.
Hirzebruch asked if $\lim\sup t(k)/k>0 $.

\begin{theorem}[G. Megyesi, and E. Szab\'o  \cite{MegyeSzaboTac}]

There exist three positive constants, $A,B$ and $C$ for which
$$Ak^{1+(B/\log\log k)}\le t(k)\le Ck^{2-(1/7633)}.$$
\end{theorem}

\subsection{Unit Distances}

Erd\H{o}s was interested in the following problem:

\begin{problem}[Unit distances]

Given an $n$-element set in the $d$-dimensional Euclidean space
$\EE^d$, how many of the distances can be the same, say equal to 1?
\end{problem}

\begin{conjecture}[Unit distances]\label{UnitDistConj}

For any $\eps>0$, there is an $n_0$ such that if $n>n_0$ and
given an $n$-element set in the
plane
$\EE^2$, then the number of unit distances is at most $n^{1+\eps}$.
\end{conjecture}

The motivation of this conjecture is -- as Erd\H{o}s observed -- that
if we arrange the $n=k\times k$ points into a $k\times k$ grid, and rescale
this grid so that the ``most popular'' distance be 1, then
this distance will occur at most $n^{1+\eps}$ times, (actually,
approximately $n^{1+(c/(\log\log n))}$ times). So Erd\H{o}s conjectured
that the number of unit distances is in the plane has an upper bound
of roughly this form.

The first upper bound was a trivial application of Theorem \ref{KovSosTurThB}:

\begin{theorem}
[Unit distances, Erd\H{o}s 1946]
\label{th:Unit}

Given $n$ points in the plane, the number of unit distances among them
is at most
$$\ext(n,\KK23 )<\left({1\over \sqrt{2}}+o(1)\right)n^{3/2}.$$
In $\EE^3$ the number of unit distances is at most
\beq{U3dist}\ext(n,\KK33 )<c_{3,3}n^{5/3}.\eeq
 \end{theorem}

\proof. Since two circles intersect in at most 2 points, the Unit
Distance Graph of $\EE^2$ contains no $\KK23 $. This implies the first
inequality. Since 3 unit balls intersect in at most 2 points, the
Unit Distance Graph of $\EE^3$ does not contain any $\KK33 $. This
implies \eqref{U3dist}.\Qed

\begin{remarks}
(a) Everything is different for the higher dimensions: $\EE^4$
contains two orthogonal circles of radii $\reci \sqrt{2} $, and these
form a $K(\infty,\infty)$ in the corresponding unit graphs of $\EE^d$,
for $d\ge 4$. (This is the so called Lenz Construction.)
(See also Section \ref{NormGraphS}.)

(b) How sharp is this application? As the reader can see, it is very
far from the conjectured upper bound. However, just to improve it to
$o(n^{3/2})$ is non-trivial (J\'ozsa-Szemer\'edi
 \cite{JozsaSzem}). Actually, for the plane an $O(n^{4/3})$ upper bound
was proved by Beck and Spencer  \cite{BeckSpencUnit} and Spencer, Szemer\'edi
and Trotter  \cite{SpencSzemTrot},
which is sharp if we do not insist on Euclidean metric, only on ``normed spaces''.
For this see the results of Peter
Brass  \cite{PBrassUnit} and Pavel Valtr  \cite{ValtrUnit}.
\end{remarks}

\subsection{Cells in line arrangements}
\def\cI{{\mathcal I}}

Let $\cI(m,n)$ denote the maximum number of edges in $m$ distinct cells
 determined by an arrangement of $n$ lines in the plane.
Canham~ \cite{Canham} showed that for an absolute constant $c>0$
\beq{eq:Canham}
  \cI(m,n) < c(m \sqrt{n} + n).
 \eeq
Indeed, if we construct a bipartite graph where one side of the
 vertex set consists of the $m$ cells (or any other family of
  $m$ convex sets with disjoint interiors), the other side of the vertex set
  consists of the $n$ (tangent) lines and two vertices are joined if the
  corresponding geometric objects are incident, then
  it is easy to see that this graph does not contain a $K_{5,2}$. \qed

This was a first nontrivial step toward the determination of the
 exact order of the magnitude of $\cI(m,n)$
by Clarkson, Edelsbrunner,  Guibas,  Sharir, and Welzl~ \cite{CEGSW};
 it is $\Theta(n^{2/3}m^{2/3}+n)$.
More about this and other geometric applications see the monograph of
 Pach and Agarwal~ \cite{PachAgarwal}.

\section{Further connections and problems}

\subsection{Connections of hypergraphs and critical graphs}\label{CriticalGrS}

We discussed Degenerate Hypergraph Extremal Problems in Section
\ref{DegenHyperSec}. Here we continue that line.

\subsub Excluding the 3-uniform hypergraph cones//

Many of the other results, problems of  \cite{BrownErdSosA} were also
degenerate ones. One of them was where $\TT$ is the family of
triangulations of the 3-dimensional sphere. This problem gave the name
to this paper  \cite{BrownErdSosA}. The crucial point was excluding
the double cones:

\begin{definition}[$r$-cones]
The vertices of the 3-uniform hypergraph $Q_{r,t}$
are $x_1,x_2,\dots,x_r$, and
 $y_1,y_2,\dots,y_t$ for some $t$,
and the hyperedges are $x_iy_jy_{j+1}$,
for all the possible $i,j$, where $y_{t+1}=y_1$. Further,
$\QQ_r:=\{Q_{r,t}~:~t=3,4,5,\dots\}$.
\end{definition}

\begin{theorem}[Brown, Erd\H{o}s, S\'os, $r=2$,  \cite{BrownErdSosA},
 Simonovits  \cite{SimCrit} $r\ge 2$]
$$\ext(n,\QQ_r):=O(n^{3-(1/r)}).$$
\end{theorem}

For $r=2,3$ there are matching lower bounds here. Actually, for $r=2$
Brown, Erd\H{o}s and S\'os gave a construction, where not only the
double-cone was excluded, but all the triangulations of the sphere. In
Simonovits' lower bound only the double cone was considered.

In  \cite{SimColumbus} Simonovits returned to this question and --
using the main idea of Brown's construction  \cite{BrownThom} -- he
proved

\begin{theorem}[Simonovits  \cite{SimColumbus}]

There are (finite geometric) 3-uniform hypergraphs without triple-cones
(i.e. without hypergraphs from $\QQ_3$)
and still having at least $cn^{3-(1/3)}$ triples.
\end{theorem}

We saw that for the family of
triangulations of the sphere, and for the family of Double Cones the
extremal number is $O(n^{3-(1/2)})$  \cite{BrownErdSosA}, (see
 \cite{SimCrit}).

Brown, Erd\H{o}s and T. S\'os arrived at their question (most probably) since
they wanted to generalize certain results from ordinary graphs to
hypergraphs. Simonovits came from a completely different direction: he used
this to disproved a conjecture of Gallai on independent vertices in
4-colour-critical graphs.

$G$ is colour-edge-critical, if deleting any edge of $G$, we get a
$(\chi(G)-1)$-chromatic graph.  The 3-colour-critical graphs are the odd
cycles, so the problem of critical graphs becomes interesting for the
4-chromatic case. Here we shall restrict ourselves to this case and
suggest the reader to read Bjarne Toft's results on this topic in general.

Erd\H{o}s asked if a 4-colour-critical graph can have $cn^2$ edges and
Bjarne Toft constructed such a 4-chromatic graph  \cite{ToftCrit} of
$\approx{n^2\over 16}$ edges.
This and some related questions can also be
found in Lov\'asz' book: Combinatorial Exercises  \cite{LovCombExerc}.

Gallai had many beautiful conjectures on 4-colour-critical graphs. One of them,
  however,
was ``completely demolished''. He conjectured that if $\Gn$ is
4-colour-critical, then $\alpha(\Gn)\le n/2$.
$G_{4m+2}=C_{2m+1}\Otimes C_{2m+1}$  is 6-critical, with $\mindeg(G_{4m+2})=2m+3$.
Simonovits -- ``blowing up'' the vertices in one of the two odd cycles, -- proved that there are 6-critical graphs $G_n$ with
$\alpha(\Gn)= n-o(n)$.

It turned out that slightly earlier Brown and Moon  \cite{BrownMoon}
already disproved Gallai's conjecture for the 4-chromatic case, with a ``clever but simple'' construction.

\begin{theorem}[Brown and Moon  \cite{BrownMoon}]
There exist 4-chromatic edge-critical graphs $\Gn$ with
$\alpha(\Gn)>n-c\sqrt{n}$, for some constant $c>0$.
\end{theorem}

Next, Bjarne Toft came up with his construction, mentioned above.
Using this and a hypergraph extremal theorem, Simonovits proved

\begin{theorem}
There exists a constant $c_2>0$ such that if $\Gn$ is 4-colour-critical, then
$\alpha(\Gn)\le n-c_2n^{2/5}.$
\end{theorem}

This was obtained as follows: Simonovits reduced the original problem
to estimating the number of independent vertices of degree 3 in a
4-colour-critical graph.  The neighborhoods of these vertices generated a
3-uniform hypergraph $\HHH3m$ on the remaining vertices.
Simonovits -- using the Sperner Lemma from Topology proved that
if $I$ is a set of independent vertices of degree 3, in $V(\Gn)$, then for $m:=n-|I|$,
$|I|<\ext_3(m,\QQ_2)=O(m^{5/2})$, see  \cite{SimCrit}.
He observed that
$\HHH3m$ cannot contain double cones. This proved that
 $|I|<n-cn^{2/5}$. \qed

(b) Lov\'asz observed that instead of excluding the graphs from
$\QQ_2$ one can exclude a larger family, $\ti\QQ$: those 3-uniform
hypergraphs which obey the conclusion of Sperner's lemma
 \cite{LovCrit}: each pair $(x,y)$ is contained in an even number of
hyperedges. This enabled him to completely settle {\emR this} Gallai
problem on colour-critical graphs.
He proved that $\ext(n,\ti\QQ_2)\le
{n\choose 2}$. So he obtained $|I|<n-c\sqrt{n}$, in Gallai's
problem. Besides proving and using a more applicable extremal graph
theorem he also generalized the Brown--Moon construction.

(c) It was an interesting feature of Lov\'asz' solution that to get
an upper bound on $\ext(n,\ti\QQ)$ he used linear algebra.

We finish this part by sketching the proof of Lov\'asz on the upper
bound. 

\begin{theorem}[Lov\'asz  \cite{LovCrit}]
Let $\EE^{(3)}$ denote the family of 3-uniform hypergraphs $H$ in
which each pair of vertices is contained in an even number of triplets
(i.e. hyperedges).  Then $\ext(n,\EE^{(3)})\le {n\choose 2}$.
\end{theorem}

\Proof (Sketch). Assume that $H^{(3)}_n$ contains no subgraphs from $\EE^{(3)}$.
Consider that vectorspace over $GF(2)$ of dimension ${n\choose 2}$ where
the coordinates are indexed by pairs from $1,\dots,n$.
Represent each triple of $H^{(3)}_n$ by such a vector,
where we have 1 in those coordinates which are pairs form our triple.
The condition that $H^{(3)}_n$ contains no subgraphs from $\EE^{(3)}$ translates into the fact,
that these vectors are linearly independent. Hence their number is at most the dimension of the vector-space.
\qed

Now, repeating the original argument of Simonovits, Lov\'asz obtained

\begin{theorem}
There exists a constant $c_3>0$ such that if $\Gn$ is 4-colour-critical, then
$\alpha(\Gn)\le n-c_3n^{1/2},$
\end{theorem}

This with the Brown-Moon construction completely settles Gallai's
original problem, providing a matching lower bound. Lov\'asz proved a
more general theorem, and extended the Brown-Moon construction as
well.
We close this part with a beautiful conjecture of Erd\H{o}s:

\begin{problem}
Is it true that if $(\Gn)$ is a sequence of 4-colour-critical graphs,
then $\mindeg(\Gn)=o(n)$?
\end{problem}

(Simonovits  \cite{SimCrit} and Toft  \cite{ToftCrit} succeeded in constructing
$4$-color-critical graphs with minimum degrees around $c\root 3 \of n$.)

\SepaRef Several related results can be found in Lov\'asz  \cite{LovCombExerc}.

\subsection{A multiplicative Sidon problem and $C_{2k}$-free
graphs}\label{AppliNumb}

As it was explained in Subsection \ref{HistorA}, the Erd\H{o}s problem
about $\ext(n, C_4)$ in  \cite{ErdTomsk} was obtained from a
multiplicative Sidon type question.  He investigated subsets of
integers of $A\subset \{1,2,\dots, n\}$ with the property that for any
four members of $A$ the pairwise products are distinct, $a_ia_j\neq
a_ka_\ell$.

A. S\'ark\"ozy, P. Erd\H{o}s, and V. T. S\'os
 \cite{ErdSarkSos} started investigating the
more general problem.

\begin{problem}\label{GenerMultSidon}

Fix an integer $k$.
How many integers can we take from $[1,n]$ if the product of no $k$ of them
is a square.
\end{problem}

Interestingly, this Problem  also lead to
 Tur\'an type questions, namely to $\ext(m,n,C_{2k})$
 with $m\gg n$.
Their conjecture (Conjecture \ref{SarkoSosConj} above) was proved by
  Gy\H{o}ri  \cite{GyoriC6}, see Theorem~\ref{th:Gyori}.
We shall not go  into the number
theoretic details; just refer the reader again to 
 \cite{GyoriC6}.

\subsection{Cycle-free subgraphs of the $d$-dimensional hypercube}\label{S:Hypercube}

The {\em d-dimensional hypercube}, $Q^d$, is the graph
whose vertex set is $\{0,1\}^d$ and whose edge set is the set of pairs
that differ in exactly one coordinate, $e(Q^d)= d2^{d-1}$.
Let $\gamma(C_{\ell})=\lim_{d\to \infty} \ext(Q^d,C_4)/e(Q^d)$.
Note that $\gamma(C_\ell)$ exists, because $\ext(Q^d,C_4)/e(Q^d)$ is a non-increasing and
 bounded function of $d$.
Considering the edges between the levels $2i$ to $2i+1$ one can see that
 $\ext(Q^d,C_4)\ge (1/2)e(Q^d)$.
The following conjecture is still open.
\begin{conjecture}[Erd{\H o}s \cite{ErdosCubeConj}]
$ \ext(Q^d,C_4)= \left(\frac{1}{2}+o(1)\right)e(Q^d).$
\end{conjecture}
The best upper bound $\gamma(C_4) \leq 0.6226$  was obtained by
 Thomason and Wagner~\cite{ThomasonW},
slightly improving the result of Chung~\cite{ChungHypercube}.

Erd{\H o}s~\cite{ErdosCubeConj} also asked whether $\ext(Q^d,C_{2k})$ is $o(d)2^d$
 for $k>2$.
This was answered negatively for $C_6$ by Chung~\cite{ChungHypercube},
 showing that $\gamma(C_6)\ge 1/4$.
The best known results for $C_6$ are $1/3 \leq \gamma(C_6) < 0.3941$ due to
Conder~\cite{ConderHypercube} and Lu~\cite{LuCube}, respectively.

On the other hand, for every $t\geq 2$ the inequalities
\begin{equation}\label{c4t}
\ext(Q^d,C_{4t}) \leq O(d^{\frac{1}{2}-\frac{1}{2t}}2^d)\quad \text{and}\quad
\ext(Q^d,C_{4t+6})= O(d^{\frac{15}{16}-\frac{1}{16t}}2^d)
\end{equation}
were proved by Chung \cite{ChungHypercube} and  F\"uredi and \"Ozkahya~\cite{FurOzkahya}, respectively.
Hence $\gamma(C_{2k})=0$, except $\gamma(C_4)\geq 1/2$, $\gamma(C_6)\geq 1/3$ and
 the problem of deciding wether $\gamma(C_{10})=0$  is still open.

Conlon~\cite{ConlonHypercube} generalized \eqref{c4t} by showing
$\ext(Q^d,H)=o(e(Q^d))$ for all $H$ that admit a $k$-{\it partite representation},
 also satisfied by each $H=C_{2k}$ except for $k\in\{2,3,5\}$.

\subsection{Two problems of Erd\H{o}s}

 Of course, we should close with two open problem of Erd\H{o}s. The first one is the general version of that problem which was solved in \cite{GrzesikC5} and \cite{HataRazboC5C3}, see Section \ref{CentralExamplesS}.

\begin{conjecture}[Erd\H{o}s  \cite{ErdosTwoProblems}] \label{ErdosOddCycle1}
 Suppose that $G$ is a graph on $(2k + 1)n$ vertices and of odd girth
 $2k + 1$. Then $G$ contains at most $n^{2k+1}$ induced cycles of
 length $2k + 1$.
\end{conjecture}

 The next conjecture is also very famous and is motivated by
  the blown up pentagon (if we restrict it to $k=2$.)

\begin{conjecture}[Erd\H{o}s  \cite{ErdosTwoProblems}] \label{ErdosOddCycle2}
 Suppose that $G$ is a graph on $(2k + 1)n$ vertices and of odd girth at least $2k + 1$.
Then $G$ can be made bipartite by omitting at most $n^2$ edges.
\end{conjecture}

For the best known results here, for $k=1$, see Erd\H{o}s, Faudree, Pach, and Spencer
\cite{ErdFauPachSpe} and
Erd\H{o}s, Gy\H{o}ri, and Simonovits  \cite{ErdGyoriSim}.

\section
{Acknowledgements}

The authors are greatly indebted for fruitful discussions and helps to
 a great number of colleagues, among others to
 R. Faudree, E. Gy\H ori, and Z. Nagy.



\begin{thebibliography}{999}

{\small  

\itemsep=0pt\parskip=0pt

\bibitem{AKSSzApprox}
M. Ajtai, J. Koml\'os, M. Simonovits, and E. Szemer\'edi:
On the approximative solution of the Erd\H{o}s-S\'os conjecture on trees,
   (manuscript).

\bibitem{AKSSzDenseBlock}
M. Ajtai, J. Koml\'os, M. Simonovits, and E. Szemer\'edi:
Some elementary lemmas on the Erd\H{o}s-T. S\'os conjecture for trees,
   (manuscript).

\bibitem{AKSSzSharp}
M. Ajtai, J. Koml\'os, M. Simonovits, and E. Szemer\'edi:
The solution of the Erd\H{o}s-S\'os conjecture for large trees,
   (manuscript, in preparation).



\bibitem{AKSLoebl95}
M. Ajtai, J. Koml\'os, and E. Szemer\'edi:
On a conjecture of Loebl,
in Graph theory, Combinatorics, and Algorithms, Vol. 1, 2 (Kalamazoo, MI, 1992),
Wiley-Intersci. Publ., pp. 1135--1146. Wiley, New York, 1995.

\bibitem{AllenKeevSudVers}
P. Allen, P. Keevash, B. Sudakov, and J. Verstra\"ete:
Tur\'an numbers of bipartite graphs plus an odd cycle, submitted. 

\bibitem{AlonExpand}
N. Alon:
Eigenvalues and expanders,
Combinatorica 6 (1983), 83--96.

\bibitem{AlonHB}
N. Alon:
Tools from higher algebra,
in : "Handbook of Combinatorics", R. L. Graham, M. Gr\"otschel and L. Lov\'asz, eds,
North Holland (1995), Chapter 32, pp. 1749--1783.


\bibitem{AHL}
N. Alon, S. Hoory, and N. Linial:
The Moore bound for irregular graphs,
Graphs Combin. 18 (2002), no. 1, 53--57.



\bibitem{AlonKriveSudaDep}
N. Alon, M. Krivelevich, and B. Sudakov:
Tur\'an numbers of bipartite graphs and related Ramsey-type questions,
Combin. Probab. Comput. 12 (2003), no. 5-6, 477--494.



\bibitem{AlonMilman}
N. Alon and V. D. Milman:
$\lambda_1$-isoperimetric inequalities for graphs and superconcentrators,
J. Combin. Theory Ser. B 38 (1985), 73--88.

\bibitem{AlonRonySzab}
N. Alon, L. R\'onyai, and T. Szab\'o:
Norm-graphs: variations and applications, J. Combin. Theory Ser. B 76 (1999), 280--290.

\bibitem{BabaiGuidulli}
L. Babai and B. Guiduli:
Spectral extrema for graphs: the Zarankiewicz problem,
Electronic J. Combin. 15 (2009), R123.

\bibitem{BaerPolar}
R. Baer:
Polarities in finite projective planes,
Bull. Amer. Math. Soc. 52 (1946), 77--93.

\bibitem{BGMV2007}
C. Balbuena, P. Garc\'{\i}a-V\'azquez, X. Marcote, and J. C. Valenzuela:
New results on the Zarankiewicz problem,
Discrete Math. 307 (2007), no. 17-18, 2322--2327.


\bibitem{BGMV2007Gyori}
C. Balbuena, P. Garc\'{\i}a-V\'azquez, X. Marcote, and J. C. Valenzuela:
Counterexample to a conjecture of Gy\H{o}ri on $C_{2l}$-free bipartite graphs,
Discrete Math. 307 (2007), no. 6, 748--749.

\bibitem{BGMV2008}
C. Balbuena, P. Garc\'{\i}a-V\'azquez, X. Marcote, and J. C. Valenzuela:
Extremal K(s,t)-free bipartite graphs,
Discrete Math. Theor. Comput. Sci. 10 (2008), no. 3, 35--48.


\bibitem{BaliBollRioSch}
P. N. Balister, B. Bollob\'as, O. M. Riordan, and R. H. Schelp:
Graphs with large maximum degree containing no odd cycles of a given length,
J. Combin. Theory B 87 (2003), 366--373.

\bibitem{BalisGyoriLehSchelp}
P. N. Balister, E. Gy\H{o}ri, J. Lehel, and R. H. Schelp:
Connected graphs without long paths,
Discrete Math 308 (2008), no. 19, 4487--4494.

\bibitem{BallPepe}
S. Ball and V. Pepe:
Asymptotic improvements to the lower bound of certain bipartite Tur\'an numbers,
Combin. Probab. Comput. 21 (2012), no. 3, 323--329.

\bibitem{BeckSpencUnit}
J. Beck and J. Spencer:
Unit distances,
J. Combin. Theory Ser. A 37 (1984), 231--238.

\bibitem{Behrend}
F. Behrend:
On sets of integers which contain no three terms in arithmetic progression,
Proc. Nat. Acad. Sci. US. 
32 (1956), 331--332. 

\bibitem{Benson}
C. T. Benson:
Minimal regular graphs of girths eight and twelve,
Canad. J. Math. 18 (1966), 1091--1094. 

\bibitem{BienGyori6}
D. Bienstock and E. Gy\H{o}ri:
An extremal problem on sparse $0$-$1$ matrices,
SIAM J. Discrete Math. 4 (1991), no. 1, 17--27.

\bibitem{BBK}
P. Blagojevi\'c, B. Bukh, and R. Karasev:
Tur\'an numbers for $K_{s,t}$-free graphs: topological obstructions and algebraic constructions,
arXiv:1108.5254v3, 3 Jun 2012. 



\bibitem{BolloCycleMod}
B. Bollob\'as:
Cycles modulo $k$,
Bull. London Math. Soc. 9 (1977), no. 1, 97--98.


\bibitem{BolloExtreBook}
B. Bollob\'as:
Extremal Graph Theory,
Academic Press, London, 1978.

\bibitem{BolloRGbook}
B. Bollob\'as:
Random Graphs, Academic Press, London, 1985.


\bibitem{BolloHB}  B. Bollob\'as: Extremal graph theory, in: R. L. Graham, M. Gr\"otschel, and L. Lov\'asz (Eds.), Handbook of Combinatorics, Elsevier Science, Amsterdam, 1995, pp. 1231--1292.




\bibitem{BolloThom}
B. Bollob\'as and A. Thomason:
Proof of a conjecture of Mader, Erd\H{o}s and Hajnal on topological subgraphs,
European J. Combin 19 (1998), 883--887.

\bibitem{BondyHB}
J. A. Bondy:
Basic graph theory: paths and circuits,
Handbook  of Combinatorics, Vol. I., pp. 3--110, Elsevier, Amsterdam, 1995.

\bibitem{BondyOnErdos}
J. A. Bondy:
Extremal problems of Paul Erd\H{o}s on circuits in graphs,
Paul Erd\H{o}s and his mathematics, II (Budapest, 1999), 135--156,
Bolyai Soc. Math. Stud., 11, J\'anos Bolyai Math. Soc., Budapest, 2002.

\bibitem{BondySim}
J. A. Bondy and M. Simonovits:
Cycles or even length in graphs,
J. Combin. Theory Ser. B 16 (1974), 97--105.

\bibitem{BondyVince}
J. A. Bondy and A. Vince:
Cycles in a graph whose lengths differ by one or two,
J. Graph Theory 27 (1998), 11--15.

\bibitem{BrandtDobson}
S. Brandt and E. Dobson:
The Erd\H{o}s-S\'os conjecture for graphs of girth $5$,
Selected papers in honour of Paul Erd\H{o}s on the occasion of his 80th
birthday (Keszthely, 1993), Discrete Math. 150 (1996), no. 1-3. 411--414.

\bibitem{PBrassUnit}
P. Brass:
Erd\H{o}s distance problems in normed spaces,
Comput. Geom. 6 (1996), no. 4, 195--214.

\bibitem{BrownThom}
W. G. Brown:
On graphs that do not contain a Thomsen graph,
Canad. Math. Bull. 9 (1966), 281--285.

\bibitem{BrownNonExist}
W. G. Brown:
On the non-existence of a type of regular graphs of girth 5,
Canad. J. Math. 19 (1967), 644--648.

\bibitem{BrownMoon}
W. G. Brown and J. W. Moon:
Sur les ensembles de sommets ind\'ependants dans les graphes chromatiques minimaux, (French),
Canad. J. Math. 21 (1969), 274--278.

\bibitem{BrownErdSosA}
W. G. Brown, P. Erd\H{o}s and V. T. S\'os:
On the existence of triangulated spheres in 3-graphs, and related problems,
Period Math. Hungar.  3 (1973), 221--228.

\bibitem{BrownErdSosB}
W. G. Brown, P. Erd\H{o}s and V. T. S\'os:
Some extremal problems on $r$--graphs,
New Directions in the Theory of Graphs (ed. F. Harary),
Academic Press, New York, 1973, pp. 53--63. 

\bibitem{BrownSimDM}
W. G. Brown and M. Simonovits:
Digraph extremal problems, hypergraph extremal problems, and the densities of graph structures,
Discrete Math. 48 (1984), no. 2-3, 147--162.

\bibitem{BrownSim99} W. G. Brown, and M. Simonovits:
Extremal multigraph and digraph problems,
Paul Erd\H{o}s and his mathematics, II (Budapest, 1999), pp. 157--203,
Bolyai Soc. Math. Stud., 11, J\'anos Bolyai Math. Soc., Budapest, 2002.

\bibitem{CV1991}
L. Caccetta and K. Vijayan:
Long cycles in subgraphs with prescribed minimum degree,
Discrete Math. 97 (1991), no. 1--3, 69--81.

\bibitem{CaenSzek}
D. de Caen and L. A. Sz\'ekely:
The maximum size of 4- and 6-cycle free bipartite graphs on $m$, $n$ Vertices,
Sets, Graphs and Numbers (Budapest, 1991),
Colloquium Mathematical Society J\'anos Bolyai, vol. 60, North-Holland, Amsterdam, 1992, pp. 135--142.

\bibitem{Canham}
R. Canham:
A theorem on arrangements of lines in the plane,
Israel J. Math. 7 (1969), 393--397.


\bibitem{ChungHypercube}
F. Chung:
Subgraphs of a hypercube containing no small even cycles,
J. Graph Theory 16 (1992), 273--286.


\bibitem{ChungGrahamErdLegacy}
F. R. K. Chung and R. L. Graham:
Erd\H{o}s on Graphs: His Legacy of Unsolved Problems,
A. K. Peters Ltd., Wellesley, MA, 1998.

\bibitem{ClapFlocShee}
C. R. J. Clapham, A. Flockart, and J. Sheehan:
Graphs without four-cycles,
J. Graph Theory 13 (1989), 29--47.

\bibitem{CEGSW}
K. Clarkson, H. Edelsbrunner, L. J. Guibas, M.  Sharir, and E. Welzl:
Combinatorial complexity bounds for arrangements of curves and  spheres,
Discrete Comput. Geom. 5 (1990), no. 2, 99--160.

\bibitem{ConderHypercube}
M. Conder:
Hexagon-free subgraphs of hypercubes,
J. Graph Theory 17 (1993), 477--479.

\bibitem{ConlonHypercube}
D. Conlon:
An extremal theorem in the hypercube,
Electron. J. Combin. 17 (2010), Research Paper 111.

\bibitem{ConFoxSudSid}
D. Conlon, J. Fox, and B. Sudakov:
An approximate version of Sidorenko's conjecture,
Geom. Funct. Anal. 20 (2010), no. 6, 1354--1366.

\bibitem{ConFoxSudSidShort}
D. Conlon, J. Fox, and B. Sudakov:
Sidorenko's conjecture for a class of graphs: an exposition,
 http://arxiv.org/abs/1209.0184

\bibitem{Cooley09} O. Cooley:
Proof of the Loebl-Koml\'os-S\'os conjecture for large, dense graphs,
Discrete Math. 309 (2009), no. 21, 6190--6228.

\bibitem{CulikZara}
K. \v Cul\'{\i}k:
Teilweise L\"osung eines verallgemeinerten Problems von K. Zarankiewicz,
Ann. Polon. Math. 3 (1956), 165--168.

\bibitem{CvetkoDoobSachs}
D. M. Cvetkovi\v{c}, M. Doob, and H. Sachs:
Spectra of Graphs,
Academic Press Inc., New York, 1980.

\bibitem{DamasdiHegerSzonyiZaran}
G. Dam\'asdi, T. H\'eger, and T. Sz\H onyi:
The Zarankiewicz problem, cages, and geometries,
manuscript 2013.

\bibitem{EdelPHajnal}
H. Edelsbrunner and P. Hajnal:
A lower bound on the number of unit distances between the vertices of a convex polygon,
J. Combin. Theory Ser. A 56 (1991), no. 2, 312--316.

\bibitem{Erdhomepage}Paul Erd\H{o}s:  Erd\H{o}s homepage (his scanned in papers up to 1989):\dori \hskip 2cm www.renyi.hu/\~{}p\_erdos.

\bibitem{ErdTomsk}
P. Erd\H{o}s:
On sequences of integers no one of which divides the product of two others, and some related problems,
Mitt. Forschungsinst. Math. u. Mech. Tomsk 2 (1938), 74--82.

\bibitem{ErdGraphProbA}
P. Erd\H{o}s:
Graph theory and probability I,
Canad. J. Math. 11 (1959), 34--38.

\bibitem{ErdGraphProbB}
P. Erd\H{o}s:
Graph theory and probability II,
Canad. J. Math. 13 (1961), 346--352.

\bibitem{ErdSmole}
P. Erd\H{o}s:
Extremal problems in graph theory,
Proc. Sympos. Smolenice, 1963, pp. 29--36, Publ. House Czechoslovak Acad. Sci., Prague, 1964.

\bibitem{ErdSmoleProb}
P. Erd\H{o}s:
Some applications of probability to graph theory and combinatorial problems,
Theory of Graphs and its Applications (Proc. Sympos. Smolenice, 1963), pp. 133--136,
 Publ. House Czech. Acad. Sci., Prague, 1964.

\bibitem{ErdHyper34}
P. Erd\H{o}s:
On some extremal problems in graph theory,
Israel J. Math.  3 (1965), 113--116.

\bibitem{ErdRome}
P. Erd\H{o}s:
Some recent results on extremal problems in graph theory,
Theory of Graphs (ed P. Rosenstiehl), (Internat. Sympos.,  Rome, 1966),
Gordon and Breach, New York, and Dunod, Paris, 1967,  pp. 117--123.

\bibitem{ErdTihany}
P. Erd\H{o}s:
On some new inequalities concerning extremal properties of graphs,
Theory of Graphs (P. Erd\H{o}s and G. Katona, Eds.), Academic Press, Nev. York, 1968, pp. 77--81.

\bibitem{ErdArt} P. Erd\H{o}s: The Art of Counting (ed. J. Spencer), The MIT Press, Cambridge, Mass., 1973.

\bibitem{ErdosKeszthely73}
P. Erd\H{o}s:
Problems and results on finite and infinite combinatorial analysis,
{\it in} Infinite and Finite Sets (Proc. Conf., Keszthely, Hungary, 1973), pp. 403--424,
Proc. Colloq. Math. Soc. J. Bolyai 10, Bolyai--North-Holland, 1975.

\bibitem{Erd1975-42}
P. Erd\H{o}s:
Some recent progress on extremal problems in graph theory,
Congr. Numerantium 14 (1975), 3--14.

\bibitem{ErdRome1976-35}
P. Erd\H{o}s:
Problems and results in combinatorial analysis,
Colloquio Internazionale sulle Teorie Combinatorie (Rome, 1973), Tomo II,
Atti dei Convegni Lincei, No. 17, pp. 3--17, Accad. Naz. Lincei, Rome, 1976.

\bibitem{Erd1979-17}
P. Erd\H{o}s:
Problems and results in graph theory and combinatorial analysis,
Graph theory and related topics (Proc. Conf., Univ. Waterloo, Waterloo, Ont., 1977), pp. 153--163,
Academic Press, New York-London.

\bibitem{ErdMost}
P. Erd\H{o}s:
On the combinatorial problems which I would most like to see solved,
Combinatorica 1 (1981), no. 1, 25--42.

\bibitem{ErdosCubeConj}
P. Erd\H{o}s:
On some problems in graph theory, combinatorial analysis and combinatorial number theory,
Graph Theory and Combinatorics (Cambridge, 1983), pp. 1--17, Academic Press, London, 1984.

\bibitem{ErdosTwoProblems}
P. Erd\H{o}s:
Two problems in extremal graph theory.
Graphs Combin. 2 (1986), no. 1, 189--190.

\bibitem{ErdFavourTheor}
P. Erd\H{o}s:
On some of my favourite theorems,
Combinatorics, Paul Erd\H{o}s is eighty, Vol. 2 (Keszthely, 1993), 97--132,
Bolyai Soc. Math. Stud., 2, J\'anos Bolyai Math. Soc., Budapest, 1996.

\bibitem{ErdFauPachSpe}
P. Erd\H{o}s, R. J. Faudree, J. Pach, and J. Spencer:
How to make a graph bipartite,
J. Combin. Theory Ser. B 45 (1988), no. 1, 86--98.

\bibitem{ErdFaudSchelpSim}
P. Erd\H{o}s, R. J. Faudree, R. H. Schelp, and M. Simonovits:
An extremal result for paths,
Graph theory and its applications: East and West (Jinan, 1986), 155--162,
Ann. New York Acad. Sci., 576, New York Acad. Sci., New York, 1989.

\bibitem{EFLS95}
P. Erd\H{o}s, Z. F\"uredi, M. Loebl, and V. T. S\'os:
Discrepancy of trees,
Studia Sci. Math. Hungar. 30 (995), no. 1-2, 47--57.

\bibitem{ErdGallai}
P. Erd\H{o}s and T. Gallai:
On maximal paths and circuits of graphs,
Acta Math. Acad. Sci. Hungar. 10 (1959), 337--356.

\bibitem{ErdGyoriSim} P. Erd\H{o}s, E. Gy\H{o}ri, and M.  Simonovits: How many
  edges should be deleted to make a triangle-free graph bipartite? Sets,
  graphs and numbers (Budapest, 1991), pp. 239--263, Colloq. Math. Soc. J\'anos
  Bolyai, 60, North-Holland, Amsterdam, 1992.

\bibitem{ErdHarcPach}
P. Erd\H{o}s, G. Harcos, and J. Pach:
Popular distances in 3-space,
Discrete Math. 200 (1999), no. 1-3, 95--99.

\bibitem{ErdMoser}
P. Erd\H{o}s and L. Moser:
Problem 11,
Canad. Math. Bull. 2 (1959), 43.

\bibitem{ErdRenyiEvol}
P. Erd\H{o}s and A. R\'enyi:
On the evolution of random graphs,
Magyar Tud. Akad. Mat. Kutat\'o Int. K\"ozl.  5 (1960), 17--61.

\bibitem{ErdRenyiDiam}
P. Erd\H{o}s and A. R\'enyi:
On a problem in the theory of graphs,
Magyar Tud. Akad. Mat. Kutat\'o Int. K\"ozl. 7 (1962), 623--641.

\bibitem{ErdRenyiSos}
P. Erd\H{o}s, A. R\'enyi, and Vera T. S\'os:
On a problem of graph theory,
Stud Sci. Math. Hung. 1 (1966), 215--235.

\bibitem{ErdSachs}
P. Erd\H{o}s and H. Sachs:
Regul\"are Graphen gegebener Taillenweite mit minimaler Knotenzahl (in German),
Wiss. Z. Martin-Luther-Univ. Halle-Wittenberg Math.-Natur. Reihe 12 (1963), 251--257.

\bibitem{ErdSarkSos}
P. Erd\H{o}s, A. S\'ark\"ozy, and V. T. S\'os:
On product representation of powers, I,
European J. Combin. 16 (1995), 567--588.

\bibitem{ErdSimLim}
P. Erd\H{o}s and M. Simonovits:
A limit theorem in graph theory,
Studia Sci. Math. Hungar. 1 (1966), 51--57. 

\bibitem{ErdSimCube}
P. Erd\H{o}s and M. Simonovits:
Some extremal problems in graph theory,
Combinatorial Theory and Its Applications, I. (Proc. Colloq. Balatonf\"ured, 1969),
North Holland, Amsterdam, 1970, pp. 377--390.

\bibitem{ErdSimOcta}
P. Erd\H{o}s and M. Simonovits:
An extremal graph problem,
Acta Math. Acad. Sci. Hungar. 22 (1971/72), 275--282.  

\bibitem{ErdSimCubeWat}
P. Erd\H{o}s and M. Simonovits:
Cube-supersaturated graphs and related problems,
Progress in Graph Theory (Waterloo, Ont., 1982), pp. 203--218,
Academic Press, Toronto, Ont., 1984.

\bibitem{ErdSimComp}
P. Erd\H{o}s and M. Simonovits:
Compactness results in extremal graph theory,
Combinatorica 2 (1982), no. 3, 275--288.

\bibitem{ErdSimSuper}
P. Erd\H{o}s and M. Simonovits:
Supersaturated graphs and hypergraphs,
Combinatorica 3 (1983), 181--192.

\bibitem{ErdStone}
P. Erd\H{o}s and A. M. Stone:
On the structure of linear graphs,
Bull. Amer. Math. Soc 52 (1946), 1087--1091. 

\bibitem{FanGDistri}
G. Fan:
Distribution of cycle lengths in graphs,
J. Combin. Theory Ser. B 84 (2002), 187--202.

\bibitem{FanXWang}
G. Fan, Xuezheng Lv, and Pei Wang:
Cycles in 2-connected graphs,
J. Combin. Theory Ser. B 92 (2004), no. 2, 379--394.

\bibitem{FaudreeSchelpJCT75}
R. J. Faudree and R. H. Schelp:
Path Ramsey numbers in multicolorings,
J. Combin. Theory Ser. B 19 (1975), no. 2,  150--160.

\bibitem{FaudreeSimCCA}
R. J. Faudree and M. Simonovits:
On a class of degenerate extremal graph problems,
Combinatorica 3 (1983), 83--93.

\bibitem{FaudreeSimDegII}
R. J. Faudree and M. Simonovits:
On a class of degenerate extremal problems II, preprint.  

\bibitem{FoxSudDRCSurv}
J. Fox and B. Sudakov:
Dependent random choice,
Random Structures Algorithms 38 (2011), no. 1-2, 68--99.

\bibitem{FirkeWilli}
F. A. Firke, P. M. Kosek, E. D. Nash, and J. Williford:
Extremal graphs without 4-cycles,
http://arxiv.org/pdf/1201.4912v1.pdf

\bibitem{FureC4JCT83}
Z. F\"uredi:
Graphs without quadrilaterals,
J. Combin. Theory Ser. B 34 (1983), 187--190.

\bibitem{FureC4Prep}
Z. F\"uredi:
Quadrilateral-free graphs with maximum number of edges,
preprint 1988,\\
http://www.math.uiuc.edu/\~{}z-furedi/PUBS/furedi\_C4from1988.pdf

\bibitem{FureDiam3}
Z. F\"uredi:
Graphs of diameter 3 with the minimum number of edges,
Graphs Combin. 6 (1990), no. 4, 333--337.

\bibitem{FurConvGeom}
Z. F\"uredi:
The maximum number of unit distances in a convex $n$-gon,
J. Combin. Theory Ser. A 55 (1990), no. 2, 316--320.

\bibitem{FurediL11}
Z. F\"uredi:
On a Tur\'an type problem of Erd\H{o}s,
Combinatorica 11 (1991), 75--79.

\bibitem{FureLondon}
Z. F\"uredi:
Tur\'an type problems,
in Surveys in Combinatorics, London Math. Soc. Lecture Note Ser. 166,
Cambridge University Press, Cambridge, UK, 1991, pp. 253--300.

\bibitem{FureDiamMin2}
Z. F\"uredi:
The maximum number of edges in a minimal graph of diameter 2,
J. Graph Theory 16 (1992), no. 1, 81--98.

\bibitem{FureZurich}
Z. F\"uredi:
Extremal hypergraphs and combinatorial geometry,
Proceedings of the International Congress of Mathematicians, Vol. 1, 2 (Z\"urich, 1994), pp. 1343--1352,
Birkh\"auser, Basel, 1995.

\bibitem{FureC4JCT96}
Z. F\"uredi:
On the number of edges of quadrilateral-free graphs,
J. Combin. Theory Ser. B 68 (1996), 1--6.

\bibitem{FureK33}
Z. F\"uredi:
An upper bound on Zarankiewicz problem,
Combin. Probab. Comput. 5 (1996), no. 1, 29--33.

\bibitem{FureK2t}
Z. F\"uredi:
New asymptotics for bipartite Tur\'an numbers,
J. Combin. Theory Ser. A 75 (1996), no. 1, 141--144.

\bibitem{FurediC5}
Z. F\"uredi:
On the number of fivecycles, 
manuscript, unpublished, superseded by  \cite{GrzesikC5}

\bibitem{FureCube}
Z. F\"uredi:
On a theorem of Erd\H{o}s and Simonovits on graphs not containing the cube,
to appear 

\bibitem{FureHajn}
Z. F\"uredi and Peter Hajnal:
Davenport--Schinzel theory of matrices,
Discrete Math. 103 (1992), 231--251.

\bibitem{FureNaorVerstra}
Z. F\"uredi, A. Naor, and J. Verstra\"ete:
On the Tur\'an number for the hexagon,
Adv. Math. 203 (2006), no. 2, 476--496.

\bibitem{FurOzkahya}
Z. F\"uredi and L. \"Ozkahya:
On even-cycle-free subgraphs of the hypercube,
J. Combin. Theory Ser. A 118 (2011), 1816--1819.

\bibitem{FurePikhSim}
Z. F\"uredi, O. Pikhurko, and M. Simonovits:
The Tur\'an density of the hypergraph $\{ abc,ade,bde,cde\}$,
Electronic J. Combin. 10 (2003), R18.

\bibitem{FureSimFano}
Z. F\"uredi and  M. Simonovits:
Triple systems not containing a Fano configuration,
Combin. Probab. Comput. 14 (2005), no. 4, 467--484.

\bibitem{FureWest}
Z. F\"uredi and D. West:
Ramsey theory and bandwidth of graphs,
Graphs and Combin. 17 (2001), 463--471.

\bibitem{GarnLazebKwong}
D. K. Garnick,  Y. H. H. Kwong, and F. Lazebnik:
Extremal graphs without three-cycles or four-cycles,
J. Graph Theory 17 (1993), no. 5, 633--645.

\bibitem{GarnickNieu}
D. K. Garnick, and N. A. Nieuwejaar,
Non-isomorphic extremal graphs without three-cycles and four-cycles,
J. Combin. Math. Combin. Comput. 12 (1992), 33--56.






\bibitem{GrzesikC5}
A. Grzesik:
On the maximum number of five-cycles in a triangle-free graph,
J. Combin. Theory Ser. B 102 (2012), no. 5, 1061--1066.

\bibitem{GriggsHo}
J. R. Griggs and  Chih-Chang Ho:
On the half-half case of the Zarankiewicz problem,
Discrete Math. 249 (2002), no. 1-3, 95--104.

\bibitem{GriggsQuyang}
J. Griggs, J. Ouyang:
$(0,1)$-matrices with no half-half submatrix of ones,
European J. Combin. 18 (1997), 751--761.

\bibitem{GriggsSim}
J. R. Griggs, M. Simonovits, and George Rubin Thomas:
Extremal graphs with bounded densities of small subgraphs,
J. Graph Theory 29 (1998), no. 3, 185--207.

\bibitem{GuyTihany}
R. K. Guy:
A problem of Zarankiewicz,
in: Theory of Graphs (Proc. Colloq., Tihany, 1966), pp. 119--150.
Academic Press, New York 1968.

\bibitem{GuyZnam}
R. K. Guy and S. Zn\'am:
A problem of Zarankiewicz,
Recent Progress in Combinatorics (Proc. Third  Waterloo Conf. on Combinatorics, 1968), pp. 237--243.
Academic Press, New York 1969.

\bibitem{GyarfasK}
A. Gy\'arf\'as:
Graphs with $k$ odd cycle lengths,
Discrete Math. 103 (1992), 41--48.

\bibitem{GyarfasKomSzem}
A. Gy\'arf\'as, J. Koml\'os, and E. Szemer\'edi:
On the distribution of cycle lengths in graphs,
J. Graph Theory 8 (1984), 441--462.

\bibitem{GyRS}
A. Gy\'arf\'as, C. C. Rousseau, and R. H. Schelp:
An extremal problem for paths in bipartite graphs,
J. Graph Theory 8 (1984), 83--95.

\bibitem{GyoriC3C5}
E. Gy\H{o}ri:
On the number of $C_5$'s in a triangle-free graph,
Combinatorica 9 (1989), 101--102.

\bibitem{GyoriC6}
E. Gy\H{o}ri:
$C_6$-free bipartite graphs and product representation of squares,
Graphs Combin. (Marseille, 1995),
Discrete Math. 165/166 (1997), 371--375.

\bibitem{Gyori2006CPC}
E. Gy\H ori:
Triangle-free hypergraphs,
Combin. Prob. Comput. 15 (2006), 185--191.

\bibitem{GyoriRothRuc}
E. Gy\H{o}ri, B. Rothschild, and A. Ruci\'nski:
Every graph is contained in a sparsest possible balanced graph,
Math. Proc. Cambridge Philos. Soc. 98 (1985), no. 3, 397--401.

\bibitem{HaggkScottArith}
R. H\"aggkvist and A. D. Scott:
Arithmetic progressions of cycles,
Tech. Rep. Mat. Inst. Ume\"a Univ. 16, (1998).

\bibitem{HaggkvScottCyc}
R. H\"aggkvist and A. Scott:
Cycles of nearly equal length in cubic graphs, Preprint.

\bibitem{HartmanRyll}
S. Hartman, J. Mycielski, C. Ryll-Nardzevski:
Syst\`emes sp\'eciaux de points \`a coordonn\'ees enti\'eres,
Colloq. Math.  3 (1954),  84--85,
(Bericht \"Uber di Tagung der Poln Math Gesellschaft, Wroclaw, am 20. September 1951.)  

\bibitem{HatamiSidor}
H. Hatami:
Graph norms and Sidorenko's conjecture,
Israel J. Math. 175 (2010), 125--150.

\bibitem{HataRazboC5C3}
H. Hatami, J. Hladk\'y, D. Kr\'al, S. Norine, and A. Razborov:
On the number of pentagons in triangle-free graphs,
J. Combin. Theory Ser. A 120 (2013), no. 3, 722--732.

\bibitem{HatamiNorine}
H. Hatami and S. Norine:
Undecidability of linear inequalities in graph homomorphism densities,
J. Amer. Math. Soc. 24 (2011), no. 2, 547--565.

\bibitem{HladKomPiguSimSteinSzem}
J. Hladk\'y, J. Koml\'os, D. Piguet, M. Simonovits, M. Stein, and E. Szemer\'edi:
The approximate Loebl-Koml\'os-S\'os Conjecture,
submitted, on arXiv:1211.3050.v1, 2012, Nov 13.

\bibitem{HladkyPiguet} J. Hladk\'y and D. Piguet:
Loebl-Koml\'os-S\'os Conjecture: dense case,
Manuscript (arXiv:0805:4834).  

\bibitem{HuxleyIwaniec}
M. N. Huxley and H. Iwaniec:
Bombieri's theorem in short intervals,
Mathematika 22 (1975), 188--194.

\bibitem{Hylten}
C. Hylt\'en-Cavallius:
On a combinatorial problem,
Colloq. Math. 6 (1958), 59--65.

\bibitem{ImrichGirth}
W. Imrich:
Explicit construction of graphs without small cycles,
Combinatorica 4 (1984), 53--59.

\bibitem{JaggerStovicekThomasonCCA1996}
C. Jagger, P. \v{S}\'tov\'{\i}\v{c}ek, and A. Thomason:
Multiplicities of subgraphs,
Combinatorica 16 (1996), no. 1, 123--141.

\bibitem{JansonLuczRuc}
S. Janson, T. \L{u}czak, and A. Ruci\'nski:
Random Graphs,
Wiley-Interscience Series in Discrete Mathematics and Optimization. Wiley-Interscience,
  New York, 2000. xii+333 pp.

\bibitem{Jiang1}
T. Jiang:
Compact topological minors in graphs,
J. Graph Theory 67 (2011), 139--152.

\bibitem{JiangSeiver}
T. Jiang and R. Seiver:
Tur\'an numbers of subdivided graphs,
SIAM J. Discrete Math. 26 (2012), no. 3, 1238--1255.

\bibitem{JozsaSzem}
S. J\'ozsa and E. Szemer\'edi:
The number of unit distance on the plane,
Infinite and finite sets (Colloq., Keszthely, 1973;
dedicated to P. Erd\H{o}s on his 60th birthday), Vol. II, pp. 939--950.
Colloq. Math. Soc. J\'anos Bolyai, Vol. 10, North-Holland, Amsterdam, 1975.

\bibitem{KatonaTurSurv}
G. O. H. Katona:
Tur\'an's graph theorem and probability theory,
Tur\'an Memorial: Number theory, Analysis and Combinatorics, de Gruyter, Berlin,
to appear.

\bibitem{KatNemSim}
Gy. Katona, T. Nemetz, and M. Simonovits:
On a problem of Tur\'an in the theory of graphs,
Mat. Lapok 15 (1964), 228--238.

\bibitem{KeevashHyperSurv}
P. Keevash:
Hypergraph Turan problems,
Surveys in Combinatorics, Cambridge University Press, 2011, 83--140.

\bibitem{KeevSud}
P. Keevash and B. Sudakov:
The Turan number of the Fano plane,
Combinatorica 25 (2005), 561--574.

\bibitem{KeevSudVersES}
P. Keevash, B. Sudakov, and J. Verstra\"ete:
On a conjecture of Erd\H{o}s and Simonovits: even cycles,
Combinatorica, to appear.

\bibitem{KlazarStW}
M. Klazar:
The F\"uredi--Hajnal conjecture implies the Stanley--Wilf conjecture,
in: D. Krob, A. A. Mikhalev, A. V. Mikhalev (Eds.),
Formal Power Series and Algebraic Combinatorics, Springer, Berlin, 2000, pp. 250--255.

\bibitem{KollRonySzab}
J. Koll\'ar, L. R\'onyai, and T. Szab\'o:
Norm graphs and bipartite Tur\'an numbers,
Combinatorica 16 (1996), 399--406.

\bibitem{KomlosSzem}
J. Koml\'os and E. Szemer\'edi:
Topological cliques in graphs,
Combin. Probab. Comput. 3 (1994), no. 2, 247--256.

\bibitem{KomlosSzemB}
J. Koml\'os and E. Szemer\'edi:
Topological cliques in graphs II,
Combin. Probab. Comput. 5 (1996), 79--90.

\bibitem{Kopylov}
G. N. Kopylov:
Maximal paths and cycles in a graph,
Dokl. Akad. Nauk SSSR 234 (1977), no. 1, 19--21.
(English translation: Soviet Math. Dokl. 18 (1977), no. 3, 593--596.)

\bibitem{KostoPyber}
A. Kostochka and L. Pyber:
Small topological complete subgraphs of "dense'' graphs,
Combinatorica 8 (1988), 83--86.

\bibitem{KovSosTur}
T. K\H{o}v\'ari, V. T. S\'os, and P. Tur\'an:
On a problem of K. Zarankiewicz,
Colloq. Math. 3 (1954), 50--57.

\bibitem{KuehnOstGyori}
D. K\"uhn and D. Osthus:
Four-cycles in graphs without a given even cycle,
J. Graph Theory 48 (2005), 147-156.

\bibitem{LamVer}
T. Lam and J. Verstra\"ete:
A note on graphs without short even cycles,
Electron. J. Combin. 12 (2005), Note 5, 6 pp.

\bibitem{LazebMubayi}
F. Lazebnik and  D. Mubayi:
New lower bounds for Ramsey numbers of graphs and hypergraphs,
Adv. in Appl. Math. 28 (2002), no. 3-4, 544--559.

\bibitem{LazUstiExa}
F. Lazebnik and V. A. Ustimenko,
New examples of graphs without small cycles and of large size,
European J. Combin. 14 (1993), no.5, 445--460.

\bibitem{LazUstWolProperties}
F. Lazebnik, V. A.  Ustimenko, and A. J. Woldar:
Properties of certain families of $2k$-cycle-free graphs,
J. Combin. Theory Ser. B 60 (1994), no. 2, 293--298.

\bibitem{LazebUstimWoldarGirth}
F. Lazebnik, V. A.  Ustimenko, and A. J. Woldar:
A new series of dense graphs of high girth,
Bull. Amer. Math. Soc. 32 (1995), no. 1, 73--79.

\bibitem{LazebUstimWoldPolar}
F. Lazebnik, V. A. Ustimenko, and A. J. Woldar,
Polarities and $2k$-cycle-free graphs,
Discrete Math. 197/198 (1999), 503--513.

\bibitem{LazebWoldarEqu}
F. Lazebnik and A. J. Woldar:
General properties of some families of graphs defined by systems of equations,
J. Graph Theory 38 (2001), no. 2, 65--86.

\bibitem{LenzMubayi}
J. Lenz and D. Mubayi:
Multicolor Ramsey numbers for complete bipartite versus complete graphs,
arXiv 1201.2123, 26 pp.  


\bibitem{LLP:forests}
B. Lidick\'y, Hong Liu, and C. Palmer:
On the Tur\'an number of forests, arXiv 1204.3102.

\bibitem{LovCrit}
L. Lov\'asz:
Independent sets in critical chromatic graphs,
Studia Sci. Math. Hungar.  8 (1973), 165--168.

\bibitem{LovCombExerc}
L. Lov\'asz:
Combinatorial Problems and Exercises, 2nd Ed.,
North-Holland, Amsterdam, 1993.

\bibitem{LovaszLimitBook}
L. Lov\'asz:
Large Networks and Graph Limits,
Colloquium Publications 2012, 475 pp.

\bibitem{LovSimBirk}
L. Lov\'asz and M. Simonovits:
On the number of complete subgraphs of a graph II,
Studies in Pure Mathematics, pp. 459-495, (dedicated to the memory of P. Tur\'an),
Akad\'emiai Kiad\'o and Birkh\"auser Verlag 1982.

\bibitem{RazborFlag}
A. A. Razborov:
Flag algebras,
J. Symbolic Logic 72 (2007), no. 4, 1239--1282.

\bibitem{LuCube}
Linyuan Lu:
Hexagon-free subgraphs in hypercube $Q_n$,
private communication.

\bibitem{LuboPhilSarA}
A. Lubotzky, R. Phillips, and P. Sarnak:
Ramanujan graphs,
Combinatorica 8 (1988), no. 3, 261--277.

\bibitem{LennanTree}
A. 
McLennan:
The Erd\H{o}s-S\'os conjecture for trees of diameter four,
J. Graph Theory 49 (2005), no. 4, 291--301.

\bibitem{Mader}
W. Mader:
Homomorphieeigenschaften und mittlere Kantendichte von Graphen,
Math. Ann. 174 (1967), 265--268.

\bibitem{MaderCCA1998}
W. Mader:
Topological subgraphs in graphs of large girth,
Combinatorica 18 (1998), no. 3, 405--412.

\bibitem{MaderStirin}
W. Mader:
Topological minors in graphs of minimum degree $n$,
Contemporary trends in discrete mathematics (\v{S}ti\v{r}\'{\i}n Castle, 1997), 199--211,
DIMACS Ser. Discrete Math. Theoret. Comput. Sci., 49, Amer. Math. Soc., Providence, RI, 1999.

\bibitem{MaderCCA2005}
W. Mader:
Graphs with $3n-6$ edges not containing a subdivision of $K_5$,
Combinatorica 25 (2005), no. 4, 425--438.

\bibitem{MarcusTardos}
A. Marcus and G. Tardos:
Excluded permutation  matrices and the Stanley-Wilf conjecture,
J. Combin. Theory Ser. A 107 (2004), no. 1, 153--160.

\bibitem{MargulisCCA}
G. A. Margulis:
Explicit construction of graphs without short cycles and low density codes,
Combinatorica 2 (1982), 71--78.

\bibitem{MarguTaskent}
G. A. Margulis:
Arithmetic groups and graphs without short cycles,
in: 6th Int. Symp. on Information Theory,   Tashkent, Abstracts 1, 1984, pp. 123--125 (in Russian).

\bibitem{MargulisExpander}
G. A. Margulis:
Explicit group-theoretical construction of combinatorial schemes and
 their application to the design of expanders and concentrators,
J. Problems of Inform. Trans. 24 (1988), 39-46;
translation from Problemy Peredachi  Informatsii 24 (January-March 1988), 51-60.

\bibitem{MegyeSzaboTac}
G. Megyesi and E. Szab\'o: On the tacnodes of configurations of conics in the projective plane,
Math. Ann. 305 (1996), no. 4, 693--703.

\bibitem{MolloyReedBook}
M. Molloy and B. Reed:
Graph Colouring and the Probabilistic Method,
Algorithms and Combinatorics, 23. Springer-Verlag, Berlin, 2002, xiv+326 pp. 

\bibitem{Montagh}
B. Mont\'agh:
Unavoidable substructures,
PHD Thesis, University of Memphis, May 2005.

\bibitem{Mors}
M. M\"ors:
A new result on the problem of Zarankiewicz,
J. Combin. Theory Ser. A 31 (1981), no. 2, 126--130.


\bibitem{MubayiTuranGy}
D. Mubayi and Gy. Tur\'an:
Finding bipartite subgraphs efficiently,
Inform. Process. Lett. 110 (2010), no. 5, 174--177.

\bibitem{NagyZL}
Z. L. Nagy: 
A multipartite version of the Tur\'an problem --- density conditions and eigenvalues,
Electron. J. Combin. 18 (2011), no. 1, Paper 46, 15 pp.

\bibitem{NikifEigenII}
V. Nikiforov: Bounds on graph eigenvalues II,
Linear Algebra Appl. 427 (2007), 183--189.

\bibitem{NikifZaran}
V. Nikiforov:
A contribution to the Zarankiewicz problem,
Linear Algebra Appl. 432 (2010), no. 6, 1405--1411.

\bibitem{PachAgarwal}
J. Pach and P. K. Agarwal:
Combinatorial Geometry,
Wiley-Interscience, New York, 1995. xiv+354 pp.

\bibitem{PalfySzalay}
P. P. P\'alfy and M. Szalay, in:
Tur\'an Memorial: Number theory, Analysis and Combinatorics, de Gruyter, Berlin,
to appear.

\bibitem{PiguetStein08}
D. Piguet and M. J. Stein:
Loebl-Koml\'os-S\'os conjecture for trees of diameter $5$,
Electron. J. Combin. 15 (2008), Research Paper 106, 11 pp.

\bibitem{PiguetStein12}
D. Piguet and M. J. Stein:
An approximate version of the Loebl-Koml\'os-S\'os conjecture,
J. Combin. Theory Ser. B 102 (2012), no. 1, 102--125.

\bibitem{PikhurkoC2k}
O. Pikhurko:
A note on the Tur\'an function of even cycles,
Proc. Amer. Math Soc. 140 (2012), 3687-3992.

\bibitem{PinchasiSharir}
R. Pinchasi and M. Sharir:
On graphs that do not contain the cube and related problems,
Combinatorica 25 (2005), no. 5, 615--623.

\bibitem{Reiher} C. Reiher: 
The clique density theorem,
arxiv1212.2454.

\bibitem{ReimanZara}
I. Reiman:
\"Uber ein Problem von K. Zarankiewicz,
Acta Math. Acad. Sci. Hungar. 9 (1958), no. 3-4, 269--273.

\bibitem{ReimanLapok}
I. Reiman:
An extremal problem in graph theory, (Hungarian).
Mat. Lapok 12 (1961), 44--53.

\bibitem{RenyiColl}
A. R\'enyi:
Selected Papers of Alfr\'ed R\'enyi, Akad\'emiai Kiad\'o, 1976 (ed. Paul Tur\'an).

\bibitem{RodlSchacht} V. R\"odl and M. Schacht:
Extremal results for random  graphs,
in this volume.

\bibitem{RuzsaSzemer}
I. Z. Ruzsa and E. Szemer\'edi:
Triple systems with no six points carrying three triangles,
Combinatorics (Proc. Fifth Hungarian Colloq., Keszthely, 1976), Vol. II, pp. 939--945,
Colloq. Math. Soc. J\'anos Bolyai, 18, North-Holland, Amsterdam-New York, 1978.




\bibitem{SacWoz}
J.-F. Sacl\`e and M. Wo\'zniak:
A note on the Erd\H{o}s-S\'os conjecture for graphs without $C_4$,
J. Combin. Theory Ser. B 70 (1997), no.2, 367--372.

\bibitem{SarkoGNC6}
G. N. S\'ark\"ozy:
Cycles in bipartite graphs and an application in number theory,
J. Graph Theory 19 (1995), 323--331.

\bibitem{ScottRegu}
A. Scott:
Szemer\'edi's regularity lemma for matrices and sparse graphs,
Combin. Probab. Comput. 20 (2011), no. 3, 455--466.

\bibitem{Shen_K33MinusEdge}
Jian Shen: 
On two Tur\'an numbers,
J. Graph Theory 51 (2006), 244--250.

\bibitem{SidorTree}
A. F. Sidorenko:
Asymptotic solution for a new class of forbidden $r$-graphs,
Combinatorica  9 (1989), no. 2, 207--215.

\bibitem{Sidorenko1993}
A. Sidorenko:
A correlation inequality for bipartite graphs,
Graphs Combin. 9 (1993), no. 2, 201--204.

\bibitem{SidorSurv}
A. F. Sidorenko:
What do we know and what we do not know about Tur\'an Numbers,
Graphs Combin. 11 (1995), no. 2, 179--199.

\bibitem{SimTihany}
M. Simonovits:
A method for solving extremal problems in graph theory,
Theory of Graphs, Proc. Colloq. Tihany, (1966), (P. Erd\H{o}s and G. Katona, Eds.), pp. 279--319,
Acad. Press, New York, 1968.  

\bibitem{SimCrit}
M. Simonovits:
On colour-critical graphs,
Studia Sci. Math. Hungar.  7 (1972), 67--81.

\bibitem{SimColumbus}
M. Simonovits:
Note on a hypergraph extremal problem,
Hypergraph Seminar, Columbus Ohio USA, 1972, (C. Berge and D. K. Ray-Chaudhuri, Eds.),
Lecture Notes in Mathematics 411,  pp. 147--151, Springer Verlag, 1974.

\bibitem{SimSymm77}
M. Simonovits:
Extremal graph problems with symmetrical extremal graphs, additional chromatic conditions,
Discrete Math. 7 (1974), 349--376.

\bibitem{SimTurInflu}
M. Simonovits:
On Paul Tur\'an's influence on graph theory,
J. Graph Theory 1 (1977), no. 2, 102--116.

\bibitem{SimBirk}
M. Simonovits:
Extremal graph problems and graph products,
Studies in Pure Mathematics, pp. 669--680, (dedicated to the memory of P. Tur\'an),
Akad\'emiai Kiad\'o and Birkh\"auser Verlag 1982.

\bibitem{SimFra}
M. Simonovits:
Extremal graph theory,
in: L. W. Beineke, R. J. Wilson (Eds.), Selected Topics in Graph Theory II., pp. 161--200,
Academic Press, London, 1983.

\bibitem{SimWat}
M. Simonovits:
Extremal graph problems, degenerate extremal problems and supersaturated graphs,
Progress in graph Theory, (Bondy and Murty, Eds.), pp. 419--438, Academic Press, 1984.

\bibitem{SimStirin}
M. Simonovits:
How to solve a Tur\'an type extremal graph problem? (linear decomposition),
Contemporary trends in discrete mathematics (Stirin Castle, 1997), pp. 283--305,
DIMACS Ser. Discrete Math. Theoret. Comput. Sci., 49, Amer. Math. Soc., Providence, RI, 1999.

\bibitem{SimErdosInflu}
M. Simonovits:
Paul Erd\H{o}s' influence on extremal graph theory,
The mathematics of Paul Erd\H{o}s, II., pp. 148--192, Algorithms Combin., 14, Springer, Berlin, 1997.

\bibitem{SimPragNew}
M. Simonovits:
Paul Erd\H{o}s' influence on extremal graph theory,
The new version of the old paper  \cite{SimErdosInflu}.  

\bibitem{SimTur2013}
M. Simonovits:
Paul Tur\'an's influence in Combinatorics,
in Tur\'an Memorial: Number Theory, Analysis, and Combinations, De Gruyter, 
to appear.

\bibitem{SimSosRT}
M. Simonovits and V. T. S\'os:
Ramsey-Tur\'an theory, Combinatorics, graph theory, algorithms and applications,
Discrete Math. 229 (2001), no. 1-3, 293--340.

\bibitem{Singleton}
R. R. Singleton:
On minimal graphs of maximum even girth,
J. Combinatorial Theory  1 (1966), 306--332.  

\bibitem{SosLincei}
V. T. S\'os:
Remarks on the connection of graph theory, finite geometry and block designs,
 Colloquio Internazionale sulle Teorie Combinatorie (Roma, 1973),
 Tomo II, pp. 223--233, Atti dei Convegni Lincei, No. 17, Accad. Naz. Lincei, Rome, 1976.

\bibitem{SpencSzemTrot}
J. Spencer, E. Szemer\'edi, and W. T. Trotter:
Unit distances in the Euclidean plane,
Graph theory and combinatorics (Cambridge, 1983), pp. 293--303,
Academic Press, London, 1984.

\bibitem{SudakVerstra}
B. Sudakov and J. Verstra\"ete:
Cycle lengths in sparse graphs,
Combinatorica 28 (2008), no. 3, 357--372.

\bibitem{SzemRegu}
E. Szemer\'edi:
Regular partitions of graphs,
Problemes Combinatoires et Theorie des Graphes (ed. I.-C. Bermond et al.), pp. 399--401,
CNRS, Paris, 1978. 

\bibitem{Tardos2005}
G. Tardos:
On $0$-$1$ matrices and small excluded submatrices,
J. Combin. Th. Ser. A 111 (2005), 266--288.

\bibitem{ThomasonDisproof}
A. G. Thomason:
A disproof of a conjecture of Erd\H{o}s in Ramsey Theory,
J. London Math. Soc. 39 (1989), 246-255.

\bibitem{ThomasonW}
A. Thomason and P. Wagner:
Bounding the size of square-free subgraphs of the hypercube,
Discrete Math. 309 (2009), 1730--1735.

\bibitem{Timmons}
C. 
 M. Timmons:
Ordered Tur\'an Problems,
Lecture no. 1086-05-1067 on the Joint Mathematics Meetings, San Diego, CA, January 9, 2013.

\bibitem{ToftCrit}
B. Toft:
Two theorems on critical $4$-chromatic graphs,
Studia Sci. Math. Hungar. 7 (1972), 83--89.

\bibitem{TuranHardyRam}
P. Tur\'an:
On a theorem of Hardy-Ramanujan,
Journal of London Math Soc.  9 (1934), 274--276.

\bibitem{TuranML}
P. Tur\'an:
On an extremal problem in graph theory, (Hungarian),
Mat. Fiz. Lapok 48 (1941), 436--452.

\bibitem{TuranColloq}
P. Tur\'an:
On the theory of graphs,
Colloq. Math. 3 (1954), 19--30.

\bibitem{TurWel}
P. Tur\'an:
A note of welcome,
J. Graph Theory 1 (1977), 7--9.

\bibitem{ValtrUnit} P. Valtr: Strictly convex norms allowing many unit
  distances and related touching questions, manuscript. 

\bibitem{VertreAPCyc}
J. Verstra\"ete:
On arithmetic progressions of cycle lengths in graphs,
Combin. Probab. Comput. 9 (2000), no.4, 369--373.

\bibitem{Wenger}
R. Wenger:
Extremal graphs with no $C_4$'s, $C_6$'s, or $C_{10}$'s,
J. Combin. Theory Ser. B 52 (1991), no. 1, 113--116.

\bibitem{WilsonExistence}
R. M. Wilson:
An existence theory for pairwise balanced designs, III. Proof of the existence conjectures,
J. Combin. Theory Ser. A 18 (1975), 71--79.

\bibitem{WoodallA}
D. R. Woodall:
Maximal circuits of graphs I,
Acta Math. Acad. Sci. Hungar. 28 (1976), no. 1--2, 77--80.

\bibitem{WoodallB}
D. R. Woodall:
Maximal circuits of graphs II,
Studia Sci. Math. Hungar. 10 (1975), no. 1--2, 103--109.

\bibitem{Wozn96}
M. Wo\'zniak:
On the Erd\H{o}s--S\'os conjecture,
J. Graph Theory, 21 (1996), no. 2, 229--234.

\bibitem{YuangRowliC4}
Y. Yuansheng and P. Rowlinson:
On extremal graphs without four-cycles,
Utilitas Math. 41 (1992), 204--210.

\bibitem{YuangRowliC6}
Y. Yuansheng and P. Rowlinson:   
On graphs without 6-cycles and related Ramsey numbers,
Utilitas Math. 44 (1993), 192--196.

\bibitem{Zarank}
K. Zarankiewicz:
Problem 101,
Colloquium Mathematicum 2 (1951), p. 301.

\bibitem{YiZhao11}
Yi Zhao:
Proof of the $(n/2$ --$n/2$ --$n/2)$ conjecture for large $n$,
Electron. J. Combin.  18 (2011), Paper 27.

\bibitem{ZnamB}
\v{S}. Zn\'am:
On a combinatorical problem of  K. Zarankiewicz,
Colloq. Math. 11 (1963), 81--84.

\bibitem{ZnamA}
\v{S}. Zn\'am:
Two improvements of a result concerning a problem of K. Zarankiewicz,
Colloq. Math. 13 (1964/1965), 255--258.


} 

\end{thebibliography}
\end{document}